\documentclass{amsart}
\usepackage{amsthm}

\usepackage[top=1in, left=1in, right=1in, bottom=1in]{geometry}
\usepackage{
amssymb,amsmath,enumerate,enumitem,stackrel,tabu,mathrsfs,verbatim}
\usepackage{xcolor} 
\usepackage{xr-hyper}
\usepackage{nicefrac}

\usepackage{url}

\usepackage{hyperref}
\usepackage[all,cmtip]{xy}
\usepackage[utf8]{inputenc} 
\chardef\cprime"7E 
\usepackage{cases}

\usepackage{theoremref}

\definecolor{labelkey}{rgb}{1,0,0}

\setlist[enumerate,1]{label=\roman*\textup{)},
style=nextline,
topsep=0mm,
itemsep=0.4mm, 
}

\makeatletter
\@addtoreset{equation}{section}

\numberwithin{equation}{section}

\makeatletter
\theoremstyle{definition}
\newtheorem{Defi}[equation]{Definition} \newcommand{\defi}{\begin{Defi}} \newcommand{\xdefi}{\end{Defi}} \newtheorem{DefiLemm}[equation]{Definition and Lemma} \newcommand{\defilemm}{\begin{DefiLemm}} \newcommand{\xdefilemm}{\end{DefiLemm}}
\newtheorem{DefiTheo}[equation]{Definition and Theorem} \newcommand{\defitheo}{\begin{DefiTheo}} \newcommand{\xdefitheo}{\end{DefiTheo}}
\newtheorem{Bsp}[equation]{Example} \newcommand{\exam}{\begin{Bsp}} \newcommand{\xexam}{\end{Bsp}}
\newtheorem{Syno}[equation]{Synopsis} \newcommand{\syno}{\begin{Syno}} \newcommand{\xsyno}{\end{Syno}}
\newtheorem{Ques}[equation]{Question} \newcommand{\ques}{\begin{Ques}} \newcommand{\xques}{\end{Ques}}

\newtheorem{Bem}[equation]{Remark} \newcommand{\rema}{\begin{Bem}} \newcommand{\xrema}{\end{Bem}}
\newtheorem{Nota}[equation]{Notation} \newcommand{\nota}{\begin{Nota}} \newcommand{\xnota}{\end{Nota}}

\theoremstyle{plain}
\newtheorem{Theo}[equation]{Theorem} \newcommand{\theo}{\begin{Theo}} \newcommand{\xtheo}{\end{Theo}}
\newtheorem{Satz}[equation]{Proposition} \newcommand{\prop}{\begin{Satz}} \newcommand{\xprop}{\end{Satz}}
\newtheorem{Lemm}[equation]{Lemma} \newcommand{\lemm}{\begin{Lemm}}
\newcommand{\xlemm}{\end{Lemm}}
\newtheorem{Coro}[equation]{Corollary} \newcommand{\coro}{\begin{Coro}} \newcommand{\xcoro}{\end{Coro}}

\makeatother

\newcommand{\refsect}[1]{§\ref{sect--#1}}
\newcommand{\refit}[1]{\ref{item--#1}}
\newcommand{\refeq}[1]{(\ref{eqn--#1})}

\newcommand{\eqn}{\begin{equation}} \newcommand{\xeqn}{\end{equation}}
\newcommand{\eqnarr}{\begin{eqnarray*}} \newcommand{\xeqnarr}{\end{eqnarray*}}
\newcommand{\eqnarra}{\begin{eqnarray}} \newcommand{\xeqnarra}{\end{eqnarray}}

\newcommand{\pf}{\begin{proof}} \newcommand{\xpf}{\end{proof}}

\newcommand{\nc}{\newcommand}
\nc{\StP}[1]{\cite[\href{http://stacks.math.columbia.edu/tag/#1}{Tag #1}]{StacksProject}} 

\nc{\on}{\operatorname}
\nc{\aff}{{\on{aff}}}
\nc{\modi}{{\on{mod}}} 
\nc{\even}{{\on{even}}}
\nc{\odd}{{\on{odd}}}
\nc{\naive}{{\on{naive}}}
\nc{\hofib}{\on{hofib}}

\nc{\str}{\on{-}}
\nc{\lan}{\langle}
\nc{\ran}{\rangle}

\nc{\bbA}{\A} 
\nc{\bbB}{{\mathbb B}}
\nc{\bbC}{{\bf C}}
\nc{\bbD}{{\mathbb D}}
\nc{\bbE}{{\mathbb E}}
\nc{\bbF}{{\bf F}}
\nc{\bbG}{{\bf G}}
\nc{\bbH}{{\mathbb H}}
\nc{\bbI}{{\mathbb I}}
\nc{\bbJ}{{\mathbb J}}
\nc{\bbK}{{\mathbb K}}
\nc{\bbL}{{\mathbb L}}
\nc{\bbM}{{\mathbb M}}
\nc{\bbN}{{\mathbb N}}
\nc{\bbO}{{\mathbb O}}
\nc{\bbP}{\P} 
\nc{\bbQ}{{\bf Q}}
\nc{\bbR}{{\mathbb R}}
\nc{\bbS}{{\mathbb S}}
\nc{\bbT}{{\mathbb T}}
\nc{\bbU}{{\mathbb U}}
\nc{\bbV}{{\mathbb V}}
\nc{\bbW}{{\mathbb W}}
\nc{\bbX}{{\mathbb X}}
\nc{\bbY}{{\mathbb Y}}
\nc{\bbZ}{{\bf Z}}

\nc{\calA}{{\mathcal A}}
\nc{\calB}{{\mathcal B}}
\nc{\calC}{{\mathcal C}}
\nc{\calD}{{\mathcal D}}
\nc{\calE}{{\mathcal E}}
\nc{\calF}{{\mathcal F}}
\nc{\calG}{{\mathcal G}}
\nc{\calH}{{\mathcal H}}
\nc{\calI}{{\mathcal I}}
\nc{\calJ}{{\mathcal J}}
\nc{\calK}{{\mathcal K}}
\nc{\calL}{{\mathcal L}}
\nc{\calM}{{\mathcal M}}
\nc{\calN}{{\mathcal N}}
\nc{\calO}{{\mathcal O}}
\nc{\calP}{{\mathcal P}}
\nc{\calQ}{{\mathcal Q}}
\nc{\calR}{{\mathcal R}}
\nc{\calS}{{\mathcal S}}
\nc{\calT}{{\mathcal T}}
\nc{\calU}{{\mathcal U}}
\nc{\calV}{{\mathcal V}}
\nc{\calW}{{\mathcal W}}
\nc{\calX}{{\mathcal X}}
\nc{\calY}{{\mathcal Y}}
\nc{\calZ}{{\mathcal Z}}

\nc{\Conv}{{\on{Conv}}}
\nc{\triv}{{\on{triv}}}

\nc{\scrG}{{\mathscr{G}}}
\nc{\scrU}{{\mathscr{U}}}
\nc{\scrB}{{\mathscr{B}}}
\nc{\scrA}{{\mathscr{A}}}
\nc{\bbf}{{\mathbf{f}}}
\nc{\bba}{{\mathbf{a}}}
\nc{\Hecke}{{\on{Hecke}}}
\nc{\Bun}{{\on{Bun}}}
\nc{\inv}{{\on{inv}}} 
\nc{\free}{{\on{free}}} 
\nc{\Rel}{{\on{Pos}}}
\nc{\be}{{\beta}}
\nc{\ad}{{\on{ad}}}
\nc{\bsl}{{\backslash}}

\nc{\Bmu}{{\boldsymbol \mu}}

\nc{\al}{\alpha}
\nc{\la}{\lambda}

\nc{\pot}[1]{ [\hspace{-0,5mm}[ {#1} ]\hspace{-0,5mm}] }
\nc{\rpot}[1]{ (\hspace{-0,7mm}( {#1} )\hspace{-0,7mm}) }

\nc{\defined}{\hspace{0.1cm}\stackrel{\text{\tiny \rm def}}{=}\hspace{0.1cm}}
\nc{\co}{\colon}

\newcommand{\category}[1]{\mathrm{#1}}
\newcommand{\Fin}{\category{Fin}} 
\newcommand{\DGCat}{\category{DGCat}} 
\newcommand{\WCat}{\category{WCat}} 
\newcommand{\cont}{\category{cont}} 
\newcommand{\Cat}{\category{Cat}} 
\newcommand{\CatL}{\Cat^{\category L}} 
\newcommand{\CatR}{\Cat^{\category R}} 
\newcommand{\Ho}{\category{Ho}} 
\newcommand{\open}{\mathrm{open}} 
\newcommand{\iso}{\mathrm{iso}} 
\newcommand{\adm}{\mathrm{adm}} 
\newcommand{\horiz}{\mathrm{horiz}} 
\renewcommand{\vert}{\mathrm{vert}} 
\newcommand{\all}{\mathrm{all}} 
\newcommand{\fd}{\mathrm{fd}} 
\newcommand{\PSh}{\category{PSh}} 
\newcommand{\AffSch}{\category{AffSch}} 
\newcommand{\Groups}{\category{Groups}} 
\newcommand{\PreStk}{\category{PreStk}} 
\newcommand{\Pro}{\category{Pro}} 
\newcommand{\Gpd}{\category{Gpd}} 
\newcommand{\codim}{\mathrm{codim}} 
\newcommand{\Sm}{\category{Sm}} 
\newcommand{\Sch}{\category{Sch}} 
\newcommand{\Sp}{\category{Sp}} 
\newcommand{\Corr}{\category{Corr}} 
\newcommand{\IndSch}{\category{IndSch}} 
\newcommand{\Mod}{\category{Mod}} 
\newcommand{\Alg}{\category{Alg}} 
\newcommand{\BarC}{{\rm Bar}} 
\newcommand{\pt}{{\rm pt}} 
\newcommand{\W}{\mathrm {W}} 
\newcommand{\pl}{\mathrm {pl}} 
\newcommand{\Mot}{\mathrm {Mot}} 

\def\wt{\mathrm {w}} 
\def\perv{\mathrm {perv}} 
\def\gr{\mathrm {gr}} 
\def\pure{\mathrm {pure}} 
\def\pH{^{\mathrm p} \H} 
\def\cl{\mathrm {cl}} 
\def\mot{\mathrm {m}} 
\def\motH{{}^\mot \H} 
\def\Gm{\mathbf {G}_\mathrm m} 
\def\CT{\mathrm {CT}} 
\def\GL{\mathrm {GL}} 
\newcommand{\GmX}[  1]{\mathbf {G}_{\mathrm {m}, #1}} 
\def\IC{\mathrm{IC}} 
\def\red{\mathrm{red}} 
\def\ft{\mathrm{ft}} 
\def\BD{\mathrm{BD}} 

\font\tencyr=wncyr10
\font\sevencyr=wncyr7
\font\fivecyr=wncyr5
\newcommand{\dual}{\vee} 

\newcommand{\colim}{\operatornamewithlimits{colim}} 
\def\id{{\rm id}} 
\def\pr{{\rm pr}} 
\def\opp{{\rm op}} 
\def\To#1#2{\mathop{\count0=#1 \loop\ifnum\count0>0 \smash-\mkern-7mu \advance\count0 -1 \repeat \mathord\rightarrow}\limits^{#2}} 
\def\CH{\mathop{\rm CH}\nolimits} 
\def\Tw{\mathop{\rm Tw}\nolimits} 
\def\Hom{\mathop{\rm Hom}\nolimits} 
\def\End{\mathop{\rm End}\nolimits} 
\def\Ind{\mathop{\rm Ind}\nolimits} 
\def\sep{\mathop{\rm sep}\nolimits} 
\def\all{\mathop{\rm all}\nolimits} 
\def\proper{{\mathop{\rm proper}\nolimits}} 
\def\Gr{\mathop{\rm Gr}\nolimits} 
\def\Fl{\mathop{\rm Fl}\nolimits} 
\def\CAlg{\mathop{\rm CAlg}\nolimits} 
\def\Ext{\mathop{\rm Ext}\nolimits} 
\def\Aut{\mathop{\rm Aut}\nolimits} 
\def\Lie{\mathop{\rm Lie}\nolimits} 
\def\Map{\mathop{\rm Map}\nolimits} 
\def\Rep{\category{Rep}} 
\def\bound{{\rm b}} 
\def\bd{\bound} 
\def\et{\mathrm{et}} 
\def\mix{\mathrm{mix}} 

\definecolor{hellgrau}{RGB}{200,200,200} 
\definecolor{dunkelgrau}{RGB}{160,160,160} 

\def\Z{{\bf Z}} 
\def\Fp{{\bf F}_p} %
\def\Fq{{\bf F}_q} %
\def\Q{{\bf Q}} 
\def\Ql{{\bf Q_\ell}} 
\def\A{{\bf A}} 
\renewcommand{\P}{\mathbf P} 
\def\Gm{\mathbf {G}_\mathrm m} 

\def\H{{\rm H}} 
\def\im{{\rm im}} 
\def\SH{\category{SH}} %
\def\Fun{\category{Fun}} 
\def\DM{\category{DM}} 
\def\DTM{\category{DTM}} 
\def\Vect{\category{Vect}} 
\def\Perv{\category{Perv}} 
\def\MTM{\category{MTM}} 
\def\ii{$\infty$}
\def\bound{{\rm b}} 
\def\Gal{{\rm Gal}} 
\def\Sat{{\rm Sat}} 
\def\PrL{\category{Pr}^{\rm L}} 
\def\PrLt{\category{Pr}^{\rm L, \t}} 
\def\stb{{\rm stb}}

\def\N{{\rm N}} 

\def\Spec{\mathop{\rm Spec}} 
\newcommand{\M}{\mathrm{M}} 
\newcommand{\comp}{\mathrm{c}} 
\newcommand{\cstr}{\mathrm{ct}} 
\newcommand{\C}{\mathcal{C}} 
\newcommand{\D}{\category{D}} 

\def\R{{\rm R}} 
\def\L{{\rm L}} 
\def\sm{{\rm sm}} 
\def\sbuildrel#1\over#2{\mathrel{\smash{\mathop{\kern0pt #2}\limits^{#1}}}}

\let\x\times
\let\ol\overline
\renewcommand{\t}{\otimes}
\newcommand{\ttw}{\widetilde \boxtimes}
\newcommand{\xtw}{\widetilde \x}
\renewcommand{\r}{\rightarrow}
\newcommand{\lr}{\longrightarrow}

\def\matrix#1{\null\,\vcenter{\normalbaselines
    \ialign{\hfil$##$\hfil&&\quad\hfil$##$\hfil\crcr
      \mathstrut\crcr\noalign{\kern-\baselineskip}
      #1\crcr\mathstrut\crcr\noalign{\kern-\baselineskip}}}\,}

\newdimen\harrowsize
\def\mapright#1{\smash{\mathop{\hbox to\harrowsize{\rightarrowfill}}\limits^{#1}}}

{\catcode`@=11
\gdef\cal{\fam\tw@}
\global\let\over\@@over
\global\let\atop\@@atop
\global\let\above\@@above
\global\let\overwithdelims\@@overwithdelims
\global\let\atopwithdelims\@@atopwithdelims
\global\let\abovewithdelims\@@abovewithdelims
\gdef\eqalign#1{\null\,\vcenter{\openup\jot\m@th
  \ialign{\strut\hfil$\displaystyle{##}$&$\displaystyle{{}##}$\hfil
      \crcr#1\crcr}}\,}
\newskip\xcentering \global\xcentering=0pt plus 1000pt minus 1000pt
\gdef\eqalignno#1{\displ@y \tabskip\xcentering
  \halign to\displaywidth{\hfil$\@lign\displaystyle{##}$\tabskip\z@skip
    &$\@lign\displaystyle{{}##}$\hfil\tabskip\xcentering
    &\llap{$\@lign##$}\tabskip\z@skip\crcr
    #1\crcr}}
\gdef\eqlabel#1{\refstepcounter{equation}\label{eqn--#1}\eqno\hbox{\@eqnnum}}
}

\def \nts#1{}

\usepackage{xcolor}
\hypersetup{
    colorlinks,
    linkcolor={red!50!black},
    citecolor={blue!50!black},
    urlcolor={blue!80!black}
}

\usepackage{ifthen}

\def\iref#1{\ifthenelse{\equal{#1}{Bar_Construction}}{3.2.17}{}\ifthenelse{\equal{#1}{Chevalley_Triple}}{4.1.1}{}\ifthenelse{\equal{#1}{coro--DTM.G.X.generators}}{3.2.24}{}\ifthenelse{\equal{#1}{coro--equivalence.torsors}}{2.2.24}{}\ifthenelse{\equal{#1}{coro--equivariant.IC}}{5.3.6}{}\ifthenelse{\equal{#1}{coro--etale.torsor}}{A.4.8}{}\ifthenelse{\equal{#1}{coro--exactness}}{3.2.10}{}\ifthenelse{\equal{#1}{coro--Fl.stratification}}{4.3.12}{}\ifthenelse{\equal{#1}{coro--Ind.Artin.co.limit}}{2.3.4}{}\ifthenelse{\equal{#1}{coro--MTM.G.trivial}}{3.2.21}{}\ifthenelse{\equal{#1}{coro--prelim.intersection.motives}}{6.3.5}{}\ifthenelse{\equal{#1}{coro--t-structure}}{3.2.6}{}\ifthenelse{\equal{#1}{decomp_element}}{4.2.14}{}\ifthenelse{\equal{#1}{defi--adm.act}}{A.3.1}{}\ifthenelse{\equal{#1}{defi--cellular}}{3.1.5}{}\ifthenelse{\equal{#1}{defi--DM.G}}{2.2.6}{}\ifthenelse{\equal{#1}{defi--DM.prestacks}}{2.2.1}{}\ifthenelse{\equal{#1}{defi--DTM.G}}{3.1.21}{}\ifthenelse{\equal{#1}{defi--flag.variety}}{4.3.1}{}\ifthenelse{\equal{#1}{defi--intersection.motive}}{6.3.4}{}\ifthenelse{\equal{#1}{defi--motive.ind.scheme}}{2.3.12}{}\ifthenelse{\equal{#1}{defi--MTM.G}}{3.2.14}{}\ifthenelse{\equal{#1}{defi--rel.pos}}{6.1.2}{}\ifthenelse{\equal{#1}{defi--Schubert.scheme}}{4.4.1}{}\ifthenelse{\equal{#1}{defi--Schubert}}{4.3.4}{}\ifthenelse{\equal{#1}{defi--stratified.dfn}}{3.1.1}{}\ifthenelse{\equal{#1}{defi--stratified.G.action}}{3.1.26}{}\ifthenelse{\equal{#1}{defi--Tate.geometry}}{3.1.8}{}\ifthenelse{\equal{#1}{defi--unipotent}}{A.4.5}{}\ifthenelse{\equal{#1}{defi--Whitney-Tate.map}}{3.1.15}{}\ifthenelse{\equal{#1}{defilemm--Whitney.Tate.condition}}{3.1.11}{}\ifthenelse{\equal{#1}{deformation}}{A.4.10}{}\ifthenelse{\equal{#1}{dom_weights}}{4.1.2}{}\ifthenelse{\equal{#1}{double_classes}}{4.2.15}{}\ifthenelse{\equal{#1}{double_quotient}}{5.3.1}{}\ifthenelse{\equal{#1}{DTM.flag}}{5.2.1}{}\ifthenelse{\equal{#1}{eqn--adjunction.shriek}}{2.1.3}{}\ifthenelse{\equal{#1}{eqn--adjunction.star}}{2.1.2}{}\ifthenelse{\equal{#1}{eqn--Bar.descent}}{2.2.8}{}\ifthenelse{\equal{#1}{eqn--base.change.1}}{2.1.6}{}\ifthenelse{\equal{#1}{eqn--base.change.2}}{2.1.7}{}\ifthenelse{\equal{#1}{eqn--Beck}}{2.1.12}{}\ifthenelse{\equal{#1}{eqn--BS.vanishing}}{3.2.2}{}\ifthenelse{\equal{#1}{eqn--claim1}}{2.2.18}{}\ifthenelse{\equal{#1}{eqn--colim.DGCat.general}}{2.2.4}{}\ifthenelse{\equal{#1}{eqn--colim.DGCat}}{2.2.3}{}\ifthenelse{\equal{#1}{eqn--Det.calX}}{2.3.10}{}\ifthenelse{\equal{#1}{eqn--DM.calX}}{2.3.9}{}\ifthenelse{\equal{#1}{eqn--DM.Cech}}{2.2.20}{}\ifthenelse{\equal{#1}{eqn--DM.G.colim}}{2.3.5}{}\ifthenelse{\equal{#1}{eqn--DM.K-theory}}{2.1.9}{}\ifthenelse{\equal{#1}{eqn--DM.tau}}{2.2.17}{}\ifthenelse{\equal{#1}{eqn--DTM.asymmetric}}{5.3.3}{}\ifthenelse{\equal{#1}{eqn--Hom.X.Y}}{3.2.13}{}\ifthenelse{\equal{#1}{eqn--iota.Fl}}{5.0.2}{}\ifthenelse{\equal{#1}{eqn--iota.v.w}}{5.1.2}{}\ifthenelse{\equal{#1}{eqn--localization}}{2.1.5}{}\ifthenelse{\equal{#1}{eqn--localizationOne}}{2.1.4}{}\ifthenelse{\equal{#1}{eqn--Loop_Grp_Dfn}}{4.2.1}{}\ifthenelse{\equal{#1}{eqn--motivic.t.structure}}{3.2.5}{}\ifthenelse{\equal{#1}{eqn--presentation.Ind.Artin}}{2.3.1}{}\ifthenelse{\equal{#1}{eqn--presentation.Ind.scheme}}{2.4.1}{}\ifthenelse{\equal{#1}{eqn--relative.purity}}{2.1.8}{}\ifthenelse{\equal{#1}{eqn--rho.f!}}{2.3.8}{}\ifthenelse{\equal{#1}{eqn--VE}}{A.4.4}{}\ifthenelse{\equal{#1}{eqn--X.S.GX.etc}}{3.1.24}{}\ifthenelse{\equal{#1}{exam--affine.proj}}{2.2.13}{}\ifthenelse{\equal{#1}{exam--basic.WT}}{3.1.17}{}\ifthenelse{\equal{#1}{exam--BS.exam}}{3.2.3}{}\ifthenelse{\equal{#1}{exam--descent.Tate.not}}{3.1.25}{}\ifthenelse{\equal{#1}{exam--double.quot}}{4.2.12}{}\ifthenelse{\equal{#1}{exam--fusion.loop}}{6.1.3}{}\ifthenelse{\equal{#1}{exam--G.Whitney.Tate}}{3.1.13}{}\ifthenelse{\equal{#1}{exam--groups}}{A.4.12}{}\ifthenelse{\equal{#1}{exam--length.function.exam}}{4.2.16}{}\ifthenelse{\equal{#1}{exam--monoidal.unit}}{2.4.3}{}\ifthenelse{\equal{#1}{exam--motive.grass}}{2.3.13}{}\ifthenelse{\equal{#1}{exam--parahoric}}{4.2.2}{}\ifthenelse{\equal{#1}{exam--simple.reflection}}{4.3.14}{}\ifthenelse{\equal{#1}{exam--stratified.dfn}}{3.1.6}{}\ifthenelse{\equal{#1}{fiber_over_Grass}}{6.1.4}{}\ifthenelse{\equal{#1}{flag_act}}{4.3.2}{}\ifthenelse{\equal{#1}{foot.separated}}{1}{}\ifthenelse{\equal{#1}{funda.diag}}{6.2.3}{}\ifthenelse{\equal{#1}{glob_loop_group}}{6.1.1}{}\ifthenelse{\equal{#1}{intersection.complex}}{3.3.6}{}\ifthenelse{\equal{#1}{invariant.map}}{6.2.1}{}\ifthenelse{\equal{#1}{IW_Indentify}}{4.2.9}{}\ifthenelse{\equal{#1}{IW_Sub_Indentify}}{4.2.10}{}\ifthenelse{\equal{#1}{lemm--adm.finite.type}}{A.3.5}{}\ifthenelse{\equal{#1}{lemm--adm.reduced}}{A.3.3}{}\ifthenelse{\equal{#1}{lemm--adm.strata}}{A.3.2}{}\ifthenelse{\equal{#1}{lemm--affine.proj}}{4.2.7}{}\ifthenelse{\equal{#1}{lemm--base.change.schubert.field}}{4.3.6}{}\ifthenelse{\equal{#1}{lemm--compact.motives.Ind.Artin}}{2.3.6}{}\ifthenelse{\equal{#1}{lemm--descent.Beck.Chevalley}}{2.1.11}{}\ifthenelse{\equal{#1}{lemm--DM.double.tau}}{5.3.2}{}\ifthenelse{\equal{#1}{lemm--DM.G.BarC}}{2.2.7}{}\ifthenelse{\equal{#1}{lemm--DM.G/H}}{2.2.21}{}\ifthenelse{\equal{#1}{lemm--double.orbit}}{4.2.11}{}\ifthenelse{\equal{#1}{lemm--functoriality.equivariant}}{2.2.9}{}\ifthenelse{\equal{#1}{lemm--ind.pres}}{4.3.3}{}\ifthenelse{\equal{#1}{lemm--intermediate.simple}}{3.3.4}{}\ifthenelse{\equal{#1}{lemm--lim.equivalence}}{2.2.12}{}\ifthenelse{\equal{#1}{lemm--loc.triv}}{2.2.23}{}\ifthenelse{\equal{#1}{lemm--Lurie.co.limit}}{2.3.2}{}\ifthenelse{\equal{#1}{lemm--middle.extension}}{3.3.3}{}\ifthenelse{\equal{#1}{lemm--monadic.Ind}}{2.1.16}{}\ifthenelse{\equal{#1}{lemm--orbit.flag}}{4.3.7}{}\ifthenelse{\equal{#1}{lemm--parahoric.defi}}{4.2.4}{}\ifthenelse{\equal{#1}{lemm--pi*.fullyfaithful.MTM}}{3.2.12}{}\ifthenelse{\equal{#1}{lemm--pro.group}}{A.2.1}{}\ifthenelse{\equal{#1}{lemm--smooth.detects.t-structure}}{3.2.11}{}\ifthenelse{\equal{#1}{lemm--smooth.detects.Tate}}{3.1.20}{}\ifthenelse{\equal{#1}{lemm--t-structure.stratum}}{3.2.4}{}\ifthenelse{\equal{#1}{lemm--t.structure.limit}}{3.2.18}{}\ifthenelse{\equal{#1}{lemm--Tate.conservative}}{3.2.8}{}\ifthenelse{\equal{#1}{lemm--Tate.proper.descent}}{3.1.19}{}\ifthenelse{\equal{#1}{lemm--Tate.up.down}}{3.1.18}{}\ifthenelse{\equal{#1}{lemm--torsor.sequence}}{A.4.3}{}\ifthenelse{\equal{#1}{lemm--torsors}}{2.2.22}{}\ifthenelse{\equal{#1}{lemm--var.action}}{6.1.5}{}\ifthenelse{\equal{#1}{length_Grass_Iwahori}}{4.2.18}{}\ifthenelse{\equal{#1}{length_Grass}}{4.2.17}{}\ifthenelse{\equal{#1}{map.intro}}{1.2.1}{}\ifthenelse{\equal{#1}{middle-extension}}{3.3.2}{}\ifthenelse{\equal{#1}{nota--BS.vanishing.ladic}}{5.0.1}{}\ifthenelse{\equal{#1}{nota--BS.vanishing}}{3.2.1}{}\ifthenelse{\equal{#1}{nota--S.nochmal}}{3.0.1}{}\ifthenelse{\equal{#1}{nota--S}}{2.0.1}{}\ifthenelse{\equal{#1}{parahoric}}{4.2.3}{}\ifthenelse{\equal{#1}{prop--ballaballa}}{3.2.22}{}\ifthenelse{\equal{#1}{prop--boxtimes}}{2.4.4}{}\ifthenelse{\equal{#1}{prop--cells.flag}}{4.3.9}{}\ifthenelse{\equal{#1}{prop--change.facet}}{4.3.13}{}\ifthenelse{\equal{#1}{prop--DM.G.homotopy.invariant}}{2.2.11}{}\ifthenelse{\equal{#1}{prop--DTM.Fl.characterization}}{5.2.2}{}\ifthenelse{\equal{#1}{prop--DTM.G}}{3.1.27}{}\ifthenelse{\equal{#1}{prop--DTM.G/H}}{3.1.23}{}\ifthenelse{\equal{#1}{prop--equivariant.MTM}}{3.2.20}{}\ifthenelse{\equal{#1}{prop--existence.functors}}{6.3.3}{}\ifthenelse{\equal{#1}{prop--f_!.ind-Artin}}{2.3.3}{}\ifthenelse{\equal{#1}{prop--generators.DTM.G}}{3.2.23}{}\ifthenelse{\equal{#1}{prop--MTM.G}}{3.2.15}{}\ifthenelse{\equal{#1}{prop--Schubert.base.change}}{4.4.3}{}\ifthenelse{\equal{#1}{prop--sheafification.iso}}{2.2.25}{}\ifthenelse{\equal{#1}{prop--unipotent}}{A.4.6}{}\ifthenelse{\equal{#1}{prop--unseparated}}{2.1.14}{}\ifthenelse{\equal{#1}{prop--vector.extension}}{A.4.9}{}\ifthenelse{\equal{#1}{Quasi_Coxeter}}{4.2.13}{}\ifthenelse{\equal{#1}{reduced.eq}}{A.3.4}{}\ifthenelse{\equal{#1}{rema--bounded_subsets}}{4.2.6}{}\ifthenelse{\equal{#1}{rema--classical.equivariant}}{3.2.19}{}\ifthenelse{\equal{#1}{rema--classical.t-structure}}{3.2.7}{}\ifthenelse{\equal{#1}{rema--DM.prestacks}}{2.2.2}{}\ifthenelse{\equal{#1}{rema--DTM.large}}{3.1.9}{}\ifthenelse{\equal{#1}{rema--explain.stratified.Ind.scheme}}{3.1.3}{}\ifthenelse{\equal{#1}{rema--general.cellular}}{4.2.8}{}\ifthenelse{\equal{#1}{rema--l.adic.prestacks}}{2.3.11}{}\ifthenelse{\equal{#1}{rema--prestacks.examples}}{2.2.5}{}\ifthenelse{\equal{#1}{rema--realization.functor}}{3.2.9}{}\ifthenelse{\equal{#1}{rema--WT.maps}}{3.1.16}{}\ifthenelse{\equal{#1}{rema--WT.properties.fusion}}{6.3.7}{}\ifthenelse{\equal{#1}{rema--WT.smooth.duality}}{3.1.14}{}\ifthenelse{\equal{#1}{Root_Groups}}{4.3.10}{}\ifthenelse{\equal{#1}{Schubert_Map_Rel}}{4.4.2}{}\ifthenelse{\equal{#1}{Schubert_Map}}{4.3.5}{}\ifthenelse{\equal{#1}{Schubert_praesi_rel}}{4.4.4}{}\ifthenelse{\equal{#1}{Schubert_praesi}}{4.3.8}{}\ifthenelse{\equal{#1}{sect--algebraic.grps}}{A.2}{}\ifthenelse{\equal{#1}{sect--chevalley.base.change}}{4.4}{}\ifthenelse{\equal{#1}{sect--DM.Artin}}{2.3}{}\ifthenelse{\equal{#1}{sect--DM.ind-schemes}}{2.4}{}\ifthenelse{\equal{#1}{sect--DM.prestacks}}{2.2}{}\ifthenelse{\equal{#1}{sect--DM.schemes}}{2.1}{}\ifthenelse{\equal{#1}{sect--DM}}{2}{}\ifthenelse{\equal{#1}{sect--DTM.definitions}}{3.1}{}\ifthenelse{\equal{#1}{sect--DTM.double}}{5.3}{}\ifthenelse{\equal{#1}{sect--DTM.Fl}}{5}{}\ifthenelse{\equal{#1}{sect--DTM}}{3}{}\ifthenelse{\equal{#1}{sect--ind.schemes}}{A.1}{}\ifthenelse{\equal{#1}{sect--intersection}}{6}{}\ifthenelse{\equal{#1}{sect--invariant}}{6.2}{}\ifthenelse{\equal{#1}{sect--loop.definitions}}{4.2}{}\ifthenelse{\equal{#1}{sect--loop.group.dfn}}{4.1}{}\ifthenelse{\equal{#1}{sect--loop.grps}}{4}{}\ifthenelse{\equal{#1}{sect--MTM}}{3.2}{}\ifthenelse{\equal{#1}{sect--pro.action}}{A.3}{}\ifthenelse{\equal{#1}{sect--realization}}{2.1.2}{}\ifthenelse{\equal{#1}{sect--shtukas.def}}{6.3}{}\ifthenelse{\equal{#1}{sect--Stratifications.flag}}{4.3}{}\ifthenelse{\equal{#1}{sect--Tate.Fl}}{5.2}{}\ifthenelse{\equal{#1}{sect--torsors}}{A.4}{}\ifthenelse{\equal{#1}{sect--WT.Fl}}{5.1}{}\ifthenelse{\equal{#1}{sheaf_condition}}{2.2.14}{}\ifthenelse{\equal{#1}{shtuka.inv}}{6.3.2}{}\ifthenelse{\equal{#1}{shtuka.pres}}{6.3.1}{}\ifthenelse{\equal{#1}{strat.map.square}}{3.1.2}{}\ifthenelse{\equal{#1}{stratified.diagram}}{3.1.7}{}\ifthenelse{\equal{#1}{stratified.Ind-pres}}{3.1.4}{}\ifthenelse{\equal{#1}{syno--motives}}{2.1.1}{}\ifthenelse{\equal{#1}{Tate.big.orbit}}{5.1.3}{}\ifthenelse{\equal{#1}{theo--D.et!.sheaf}}{2.1.15}{}\ifthenelse{\equal{#1}{theo--descent.prestacks}}{2.2.16}{}\ifthenelse{\equal{#1}{theo--DM.descent}}{2.1.13}{}\ifthenelse{\equal{#1}{theo--equivariant.Chow}}{2.2.10}{}\ifthenelse{\equal{#1}{theo--equivariant.DTM.flag}}{5.3.4}{}\ifthenelse{\equal{#1}{theo--f!.Artin}}{2.3.7}{}\ifthenelse{\equal{#1}{theo--Fl.WT}}{5.1.1}{}\ifthenelse{\equal{#1}{theo--generators.DTM.flag}}{5.2.3}{}\ifthenelse{\equal{#1}{theo--motives.Ind-schemes}}{2.4.2}{}\ifthenelse{\equal{#1}{theo--simple.objects}}{3.3.8}{}\ifthenelse{\equal{#1}{tocindent-1}}{0pt}{}\ifthenelse{\equal{#1}{tocindent0}}{0pt}{}\ifthenelse{\equal{#1}{tocindent1}}{66.11127pt}{}\ifthenelse{\equal{#1}{tocindent2}}{0pt}{}\ifthenelse{\equal{#1}{tocindent3}}{0pt}{}\ifthenelse{\equal{#1}{torsor.lim}}{A.4.1}{}\ifthenelse{\equal{#1}{torsor.lim1}}{A.4.2}{}\ifthenelse{\equal{#1}{unif.unab}}{6.1.6}{}\ifthenelse{\equal{#1}{vector_bundle}}{A.4.11}{}\ifthenelse{\equal{#1}{Weil_resitriction}}{4.2.5}{}}
\def \inter#1{\cite[#1]{RicharzScholbach:Intersection}}

\newif\ifusebiber 
\usebiberfalse

\ifusebiber
  \usepackage[backref,style=alphabetic,sorting=anyt,language=english,backend=biber]{biblatex}

  \DeclareSourcemap{ 
    \maps[datatype=bibtex]{
      \map{
	\step[fieldset=issn, null]
      }
      \map{
	\step[fieldset=doi, null]
      }
      \map{
	\step[fieldset=isbn, null]
      }
    }
  }
  \bibliography{bib}

\else
  \def\nopp{} 
\fi

\begin{document}

\title{The motivic Satake equivalence
}
\author{Timo Richarz, Jakob Scholbach*}

\begin{abstract}
We refine the geometric Satake equivalence due to Ginzburg, Beilinson--Drinfeld, and Mirkovi\'c--Vilonen to an equivalence between mixed Tate motives on the double quotient $L^+ G \bsl LG / L^+ G$ and representations of Deligne's modification of the Langlands dual group $\widehat G$.\bigskip\\
\end{abstract}

\thanks{*Research of T.R.~partially funded by the Deutsche Forschungsgemeinschaft (DFG, German Research Foundation) - 394587809. Resarch of J.S.~funded by DFG, Sonderforschungsbereich 878 and DFG Cluster of Excellency ``Mathematics Münster''.}

\maketitle

\setcounter{tocdepth}{1}
\tableofcontents

\section{Introduction}

\subsection{Motivation and goals}
Split reductive groups are classified by their root data.
These come in pairs, consisting of a root datum and its associated dual root datum.
Accordingly, to every split reductive group $G$, there is associated its (Langlands) dual group $\widehat G$.

The work of Kazhdan and Lusztig \cite{KazhdanLusztig:Hecke, KazhdanLusztig:Poincare} shows that the representation theory of $\widehat G$ is closely related to the singularities arising in certain orbit closures inside (affine) flag varieties associated to $G$.
Building upon \cite{Lusztig:Singularities}, the work of Ginzburg \cite{Ginzburg:Perverse}, Belinson-Drinfeld \cite{BeilinsonDrinfeld:Quantization} and Mirkovi\'c--Vilonen \cite{MirkovicVilonen:Geometric} revealed an equivalence of symmetric tensor categories between the category of finite-dimensional $\widehat G$-representations and the category of certain sheaves on an infinite-dimensional variety $\Gr_G$ known as the \emph{affine Grassmannian} of $G$.
This categorical equivalence is called the geometric Satake equivalence.
It is an important tool in geometric representation theory which appears in different contexts and has a wide range of applications.
For further details on the subject, the reader may refer to the notes of Baumann and Riche \cite{BaumannRiche:Satake} and of Zhu \cite{Zhu:Introduction}, to \cite{Richarz:New, RicharzZhu:Ramified} for the relation with the classical Satake isomorphism (for which see \cite{Gross:Satake}), and to \cite{Zhu:Affine} for a Satake equivalence in the case of mixed characteristic.

The goal of the present manuscript is to provide a motivic refinement of the geometric Satake equivalence.
This has both philosophical and concrete consequences:
the above papers devoted to the Satake equivalence use different base schemes, and also use different cohomology theories.
It is therefore desirable to describe the common content of such different approaches, which is a goal accomplished in this paper.
As far as concrete applications are concerned, let us point out that one of our main motivations is the work of V.~Lafforgue \cite{Lafforgue:Chtoucas} on the Langlands parametrization for global function fields. V.~Lafforgue in particular conjectures \cite[Conj.~12.12]{Lafforgue:Chtoucas} that this parametrization is of motivic origin independent of an auxiliary prime number $\ell$ coming from the use of $\ell$-adic \'etale cohomology.
A first evidence for Lafforgue's conjecture is the construction of intersection cohomology motives on moduli stacks of $G$-shtukas alias $\IC$-Chow groups in \cite{RicharzScholbach:Intersection}.
The motivic Satake equivalence established in this paper is a second step of an ongoing project whose goal is to provide a motivic approach to V.~Lafforgue's Langlands parametrization.

\subsection{Results}
Let $G$ be a Chevalley group over $\bbZ$ (=split reductive group scheme \cite{Conrad:Groups}), and fix $T\subset B\subset G$, a split maximal torus contained in a Borel subgroup over $\bbZ$.
The \emph{loop group} of $G$ is the group-valued functor on the category of rings $R$ given by $LG(R)=G(R\rpot{\varpi})$.
Its subgroup functor $L^+G(R)=G(R\pot{\varpi})$ is the positive loop group.
Here $R\pot{\varpi}\subset R\rpot{\varpi}$ denotes the ring of power series in a formal variable $\varpi$, contained in its Laurent series.
For every finite field $\bbF_q$, the classical Satake isomorphism \cite{Gross:Satake} is an isomorphism of $\bbQ(\sqrt q)$-algebras
\begin{equation}\label{classical.satake.iso}
\calC_c\left (L^+G(\bbF_q)\bsl LG(\bbF_q)/L^+G(\bbF_q);\bbQ(\sqrt q)\right )\;\simeq\; R_{\widehat G}\otimes_\bbQ\bbQ(\sqrt q),
\end{equation}
where $\sqrt q$ is a fixed square root of $q$ needed in the construction.
The left hand side of \eqref{classical.satake.iso} are thus $\bbQ(\sqrt q)$-valued functions supported on finitely many double cosets.
The convolution of such functions turns the left hand side into an algebra known as the {\it spherical Hecke algebra}.
On the right hand side of \eqref{classical.satake.iso} the group $\widehat G$ is the Langlands dual group of $G$ formed over $\bbQ$ (with respect to a fixed pinning).
Then $R_{\widehat G}$ is the Grothendieck $\bbQ$-algebra of the category of representations of $\widehat G$ on finite-dimensional $\bbQ$-vector spaces.
Its ring structure is given by the tensor product of representations.
Writing $V_\mu$ for the simple $\widehat G$-representation of highest weight $\mu$, where $\mu \in X_*(T)^+$ is a dominant cocharacter, their classes $[V_\mu]$ form a $\Q$-basis of $R_{\widehat G}$.
Under \eqref{classical.satake.iso}, these correspond to functions which are related to the singularities of an infinite-dimensional space as follows.
The affine Grassmannian is the \'etale sheaf quotient
\[
\Gr_G\defined \big(LG/L^+G\big)^\et,
\]
 which is representable by an ind-projective ind-scheme (=infinite union of projective $\bbZ$-schemes) equipped with a left action of $L^+G$.
For each dominant cocharacter $\mu\in X_*(T)^+$, we denote by $\Gr^{\leq \mu}_G$ the scheme-theoretic image of the orbit map $L^+G\to \Gr_G, g\mapsto g\cdot \varpi^\mu\cdot e$ where $e\in \Gr_G(\bbZ)$ is the base point. Then $\Gr^{\leq \mu}_G\to \Spec(\bbZ)$ is a projective scheme, usually singular, which contains the open smooth $L^+G$-orbit $\Gr_G^\mu\subset \Gr_G^{\leq \mu}$ as a fiberwise dense open subscheme.
There is a presentation on the underlying reduced locus
\[
\big(\Gr_G\big)_\red=\colim_{\mu \in X_*(T)^+}\Gr_G^{\leq \mu}.
\]
For a finite field $\bbF_q$ and each auxiliary prime $\ell\nmid q$, let $\IC_{\mu, q, \ell}$ be the $\ell$-adic intersection complex of $\Gr_G^{\leq \mu}\otimes_\bbZ\bbF_q$ in the sense of Goresky-MacPherson-Deligne.
A surprising observation of \cite{Lusztig:Singularities} is that the class $[V_\mu]$ corresponds under \eqref{classical.satake.iso} (up to a power of $\sqrt q$) to the trace of Frobenius function of $\IC_{\mu, q, \ell}$ given by Grothendieck's sheaf function dictionary.

The geometric Satake equivalence is a categorification of \eqref{classical.satake.iso}. It is known in several settings using different cohomology theories:
in \cite{Ginzburg:Perverse, MirkovicVilonen:Geometric, BaumannRiche:Satake} the authors work with $\Gr_G\otimes_\bbZ \bbC$ using Betti cohomology (the latter two with more general coefficients however), whereas \cite{BeilinsonDrinfeld:Quantization} works with $\Gr_G\otimes_\bbZ \bbC$ using $D$-modules, and \cite{Richarz:New, RicharzZhu:Ramified} works with $\Gr_G\otimes_\bbZ k$ for general fields $k$ using $\ell$-adic \'etale cohomology.
Here we provide a motivic refinement.

In analogy with the left hand side of \eqref{classical.satake.iso} we consider the double quotient $L^+G\bsl LG/L^+G\to \Spec(\bbZ)$ viewed as a groupoid-valued functor on the category of rings.
For each such functor we have constructed in \cite{RicharzScholbach:Intersection} a category of motives (with rational coefficients)
\begin{equation}\label{DM.intro}
\DM\big(L^+G\bsl LG/L^+G\big)=\DM\big(L^+G\bsl LG/L^+G;\bbQ\big).
\end{equation}
The collection of all such categories is equipped with a Grothendieck six functor formalism (with certain restrictions on the $*$-pullback).
The construction in \emph{op.~cit.~}builds upon the recent advances in the theory of motivic sheaves due to Ayoub \cite{Ayoub:Six1, Ayoub:Six2, Ayoub:Realisation} and Cisinski-Déglise \cite{CisinskiDeglise:Triangulated, CisinskiDeglise:Etale} as envisioned by Beilinson.

In there we consider a much smaller full subcategory of {\it stratified Tate motives}
\[
\DTM\big(L^+G\bsl LG/L^+G\big)\subset \DM\big(L^+G\bsl LG/L^+G\big).
\]
This category is generated by all Tate twists $1_{\Gr^\mu_G}(n)$, $n\in \bbZ$ of the characteristic relative motives of the orbits $\Gr_G^\mu$, $\mu\in X_*(T)^+$, and is well suited for applications to Hecke algebras.
It is equipped with a convolution product: namely, for two motives $A,B$ on the double quotient their convolution is defined as the motive
$$A\star B=m_!p^!(A\boxtimes B)\eqlabel{define.convolution}$$
using the maps
$$
\xymatrix{
L^+G \backslash LG / L^+G \x L^+G \backslash LG / L^+G &
L^+G \backslash LG \x^{L^+G} LG / L^+G \ar[l]_(.45)p \ar[r]^(.55)m &
L^+G \backslash LG / L^+G,}
$$
where $p$ is the canonical projection and $m$ is induced by the multiplication $LG\x LG\to LG$, $(g_1,g_2)\mapsto g_1\cdot g_2$.
The convolution of motives in \eqref{DM.intro} is modeled on the convolution in the spherical Hecke algebra \eqref{classical.satake.iso}.
(The use of the functor $p^!$ in \refeq{define.convolution}, as opposed to $p^*$, is related to the construction of \eqref{DM.intro}. Elements of this category should be more appropiately thought of as ``measures'' instead of ``functions'' which leads in a categorical setting to the use of $!$-pullback instead of $*$-pullback.)

The fibers of convolution morphisms are paved by cells \cite{Haines:Purity} which leads to the following result.
\medskip\\
 {\bf Theorem A.} {\it If $A,B$ are stratified Tate motives \textup{(}resp.~and pure of some weight\textup{)}, then $A\star B$ is again stratified Tate  \textup{(}resp.~and pure\textup{)}.} (\ref{convolution.Fl.Tate}, \ref{convolution.DTM.weights})
 \medskip\\
In fact, we show a more general version of Theorem A where $L^+G$ is replaced by an arbitrary parahoric subgroup of $LG$.
As a consequence, the category of stratified Tate motives is equipped with a monoidal structure with respect to the convolution.
We now cut out an abelian subcategory as follows.

By \cite{RicharzScholbach:Intersection}, which extends the work of Soergel and Wendt \cite{SoergelWendt:Perverse}, for each $\mu\in X_*(T)^+$, $n\in \bbZ$ there exists an {\it intersection motive}
\begin{equation}\label{intersection.intro}
\IC_{\mu,\bbZ}(n)\in \DTM\big(L^+G\bsl LG/L^+G\big)
\end{equation}
which is supported on $\Gr_G^{\leq \mu}$ and such that $\IC_{\mu,\bbZ}(n)|_{\Gr_G^\mu}=1_{\Gr_G^\mu}(n)$.
We emphasize that the non-trivial construction of $\IC_{\mu,\bbZ}(n)$ in \cite{RicharzScholbach:Intersection} bypasses the use of standard conjectures on $t$-structures of triangulated categories of motives.
For any finite field $\bbF_q$ and each prime $\ell\nmid q$, its base change $\IC_{\mu,\bbF_q}:=\IC_{\mu,\bbZ}|_{\bbF_q}$ maps under the $\ell$-adic realization to the intersection complex $\IC_{\mu, q,\ell}$ on $\Gr_G^{\leq \mu}\otimes_\bbZ\bbF_q$ as above.
For the field of rational numbers $\bbQ$, its base change $\IC_{\mu,\bbQ}:=\IC_{\mu,\bbZ}|_\bbQ$ maps under the Betti realization to the intersection complex on the stratified topological space $\Gr_G^{\leq \mu}(\bbC)$.
Thus, the motives $\IC_{\mu, \Z}$ interpolate between various intersection sheaves arising in the literature on the Satake equivalence, both in the sense of letting the base scheme vary, and also in the sense of varying the cohomology theory.

Again by \cite{RicharzScholbach:Intersection}, the category $\DTM(L^+G\bsl LG/L^+G)$ is equipped with a non-degenerate (perverse) motivic $t$-structure.
Its heart, the abelian subcategory of {\it mixed \textup{(}stratified\textup{)} Tate motives}
\[
\MTM\big(L^+G\bsl LG/L^+G\big)\subset \DTM\big(L^+G\bsl LG/L^+G\big),
\]
is generated by the intersection motives \eqref{intersection.intro}. 
Using motivic (global) cohomology provides a $\bbQ$-linear functor $\omega\co \MTM(L^+G\bsl LG/L^+G)\to \Vect_\bbQ$, see \thref{fiber.functor}.
\medskip\\
 {\bf Theorem B.} {\it
 i\textup{)} If $A, B$ are mixed Tate motives, then their convolution $A\star B$ is mixed Tate as well. \textup{(\ref{convolution.Tate.perverse})}\smallskip\\
 ii\textup{)} The motivic cohomology functor $\omega$ is $\bbQ$-linear, exact, faithful and equipped with functorial isomorphisms $\omega(A\star B)\simeq \omega(A)\otimes_\bbQ \omega(B)$ for all objects $A, B$. \textup{(\ref{pushforward.symmetric.monoidal})}\smallskip\\
 iii\textup{)} There exists a unique symmetric monoidal structure on $ \MTM(L^+G\bsl LG/L^+G)$ with respect to the convolution product $\star$ characterized by the property that $\omega$ is a tensor functor with respect to ii\textup{)} and the canonical symmetric monoidal structure on $(\Vect_\bbQ,\otimes)$. \textup{(\ref{constraints}, \ref{Tannaka_Cat})}
 }
\medskip\\

Already for $G=\GL_2$, the subcategory of {\it non-twisted} intersection motives $\IC_{\mu}$, $\mu\in X_*(T)^+$ is not stable under convolution.
As observed in \cite[Rmk.~2.10]{HeinlothNgoYun:Kloosterman} and \cite{RicharzZhu:Ramified} this phenomenon is linked to the presence of $\sqrt q$ in \eqref{classical.satake.iso} and leads to the appearance of a $\bbG_{m,\bbQ}$-extension on the dual side.

We consider the central subgroup {$\Bmu_2\subset \widehat G\x\bbG_{m,\bbQ}$} of order $2$ generated by the element $(\epsilon,-1)$ where $\epsilon:=(2\rho)(-1)\in \widehat G(\bbQ)$ and $2\rho$ denotes the sum of positive roots viewed as a cocharacter of $\widehat G$.
\emph{Deligne's modified Langlands dual group} (see \cite{Deligne:Letter2007}, \cite[\S 2]{FrenkelGross:Rigid}, \cite[\S 5]{BuzzardGee:Conjectures}, and \cite[\nopp 5.5.14]{Zhu:Introduction}) is defined as the split reductive $\bbQ$-group
\[
\widehat G_1\defined \widehat G\x^{\displaystyle \Bmu_2} \bbG_{m,\bbQ}.
\]
The extra $\bbG_{m,\bbQ}$ factor corresponds to the occurrence of Tate twists when forming the convolution of intersection motives.
We denote by $\Rep_\bbQ(\widehat G_1)$ the category of algebraic $\widehat G_1$-representations on $\bbQ$-vector spaces.
This category is semi-simple. Its simple objects are labelled by $V_\mu(n)$ with $\mu\in X_*(T)^+$, $n\in \bbZ$.

To make the connection with \eqref{classical.satake.iso} we base change the groups $LG_{\bbF_q}:=LG\otimes_\bbZ\bbF_q$ and $L^+G_{\bbF_q}:=L^+G\otimes_\bbZ\bbF_q$ to a finite field.
Then the analogue of Theorem B holds with $\MTM(L^+G\bsl LG/L^+G)$ replaced by $\MTM(L^+G_{\bbF_q}\bsl LG_{\bbF_q}/L^+G_{\bbF_q})$.
\medskip\\
 {\bf Theorem C.} {\it
For each finite field $\bbF_q$, there is an equivalence of symmetric monoidal categories
 \[
\MTM\big(L^+G_{\bbF_q}\bsl LG_{\bbF_q}/L^+G_{\bbF_q}\big) \;\simeq\; \Rep_\bbQ(\widehat G_1), \;\; \IC_{\mu,\bbF_q}(n)\mapsto V_\mu(n),
 \]
 using the tensor structure from Theorem B.
 Under this equivalence, the motivic cohomology functor $\omega$ corresponds to the forgetful functor $\Rep_\bbQ(\widehat G_1)\to\Vect_\bbQ$.
 For each prime $\ell\nmid q$, this equivalence gives under the $\ell$-adic \'etale realization the geometric Satake equivalence as explained in \textup{\cite[\nopp 5.5.14]{Zhu:Introduction}}. \textup{(\ref{Satake}, \ref{finite.field.exam})}
 }
 \medskip\\
Among other things, Theorem C asserts that the left hand category is semi-simple.
This semi-simplicity is inferred, via Lusztig's parity vanishing, from the semi-simplicity of the abelian category $\MTM(\Fq)$.
The latter semi-simplicity holds since higher algebraic $K$-theory of $\Fq$ is torsion by Quillen's computation, see \thref{finite.field.exam}.
This semi-simplicity is then lifted to the mixed Tate motives on the double quotient over a finite field.
Passing to the trace of the Frobenius function as in, say, \cite{Cisinski:SurveyCoho} one recovers the Satake isomorphism similar to $\eqref{classical.satake.iso}$ where one now considers $\bbQ$-valued functions and a quotient of the representation ring $R_{\widehat G_1}$, cf.~\S\ref{functions.sec}.

In contrast to $\MTM(\Fq)$, the categories $\MTM(\bbZ)$ and, a fortiori, $\MTM(L^+G\bsl LG/L^+G)$ are no longer semi-simple.
More generally, if $S$ is a sufficiently nice scheme which satisfies the Beilinson--Soul\'e vanishing (e.g.~the spectrum of finite fields as above; number fields or their rings of integers; function fields over a finite field or their rings of integers; or filtered colimits of these rings) the category of mixed Tate motives
\begin{equation}\label{MTM.general.intro}
\MTM\big(L^+G_S\bsl LG_S/L^+G_S\big)
\end{equation}
is well-defined and satisfies Theorem B where we denote $L^{(+)}G_S:=L^{(+)}G\x_{\Spec(\bbZ)}S$.
In the category \eqref{MTM.general.intro} we also have the intersection motives $\IC_{\mu,S}(n)$ for $\mu\in X_*(T)^+$, $n\in \bbZ$.
We denote by $\Sat_{G,S}$ the full semi-simple subcategory of \eqref{MTM.general.intro} generated by the intersection motives by means of direct sums.
This subcategory $\Sat_{G,S}$ is stable under convolution, and hence inherits a symmetric monoidal structure.
\medskip\\
 {\bf Theorem D.} {\it Let $p\co S\to \Spec(\bbZ)$ be a base scheme as above. \smallskip\\
 i\textup{)} The pullback of motives induces an equivalence of symmetric monoidal categories
 \[
 \Sat_{G,\bbZ}\to \Sat_{G,S}, \;\; \IC_{\mu,\bbZ}(n)\mapsto p^*\IC_{\mu,\bbZ}(n)=\IC_{\mu,S}(n),
 \]
 and hence $\Sat_{G,S}\simeq \Rep_\bbQ(\widehat G_1)$ by Theorem C independently of $S$. \textup{(\ref{Tannaka.Sat})}\smallskip\\
 ii\textup{)} Let $\scrU_S$ be the pro-unipotent algebraic $\bbQ$-group arising from extensions in the category $\MTM(S)$.
 Then there is an equivalence of symmetric monoidal categories
 \[
 \MTM\big(L^+G_S\bsl LG_S/L^+G_S\big)\;\simeq\; \Rep_\bbQ(\scrU_S\rtimes \widehat G_1)
 \]
 where on the right-hand side is the category of representations of the pro-algebraic group $\scrU_S\rtimes \widehat G_1$ on $\bbQ$-vector spaces.
 \textup{(\ref{full.Tannaka}, \ref{Satake.Ind.objects})}
 }
 \medskip\\
Part i) of Theorem D precisely formulates the experimental fact that under the geometric Satake isomorphism the dual side does not depend on the base scheme over which the affine Grassmannian is defined.
Part ii) is, in part, an extension of Levine's work \cite{Levine:Tate} which one recovers in the special case where $G$ is the trivial group.

\subsection{Related and future work}

Zhu \cite{Zhu:Geometric} has sketched the construction of a motivic Satake equivalence over $\bbF_q$ using the category of numerical motives of Jannsen.
Zhu's approach is based on an explicit enumeration of algebraic cycles on affine Grassmannians. By comparison, the approach taken in this paper is more strongly relying on the general framework of motives, which we expect to be fruitful also for our upcoming work.

One may imagine using the theory of Nori motives to produce an abelian category of motives related to the Satake equivalence.
Nori motives, however, depend upon the cohomology theory chosen at the outset. In the case of motives over $\Fq$, say, this would in practice mean choosing $\ell$-adic cohomology for some $\ell$ prime to $q$.
Again, the choice of working with motives as developed by Ayoub and Cisinski--Déglise is based on the desire to apply it to a Langlands parametrization over function fields, where we precisely seek to avoid a reference to $\ell$-adic cohomology.

Throughout this paper, motives have rational coefficients. 
Using upcoming work of Spitzweck on t-structures on Tate motives with integral coefficients, it would be very interesting to establish a Satake equivalence in this situation.
The reader is referred to \cite{Zhu:IntegralSatake} for a result on the level of functions. 

There are versions of the geometric Satake equivalence using different affine Grassmannians such as the Witt vector (or $p$-adic) affine Grassmannian of Zhu \cite{Zhu:Affine}, the $B_{\on{dR}}$-affine Grassmannian of Fargues-Scholze \cite{FarguesScholze:Geometrization}.
In subsequent work \cite{RicharzScholbach:MotivicWitt}, we extend the methods of this paper to cover a Satake equivalence for Witt vector Grassmannians.

As was stated above, we conceive the results in \cite{RicharzScholbach:Intersection} and the Satake equivalence in this paper to be two steps in a long-term program aiming to prove a motivic version of V.~Lafforgue's Langlands parametrization over function fields.
The immediate next step, to be addressed in a subsequent paper, is to improve on Theorem C by proving a motivic version of Gaitsgory's factorization (or fusion) version of the geometric Satake equivalence \cite{Gaitsgory:deJong}.
This will require suitable Whitney--Tate properties of Beilinson--Drinfeld Grassmannians, as opposed to the affine Grassmannian $\Gr_G$ encountered above.
Here the six functor formalism for the categories of motives mentioned in Theorem C will be crucial.
Further steps in this program include a motivic Drinfeld lemma, a motivic construction of excursion operators, and their identification with Hecke operators. All these remain to be done as well.

\bigskip
\noindent \textit{Acknowledgements.} We thank Dennis Gaitsgory, Fritz Hörmann, and Thomas Nikolaus for helpful discussions, and the anonymous referee for many suggestions that improved the quality of the manuscript.
The authors thank the University of M\"unster, Harvard University, Deutsche Forschungsgemeinschaft (DFG, German Research Foundation under Germany's Excellence Strategy EXC 2044–390685587, Mathematics Münster: Dynamics–Geometry–Structure), the Institut de Math\'ematiques de Jussieu and the Technische Universit\"at Darmstadt for financial and logistical support which made this research possible.

\section{Motives on affine flag varieties}
\label{sect--motive.affine.flag.basics}
In this section, we recollect and extend some material from \cite{RicharzScholbach:Intersection} as is needed throughout this manuscript.
In \S\ref{sect--loop.grps.Satake}, we state some facts on loop groups and their affine flag varieties.
The next \S\ref{sect--motives.prestack} treats motives on prestacks which is applied in \S\S\ref{sect--stratified.motives}--\ref{sect--changing.base} to affine flag varieties.
\refsect{parity.vanishing} gathers some facts pertaining to Kazhdan--Lusztig parity vanishing.

\nota\thlabel{base.scheme}
Throughout this manuscript, $S$ is an irreducible, regular scheme which is separated of finite type over a Noetherian, excellent, separated and at most $2$-dimensional scheme. Further, we assume that $S$ satisfies the Beilinson-Soulé vanishing conjecture (cf.~\inter{(\iref{eqn--BS.vanishing})}), and admits an $\ell$-adic realization functor in the sense of \inter{\S\iref{sect--realization}, Rmk.~\iref{rema--realization.functor}}.

Examples include finite fields, number fields and function fields of curves over finite fields, their rings of algebraic integers and filtered colimits of these rings.
\xnota

\subsection{Loop Groups and their affine flag varieties}
\label{sect--loop.grps.Satake} 
We refer the reader to \inter{\S\iref{sect--loop.grps}} for further details and references on the following material.

We denote by $\AffSch_S$ the category of affine schemes $\Spec(R)\to S$ equipped with a map to $S$. Let $G$ be a split reductive group scheme over $S$, for example $G = \GL_{n,S}$.
The \emph{loop group} $LG$ is the presheaf
$$LG\co \AffSch_S^\opp \r \Groups,\;\; \Spec(R) \mapsto G(R\rpot{\varpi}),$$
where $R\rpot{\varpi}$ is the ring of Laurent series in the formal variable $\varpi$.
It is represented by an ind-affine ind-scheme over $S$, and in particular $LG$ is an fpqc sheaf on $\AffSch_S^{\on{op}}$.

We fix $T\subset B\subset G$ over $S$, a split maximal torus contained in a Borel subgroup.
Let $\scrA=\scrA(G, B,T)$ be the standard apartment with origin $0$ defined by $G$ and standard alcove $\bba$ defined by $B$.
We only consider facets $\bbf\subset \scrA$ which are contained in the closure of $\bba$.
Attached to $\bbf$ is the parahoric subgroup $\calP_\bbf \subset LG$ which is an $S$-affine, $S$-flat closed subgroup scheme.
For this paper, the most important case is $\bbf = 0$, in which case $\calP_0 =: L^+ G$ is the \emph{positive loop group} given by the presheaf
$$L^+G\co \AffSch_S^\opp \r \Groups,\;\; \Spec(R) \mapsto G(R\pot{\varpi}).$$
If $\bbf=\bba$, then $\calP_\bba=:\calB$ is the standard Iwahori subgroup defined as the preimage of $B$ under the map $L^+G\to G$, $\varpi\mapsto 0$.

The étale sheafification of the quotient $\Fl_\bbf:=(LG / \calP_\bbf)^\et$ is called the {\it partial affine flag variety associated with $\bbf$}.
It is represented by an ind-projective ind-scheme over $S$.
For $\bbf=0$, it is denoted by $\Gr=\Gr_G$, and called the \emph{affine Grassmannian}.

Given two facets $\bbf', \bbf\subset \bar{\bba}\subset \scrA$, the orbits of the $\calP_{\bbf'}$-left-action on $\Fl_\bbf$ are enumerated by the double quotient $W_{\bbf'}\bsl W/W_\bbf$ of the Iwahori-Weyl (or extended affine Weyl) group $W=W(G,T)$ by the subgroups $W_{\bbf'}, W_\bbf\subset W$ generated by the reflections preserving $\bbf'$ resp.~$\bbf$. The choice of $\bba$ defines a length function $l=l(\bbf',\bbf)\co W_{\bbf'}\bsl W/W_\bbf \to \bbZ_{\geq 0}$, and a Bruhat partial order $\leq$ on the double coset. For each $w\in W_{\bbf'}\bsl W/W_\bbf$, the locally closed immersion of the $\calP_{\bbf'}$-orbit of $w$ is denoted by
\begin{equation}\label{orbit.inclusion}
\iota_w\co \Fl_\bbf^w \;\stackrel {j_w} {\hookrightarrow} \; \Fl_\bbf^{\leq w}\;\stackrel {i_w}{\hookrightarrow}\; \Fl_\bbf.
\end{equation}
Then $\Fl^{\leq w}_\bbf\to S$ is a proper scheme called the (affine) Schubert scheme.
It contains $\Fl_\bbf^w$ as an open $S$-smooth subscheme which is fibrewise dense and which is called the (affine) Schubert cell.
For each map $\Spec(k)\to S$ from a field, the base change $\Fl^{w}_\bbf\x_{S}\Spec(k)\subset \Fl^{\leq w}_\bbf\x_{S}\Spec(k)$ identifies on the underlying reduced locus with the Schubert cell, resp.~Schubert variety
over $k$ attached to the class $w$ and the $k$-group scheme $G\x_S\Spec(k)$.
Further, each $\Fl^w_\bbf\to S$ is pure of relative dimension $l(w)$ and $\Fl^{v}_\bbf\subset \Fl^{\leq w}_\bbf$ if and only if $v\leq w$ for $v,w \in W_{\bbf'}\bsl W/W_\bbf$. If $\bbf'=\bbf=0$, then $W_{0}\bsl W/W_0=X_*(T)^+$ is the partially ordered set of dominant cocharacters with length function $l=l(0,0)\co X_*(T)^+\to \bbZ_{\geq 0}, \mu\mapsto \lan 2\rho,\mu\ran$ where $\rho$ denotes the half sum of the $B$-positive roots and $\lan\str,\str\ran\co X^*(T)\x X_*(T)\to \bbZ$ is the natural pairing. In the case of the affine Grassmannian, we denote the orbits (resp.~orbit closures) by $\Gr^\mu\subset \Gr^{\leq \mu}$ for $\mu\in X_*(T)^+$.

\subsection{Motives on prestacks} \label{sect--motives.prestack}
We refer the reader to \inter{\S\iref{sect--DM}} for further details and references on the following material.

We consider the triangulated category of motives with rational coefficients
\[
\DM(X)\defined \DM(X,\bbQ), \;\; X\in \Sch_S^{\on{ft}},
\]
where $\Sch_S^{\on{ft}}$ is the category of finite type schemes over $S$.
This category is denoted by $\on{DA}(X, \bbQ)$ in \cite{Ayoub:Realisation} and by $\D_{\A^1, \et}(X, \Q)$ in \cite{CisinskiDeglise:Triangulated}.
Categories of motives with rational coefficients admit a full six functor formalism: there are pairs of adjoint functors $(f^*,f_*)$, $(f_!,f^!)$ for a map $f\in \Sch_S^{\on{ft}}$ and $(\str\otimes\str, \Hom(\str,\str))$ satisfying the usual compatibilities such as smooth/proper base change, Poincar\'e duality, K\"unneth/projection formula etc.
Following Hoyois \cite{Hoyois:Six} and Khan \cite{Khan:Motivic}, this can be upgraded to a presheaf of \ii-categories
\[
\DM^!\co (\Sch_S^{\ft})^\opp\to \DGCat_\cont, \;\; X\mapsto \DM(X), \;\; f\mapsto f^!,\eqlabel{DM.Sch.ft}
\]
where $\DGCat_\cont$ is the category of presentable, stable, $\bbQ$-linear, dg-\ii-categories with {\em colimit-preserving} functors.
The \ii-category $\DGCat_\cont$ is complete and cocomplete, i.e., admits all (homotopy) limits and (homotopy) colimits, so that the following Kan extensions are available.

\defi \thlabel{prestack.DM}
\begin{enumerate}
\item
Let $\AffSch_S^\ft\subset \AffSch_S$ be the full subcategory of objects of finite type over $S$.
Throughout, we will replace this category by a small skeleton containing the objects of interest to us.

\item
Fix some regular cardinal $\kappa$, and let $\AffSch_S^\kappa := \Pro_{\kappa\text{-small}}(\AffSch_S^\ft)$ be the category of $\kappa$-small pro-objects in $\AffSch_S^\ft$.

\item
The \ii-category of \emph{prestacks} is defined as $\PreStk_S^\kappa := \Fun((\AffSch_S^\kappa)^\opp,\text{\ii-}\Gpd)$
where \ii-$\Gpd$ is the \ii-category of \ii-groupoids (also called spaces).
Hereafter, we will usually drop the $\kappa$ from the notation so that $\PreStk := \PreStk_S^\kappa$, $\AffSch_S := \AffSch_S^\kappa$.

\item
Define the functor
$$\DM^! : \AffSch_S^\opp \r \DGCat_\cont$$
to be the left Kan extension of the functor $\DM^!$ in \refeq{DM.Sch.ft} along the inclusion $\AffSch_S^\ft \subset \AffSch_S$.

\item
Define the functor
$$\DM^! : \PreStk_S^\opp \r \DGCat_\cont\eqlabel{DM.prestacks}$$
to be the right Kan extension of the preceding functor along the Yoneda embedding $\AffSch_S \subset \PreStk_S$.
\end{enumerate}
\xdefi

We emphasize that $\DM^!$ in \refeq{DM.prestacks} encodes a category of motives $\DM(X)$ (with rational coefficients) for each prestack $X$, and for each map $f\co X\to Y$ in $\PreStk_S$ a colimit-preserving functor $f^!\co \DM(Y)\to \DM(X)$.
This definition follows the approach of Gaitsgory-Rozenblyum and Raskin.
We refer to \inter{§\iref{sect--DM.prestacks}} for references and also for further discussion of the definition.

\theo \thlabel{properties.DM}
i\textup{)} The presheaf $\DM^!\co \PreStk_S^\opp \to \DGCat_\cont$ is a sheaf in the \'etale topology. For each prestack $X\in \PreStk_S$ the \ii-sheafification $X\to X^\et$ induces an equivalence on categories of motives $\DM(X^\et)\stackrel \simeq \lr \DM(X)$, cf.~\inter{Thm.~\iref{theo--descent.prestacks}}.\smallskip\\
ii\textup{)} The restriction of the presheaf $\DM^!$ to the category of strict ind-schemes of ind-finite type over $S$ admits a full six functor formalism \textup{(}with certain restrictions on $f^*$\textup{)}, cf.~\inter{Thm.~\iref{theo--motives.Ind-schemes}}.
\xtheo

\subsection{Stratified Tate motives on affine flag varieties} \label{sect--stratified.motives}
We recall some material pertaining to stratified Tate motives, referring to \inter{\S\iref{sect--DTM.Fl}} for further details.

By virtue of \thref{prestack.DM}, there is the category of motives $\DM(\calP_{\bbf'}\bsl LG/\calP_\bbf)$ for each pair of facets $\bbf',\bbf \subset \bar\bba \subset \scrA$. Using \thref{properties.DM} i), we have a forgetful functor
$$ \DM(\calP_{\bbf'}\bsl LG / \calP_\bbf)\simeq \DM(\calP_{\bbf'}\bsl\Fl_\bbf) \to \DM(\Fl_\bbf),$$
which associates to each motive on the double quotient its underlying non-equivariant motive.

\defi
i) The {\em category of $(\bbf',\bbf)$-stratified Tate motives} $\DTM( \Fl_\bbf) \subset \DM(\Fl_\bbf)$ is the full subcategory consisting of objects $M\in \DM(\Fl_\bbf)$ such that for all $w\in W_{\bbf'}\bsl W/W_\bbf$
$$\iota_w^*M \in \DTM(\Fl_\bbf^w),$$
where $\iota_w$ is as in \eqref{orbit.inclusion} and $\DTM(\Fl_\bbf^w) \subset \DM(\Fl_\bbf^w)$ denotes the subcategory generated by $1(n)$, $n \in \Z$ under arbitrary shifts and colimits.
(This condition is equivalent to requiring $\iota_w^! M \in \DTM(\Fl_\bbf^w)$ for all $w$, cf.~\inter{Def.~3.1.11, Cor.~4.3.12}.)
\smallskip\\
ii) The \emph{category of $(\bbf',\bbf)$-stratified equivariant Tate motives}
$$\DTM(\calP_{\bbf'} \bsl LG / \calP_\bbf) \subset \DM(\calP_{\bbf'} \backslash LG / \calP_\bbf)$$
is the full subcategory consisting of objects $M$ whose underlying non-equivariant motive lies in $\DTM(\Fl_\bbf)$. In particular, there is a forgetful functor $\DTM(\calP_{\bbf'} \bsl LG / \calP_\bbf)\to \DTM(\Fl_\bbf)$.
\xdefi

The category $\DTM(\Fl_\bbf)\subset \DM(\Fl_\bbf)$ agrees by \inter{Thm.~\iref{theo--Fl.WT}} with the full subcategory generated by all $\iota_{w,!}1(n)$ (resp.~by all $\iota_{w,*}1(n)$), and so is well-suited for applications to Hecke algebras.

\theo \thlabel{motivic.t.structure}
\textup{(}\inter{Thm.~\iref{theo--equivariant.DTM.flag}}\textup{)}
Let $S$ be a scheme as in \thref{base.scheme}. \smallskip\\
i\textup{)} The category $\DTM(\Fl_\bbf)$ carries the so-called \emph{perverse motivic $t$-structure} whose heart is denoted by $\MTM(\Fl_\bbf)$. The subcategory $\MTM(\Fl_\bbf)^\comp$ of compact objects is Artinian and Noetherian. Its simple objects are precisely the intersection motives on the orbit closures
$$\IC_w(n)\defined (\iota_w)_{!*}1(n)[d_w] \defined i_{w,!}j_{w,!*}1(n)[d_w], \;\;\; n\in \bbZ, w\in W_\bbf\bsl W/W_\bbf$$
where $\iota_w=i_w\circ j_w$ is as in \eqref{orbit.inclusion} and $d_w$ is the relative dimension of $\Fl_\bbf^w$ over $S$.\smallskip\\
ii\textup{)} If $\bbf'=\bbf$, then the forgetful functor
\[
\DTM(\calP_\bbf \bsl LG / \calP_\bbf)\to \DTM(\Fl_\bbf),
\]
and the $t$-structure in i\textup{)} create a t-structure on the left hand category. The induced functor on the hearts $\MTM(\calP_\bbf \bsl LG / \calP_\bbf)\to \MTM(\Fl_\bbf)$ is fully faithful and induces a bijection on simple objects.
\xtheo

We point out the following results which are needed in \S\ref{sect--tate.aff.grass} below.

\lemm
\thlabel{orbit.simply.conn}
For $w\in W_{\bbf'}\bsl W/W_\bbf$, there is an equivalence
\[
\MTM(S)\overset{\simeq}{\lr} \MTM(\Fl_\bbf^w).
\]
\xlemm

\pf
The structure map $\Fl_\bbf^w\to S$ is smooth, and the Schubert cell admits a stratification into affine spaces, by virtue of the stratification in Iwahori orbits, cf.~\inter{Prop.~\iref{prop--cells.flag}}.
Hence, this lemma follows from the general \thref{pi*.equivalence.MTM} below.
\xpf

\lemm
\thlabel{pi*.equivalence.MTM}
Let $\pi\co X \r S$ be a smooth surjective map of schemes of relative dimension $d$ with connected fibers.
We assume that $X$ admits a stratification in the sense of \inter{Def.~\iref{defi--stratified.dfn}} by schemes of the form $\bbV(\calE)$ \textup{(}where $\calE$ is a vector bundle over $S$\textup{)}, e.g., by affine spaces over $S$.
Then there is an equivalence of categories
$$\pi^![-d] = \pi^*[d](d)\co \MTM(S) \overset{\simeq}{\lr} \MTM(X),$$
where the Tateness of motives on $X$ is with respect to the stratification by a single stratum.
\xlemm

\pf
By the conventions in \thref{base.scheme}, $S$ is connected and hence so is $X$ \StP{0378}.
The functor is fully faithful by \inter{Lem.~\iref{lemm--pi*.fullyfaithful.MTM}}.
For essential surjectivity,
we first claim that $\Hom_S(M, N[1]) = \Hom_X(\pi^* M, \pi^* N[1])$ for $M, N \in \MTM(S)$.
We prove this by induction on the number of strata in $X$.
If $X=\bbV(\calE)$ is a single stratum, then this holds even for all $N \in \DTM(S)$.
For the inductive step we use as in loc.~cit.~the localization sequence for a minimal stratum
$$Z = \bbV(\calE) \stackrel i \r X \stackrel j \gets U := X \bsl Z.$$
Let $\pi_Z := \pi \circ i$, $\pi_U := \pi \circ j$.
Since $X$ is connected and is assumed to have at least two strata, the codimension $c := \codim_X Z$ is positive.
By induction, the composite
$$\Hom_S(M, N[1]) \r \Hom_X(\pi^* M, \pi^*N[1]) \r \Hom_U(\pi_U^* M, \pi_U^* N[1])$$
is an isomorphism.
By the localization sequence, the kernel of the right hand map is mapped onto by $\Hom_Z(\pi_Z^* M, \pi_Z^* N(-c)[1-2c])$ which vanishes by the Beilinson--Soulé condition for $Z$ (equivalently, for $S$).
Hence the left hand map above is an isomorphism as well, showing our claim.

The generators $1(n)[d]$, $n\in \bbZ$ of $\MTM(X)$ trivially lie in the image of our functor, so we are done by using that $\pi^*$ is an isomorphism on the level of extensions of mixed Tate motives by the above claim.
\xpf

\subsection{Changing the base scheme}\label{sect--changing.base}
For facets $\bbf',\bbf\subset \bar\bba\subset \scrA$, we show that the category $\DTM(\Fl_\bbf)$ for the stratification in left-$\calP_{\bbf'}$-orbits is, to a certain extent, insensitive to the choice of the base scheme $S$, cf.~\thref{f!.nice} below.
In \refsect{Satake.category}, we will sharpen this idea by introducing the (abelian) Satake category $\Sat_G \subset \DTM(\Gr_G)$ and showing that this category is completely independent of the base scheme $S$.

Let $f\co T \r S$ be a map of schemes, where $T$ is Noetherian and of finite Krull dimension, so that $f^*\co \DM(S)\to \DM(T)$ is well-defined.
(An important example to have in mind is $T = \Spec \Fp \r S = \Spec \Z$.)
We indicate base changes to $T$ by a subscript, e.g., $G_T := G \x_S T$.
We still write $f$ for all maps obtained using such base changes, e.g., $f\co \Fl_{\bbf, T} \r \Fl_{\bbf, S}$.
The condition in \inter{Thm.~\iref{theo--motives.Ind-schemes}} is satisfied, so that we obtain a functor
$$f^*\co \DM(\Fl_{\bbf, S}) \r \DM(\Fl_{\bbf, T}).$$

As before, write $\iota_w\co \Fl_\bbf^w \r \Fl_\bbf$ for the inclusion of the $\calP_{\bbf'}$-orbits, both over $S$ and over $T$.
We clearly have an equivalence $(\iota_w)^* f^* \stackrel \simeq \r f^* (\iota_w)^*$ by functoriality, so that $f^*$ restricts to a functor
$$
f^*\co \DTM(\Fl_{\bbf,S})\to \DTM(\Fl_{\bbf,T}).
\eqlabel{base.change.DTM}$$

Here is the key lemma concerning the change of the base scheme.

\lemm\thlabel{base.change.nice}
Let $w \in W_{\bbf'} \bsl W / W_\bbf$. The following natural transformations of functors,  \emph{when restricted to the indicated categories of Tate motives},
$$\eqalign{
f^* (\iota_w)_* \r & (\iota_w)_* f^* \ \ \text{on }\;\, \DTM(\Fl^w_{\bbf, S})\cr
(\iota_w)_! f^* \r & f^* (\iota_w)_! \ \ \text {on }\;\, \DTM(\Fl^w_{\bbf, S}) \cr
f^* (\iota_w)^! \r & (\iota_w)^! f^* \ \ \text {on } \;\,\DTM(\Fl_{\bbf, S}) }
$$
are equivalences.
\xlemm

\pf
The claim for $(\iota_w)_!$ results from base change.
The claim for $(\iota_w)^!$ will follow from the others using an induction argument based on the localization fiber sequence
$$i^! \r i^* \r i^* j_* j^*$$
for any complementary closed (resp.~open) embedding $i$ (resp.~$j$).

In order to show that $f^*$ commutes with $(\iota_w)_*$, we may assume that $\bbf' = \bba$ (the base alcove), since the stratification by $\calP_{\bbf'}$-orbits is coarser than the one by Iwahori orbits so that the claim for the Iwahori stratification together with a localization argument implies the one for the stratification by $\calP_{\bbf'}$-orbits.

We first show the claim for $\bbf = \bba$.
By \inter{Prop.~\iref{prop--DTM.Fl.characterization}}, $\DTM(\Fl_{\bba})$ is the smallest cocomplete full subcategory of $\DM(\Fl_{\bba})$ which contains the twists of the unit motives supported at the base points $\{\tau\}$ for each $\tau\in \on{Stab}_{\bba}\subset W$ and which is stable under the operation $\pi_{s}^!\pi_{s,!}$ along the smooth proper projection maps $\pi_s\co \Fl_{\bba} \r \Fl_s := \Fl_{\bbf_s}$ for all simple reflections $s\in\bbS$. We proceed by induction on the length of $w$, the case $l(w)=0$ being trivial since $\iota_w$ is a closed embedding of a base point $\tau$ in this case. For $l(w)>0$, let $w=v\cdot s$ be a reduced expression with $s\in \bbS$. We obtain a fibre sequence
\begin{equation} \label{iota.vw}
(\iota_v)_* 1 \r \pi_s^! \pi_{s,!} (\iota_v)_* 1\r (\iota_w)_* 1,
\end{equation}
which is the dual of the fibre sequence \inter{(\iref{eqn--iota.v.w})}.
Using induction $l(v)<l(w)$, the functor $f^*$ commutes with $(\iota_v)_*$. Since $\pi_{s}$ is smooth, $f^*$ also commutes with $\pi_s^!$, and hence with $(\iota_w)_*$ by \eqref{iota.vw}. This finishes the case $\bbf = \bba$.

Now, for a general facet $\bbf \subset \bar{\bba}$, we reduce the claim to the one previously considered using the map $\pi\co \Fl_{\bba}\to\Fl_\bbf$.
This map is smooth, proper, surjective and a stratified map with respect to the Iwahori stratification on both ind-schemes \inter{Lem.~\iref{prop--change.facet}}.
Therefore, in the cartesian diagram
$$\xymatrix{
\bigsqcup_{v\in wW_\bbf} \Fl_{\bba}^v \ar[r] &  \pi^{-1}(\Fl_{\bbf}^w) \ar[r]^{\tilde \iota_w} \ar[d]^{\tilde \pi} & \Fl_{\bba} \ar[d]^\pi \\
& \Fl_\bbf^w \ar[r]^{\iota_w} & \Fl_\bbf
}$$
the preimage $\pi^{-1}(\Fl_\bbf^w)$ is itself stratified by some Iwahori strata on $\Fl_{\bba}$ as indicated above.
Using that $\pi$ is smooth, $(\tilde \iota_w)_* \tilde \pi^* = \pi^* (\iota_w)_*$.
Moreover, $\pi^*$ commutes with $f^*$.
Finally, $\pi^*$ is conservative.
Thus, to show that $f^*$ commutes with $(\iota_w)_*$, we may replace the inclusion $\iota_w$ (both on $T$ and on $S$) by $\tilde \iota_w$.
Using again a localization argument, we then reduce this statement to the one for the inclusions $\iota_v\co \Fl_{\bba}^v \r \Fl_{\bba}$ of the Iwahori strata in the flag variety refining the preimage stratification under $\pi$.
\xpf

Recall from \inter{Thm.~\iref{theo--motives.Ind-schemes}} that both categories carry weight structures.
The aim of this section is to prove:

\theo
\thlabel{f!.nice}
Let $f\co T\to S$ be a map of schemes both satisfying the conditions in \thref{base.scheme}.
Then the functor \refeq{base.change.DTM} has the following properties:
\begin{enumerate}
\item
\label{item--f!.conservative}
it is conservative.

\item
\label{item--f!.weights}
it creates weights, i.e., $M \in \DTM(\Fl_{\bbf, S})$ is of weights $\ge n$ \textup{(}resp. $\le n$\textup{)} iff $f^* M$ has the corresponding property.

\item
\label{item--f!.creates.t}
it creates the $t$-structure, i.e., $M \in \DTM(\Fl_{\bbf, S})$ is of in the ``$\ge n$'' \textup{(}resp. ``$\le n$''\textup{)} part of the $t$-structure iff $f^* M$ has the corresponding property.

\end{enumerate}
\xtheo


We need some preparation for the proof.

\prop
\thlabel{base.change.t.exact}
In the situation of \thref{f!.nice}, the functor \refeq{base.change.DTM} is $t$-exact with respect to the perverse motivic $t$-structures \textup{(}cf.~\thref{motivic.t.structure}\textup{)}, and commutes with the intermediate extension functors $(j_w)_{!*}$ defined in~\eqref{orbit.inclusion}. In particular,
$$f^* \big(\IC_{w, S}(n)\big) = \IC_{w, T}(n), \;\; n\in\bbZ, w\in W_{\bbf'}\bsl W/W_\bbf. \eqlabel{f*.IC}$$
\xprop

\pf
For both base schemes $S$ and $T$, the subcategory $\DTM(\Fl_{\bbf})^{\le 0}$ consists by definition precisely of those objects $M$ such that $\iota^* M \in \DTM(\Fl_\bbf^+)^{\le 0}$ where $\iota\co \Fl_{\bbf}^+=\sqcup_{w}\Fl^w_\bbf\to \Fl_\bbf$ denotes the disjoint union of the inclusions of all strata.
Likewise with ``$\ge 0$'' and $\iota^!$ instead.
To show the exactness of $f^*$ we may by \thref{base.change.nice} replace $f$ by the induced map $f^+\co \Fl_{\bbf, T}^+ \r \Fl_{\bbf, S}^+$.
It then remains to observe that the following diagram is cartesian and has smooth vertical maps
$$\xymatrix{
\Fl_{\bbf, T}^w \ar[d] \ar[r]^f & \Fl_{\bbf, S}^w \ar[d] \\
T \ar[r]^f & S.
}$$
Thus the $t$-exactness of $f^*$ on the base implies the one for $f^*\co \DTM(\Fl_{\bbf, S}^w) \r \DTM(\Fl_{\bbf, T}^w)$ since the hearts of these $t$-structures are generated by the objects $1(n)[d_w]$ where $d_w=l(w)$ is the dimension of $\Fl^w_\bbf$ relative to the base scheme (which is the same for $S$, resp.~$T$). The remaining claim now follows from Lemma \ref{base.change.nice} which ensure that $f^*$ commutes with all functors involved in the formation of $j_{w,!*} := \im (\pH^0 (j_{w,!}) \r {} \pH^0 (j_{w,*}))$.
\xpf

\begin{proof}[Proof of \thref{f!.nice}]
For \refit{f!.conservative},
let $M \in \DTM(\Fl_{\bbf, S})$. For the conservativity, we have to show $f^* M = 0$ implies $M = 0$.
Using the non-degeneracy of the perverse motivic t-structure \inter{Cor.~\iref{coro--t-structure}} and the t-exactness of $f^*$, it is enough to show the conservativity of $f^*|_{\MTM(\Fl_{\bbf,S})}$.
Any $M \in \MTM(\Fl_{\bbf, S})$ is the filtered colimit of its compact subobjects, so we may assume $M$ is also compact.
Then, $M$ has a Jordan-Hölder series with simple constituents given by twisted intersection motives \inter{\iref{theo--simple.objects}}.
We may thus assume that $M$ is an intersection motive, so we are done by \refeq{f*.IC}.
For \refit{f!.weights}, we need to show that $f^*$ is weight-exact and detects weights. As in the proof of \thref{base.change.t.exact}, to show that $f^*$ is weight-exact, we may replace $\Fl_{\bbf}$ by $\Fl_{\bbf}^+$ over both base schemes $S$ and $T$, which is again clear by definition of the weight structures.
The detection of weights then follows from \thref{weight.creation} below using part \refit{f!.conservative}.
For \refit{f!.creates.t}, we use likewise the conservativity and $t$-exactness of $f^*$.
\end{proof}

\lemm
\thlabel{weight.creation}
A conservative, weight-exact functor $F: \calC \r \calD$ between triangulated categories with weight structures detects weights: if $F(M)$ has weights $< n$ \textup{(}resp. $\ge n$\textup{)} for some $M \in \C$, then the same is true for $M$.
\xlemm

\pf
We use that $M$ has weights $\ge n$ (resp. $< n$) iff for any weight truncation triangle $E : M_{<n} \stackrel {s_{< n}} \r M \stackrel {s_{\ge n}} \r M_{\ge n}$, the maps $s_{\ge n}$ (resp. $s_{<n}$) are isomorphisms.
Indeed, the ``$\Leftarrow$'' direction holds by definition, the converse also follows from elementary applications of the axioms, see \cite[Cor.~2.2.6, Cor.~2.2.7]{Fontes:Weighty}.
Given a weight truncation triangle $E$ for $M$, $F(E)$ is a weight truncation triangle for $F(M)$ by assumption.
Hence our claim follows since $F$ is conservative.
\xpf

\subsubsection{Pullback functoriality for equivariant motives}
\lemm
The functor $f^*$ descends to a functor
$$\ol f^*\co \DM(\calP_{\bbf', S} \bsl LG_S / \calP_{\bbf, S}) \r \DM(\calP_{\bbf', T} \bsl LG_T / \calP_{\bbf, T}).$$
\nts{More formally, $\ol f^*$ commutes with $f^*$ via the forgetful functors, which are given by !-pullback along the projection map $LG / \calP_\bbf \r \calP_{\bbf'} \bsl LG / \calP_\bbf$.}
This functor $\ol f^*$ preserves the subcategories of equivariant Tate motives.
\xlemm

\pf
By construction, $f^*$ is the unique functor which is given by the usual $f^*$ on the level of finite type $S$-schemes and compatible with the insertion functors $\DM(\Fl_{\bbf}^{\le w}) \r \DM(\Fl_\bbf)$ (both over $S$ and $T$).
By \inter{Cor.~\iref{coro--Ind.Artin.co.limit}}, it is therefore enough to construct a functor
$$\ol f^*\co \DM(\calP_{\bbf', S} \bsl \Fl^{\le w}_{\bbf, S}) \r \DM(\calP_{\bbf', T} \bsl \Fl^{\le w}_{\bbf, T}).$$
There is a split pro-unipotent subgroup $U \subset \calP_{\bbf'}$ such that the quotient $K := \calP_{\bbf'} / U$ is smooth and of finite type and the $\calP_{\bbf'}$-action on $\Fl^{\le w}_\bbf$ factors over $K$.
By \inter{Prop.~\iref{prop--DM.G.homotopy.invariant}}, $\DM(\calP_{\bbf'} \bsl \Fl^{\le w}_{\bbf}) \cong \DM(K \bsl \Fl^{\le w}_\bbf)$.
Finally, in order to check the existence of $\ol f^*$ on this level, it is enough to observe that the maps in the bar construction $\BarC(K, \Fl^{\le w})$ are all smooth, and hence !-pullback along them commutes with $f^*$.
Hence the $f^*$-functors in all levels of the diagram $\DM^!(\BarC(K, \Fl^{\le w}_\bbf))$ glue to a functor on the limit of this diagram, which is $\DM (K \bsl \Fl^{\le w}_\bbf)$.

Given that the underlying non-equivariant functor of $\ol f^*$ is just $f^*$, the preservation of equivariant Tate motives follows from \refeq{base.change.DTM}.
\xpf

\subsection{Kazhdan-Lusztig parity vanishing}\label{sect--parity.vanishing}
We now apply \thref{base.change.t.exact} to prove the Kazhdan-Lusztig parity vanishing \cite[Thm.~5.5]{KazhdanLusztig:Poincare} (see also \cite[Thm.~11.c)]{Lusztig:Singularities}) for the intersection motives.
Our main tool is the $\ell$-adic realization functor which exists by assumption on $S$ (\thref{base.scheme}).
We continue with the notation and assumptions from \S\ref{sect--changing.base}.
In particular, we fix two facets $\bbf, \bbf'$ contained in the closure of the standard alcove $\bba$, and denote by $\DTM(\Fl_\bbf)$ the category of $(\bbf',\bbf)$-stratified Tate motives whose heart is the abelian category $\MTM(\Fl_\bbf)$ (\thref{motivic.t.structure}).

\theo\thlabel{realization.functor}
\textup{(}\inter{Thm.~\iref{theo--generators.DTM.flag}}\textup{)}
The restriction of the $\ell$-adic realization functor
\[
\rho_\ell\co \DM(\Fl_\bbf)\to \D_{\et}(\Fl_\bbf,\bbQ_\ell),
\]
to the subcategory $\DTM(\Fl_\bbf)$ is conservative.
Moreover, for $M \in \DTM(\Fl_\bbf)$ the following are equivalent: a\textup{)} $M$ lies in $\MTM(\Fl_\bbf)$, and b\textup{)} $\rho_\ell(M)$ is a perverse sheaf.
\xtheo

The following corollary is useful in lifting results from the $\ell$-adic to the motivic setting.

\coro
\thlabel{realization.functor.coro}
For each geometric point $f\co \bar s\to S$, the composition of functors
\[
f^*\circ \rho_\ell\co \MTM(\Fl_\bbf)^\comp\to \Perv(\Fl_{\bbf,\bar s},\Ql)
\]
is well-defined, exact, conservative and faithful.
\xcoro

\pf
Each object in $\MTM(\Fl_\bbf)^\comp$ admits a Jordan-Hölder series (\thref{motivic.t.structure} i)) whose simple constituents are the intersection motives $\IC_w(n)$ for $w\in W_{\bbf'}\bsl W/W_\bbf$ and $n \in \Z$.
Using the same method as in \thref{base.change.t.exact} we deduce that these are mapped under $\rho := f^* \circ \rho_\ell$ to the corresponding $\ell$-adic intersection complex on $\Fl_{\bbf,\bar s}^{\leq w}$.
Since the subcategory $\Perv(\Fl_{\bbf,\bar s},\Ql)\subset  \D_{\et}(\Fl_\bbf,\bbQ_\ell)$ is closed under extensions, it follows that $f^*\circ \rho_\ell$ is well-defined.

Being the restriction of an exact functor between triangulated categories, $\rho$ is exact as well.
For the conservativity of $\rho$ it is therefore enough to show that the simple objects, namely the $\IC_w(n)$ are not mapped to 0, which holds true by the above.

Being an exact conservative functor between abelian categories, $\rho$ is also faithful:
If a morphism $p\co A \r B$ maps to $0$ under $\rho$, then $\ker \rho(p) = \rho(\ker p)=A$ by exactness.
Hence, the natural map $\ker p \r A$ is mapped to an isomorphism, and therefore is an isomorphism by conservativity.
This shows $p=0$.
\xpf

We can now prove the Kazhdan-Lusztig parity vanishing for intersection motives.
Recall that for each class $w\in W_{\bbf'}\bsl W/W_\bbf$ the relative dimension of $\Fl^w\to S$ is given by $l(w)\in \bbZ_{\geq 0}$ where $l=l(\bbf',\bbf)$ denotes the length function as in \S\ref{sect--loop.grps.Satake}.
Let $e_w\co S=\{w\}\to \Fl_\bbf$ be the canonical inclusion.

\theo
\thlabel{KL.vanishing.theo}
For each $v,w \in W_{\bbf'}\bsl W/W_\bbf$, $n\in \bbZ$, one has
\begin{equation}\label{parity.eqn}
\H^i(e_v^*\IC_w(n))=0\;\;\;\;\text{whenever \;\; $i\not \equiv l(w)\mod 2$,}
\end{equation}
where $\H^i$ denotes the truncation with respect to the classical motivic $t$-structure \textup{(}cf.~\inter{Rem.~\iref{rema--classical.t-structure}}\textup{)} which agrees on $S$ also with the perverse motivic $t$-structure \textup{(}\thref{motivic.t.structure}\textup{)}.
\xtheo
\pf
We may assume $n=0$.
By \thref{realization.functor.coro}, it is enough to show that, for $S$ being the spectrum of a separably closed field, the $\ell$-adic intersection complex $\IC_{w,\ell}$ on $\Fl_\bbf^{\leq w}$ satisfies the parity vanishing \eqref{parity.eqn} where $\H^i$ denotes the classical cohomology functor.
This case is certainly well-known; we recall the part of the argument where we did not find a reference for the reader's convenience.\smallskip\\
{\it Reduction to the case $\bbf'=\bbf=\bba$.} By refining the orbit stratification on $\Fl_\bbf$ we may assume that $\bbf'=\bba$ is the standard alcove.
Now consider the projection $\pi\co \Fl_\bba\to \Fl_\bbf$ which is a smooth surjective map of relative dimension $d:=\dim(\calP_\bbf/\calP_\bba)$ by \inter{Prop.~\iref{prop--change.facet}}.
The preimage $\pi^{-1}(\Fl_\bbf^{\leq w})$ is a Schubert scheme in $\Fl_\bba$, and it follows from e.g.~\inter{Lem.~\iref{lemm--orbit.flag} iii)} that
$$\pi^{-1}(\Fl_\bbf^{\leq w})=\Fl_\bba^{\leq w_{\on{\max}}},$$
where $w_{\on{max}}$ is the unique representative of right maximal length with respect to $l_\bba:=l(\bba,\bba)$ in $w\cdot W_\bbf$.
Its length is $l_\bba(w_{\on{max}})=\dim(\Fl_\bba^{\leq w_{\on{\max}}})=l(w)+d$ by {\it loc.~cit.}.
As taking intermediate extensions commutes with smooth pullback, we have $\pi^*[d]\IC_{\ell,w}=\IC_{\ell,w_{\on{max}}}$.
Taking the cohomological shift into account and using the conservativity of pullback of surjective maps, we see that it is enough to prove \eqref{parity.eqn} in the case $\bbf'=\bbf=\bba$.\smallskip\\
{\it Proof for $\bbf'=\bbf=\bba$.} Here we refer to the classical sources \cite{KazhdanLusztig:Poincare}, \cite[A.7]{Gaitsgory:Central} and \cite{Haines:Purity}.

\xpf

\section{The convolution product}
\label{sect--convolution.product}

In this section, we will discuss the tensor structure on the category $\DM(\calP_\bbf \backslash LG / \calP_\bbf)$ given by convolution. We start with the definition and basic properties in \S\ref{sect--defi.associ}.
In \S\ref{sect--DTM.Tate}, we show that the convolution product preserves stratified Tate motives.

\subsection{Definition and associativity}\label{sect--defi.associ}

\defi
\thlabel{convolution.product}
Let $\bbf', \bbf, \bbf''$ be three facets in the closure of the standard alcove, see \S\ref{sect--loop.grps.Satake}. The \emph{convolution product} is the functor
$$
\star \colon \DM(\calP_{\bbf'} \bsl LG / \calP_\bbf) \x \DM(\calP_\bbf \bsl LG / \calP_{\bbf''}) \lr \DM(\calP_{\bbf'} \bsl LG / \calP_{\bbf''})
$$
defined by $(M_1, M_2) \mapsto M_1 \star M_2 := m_! p^! (M_1 \boxtimes M_2)$. Here the maps
$$\xymatrix{
\calP_{\bbf'} \backslash LG / \calP_\bbf \x \calP_\bbf \backslash LG / \calP_{\bbf''} &
\calP_{\bbf'} \backslash LG \x^{\calP_\bbf} LG / \calP_{\bbf''} \ar[l]_(.45)p \ar[r]^(.55)m &
\calP_{\bbf'} \backslash LG / \calP_{\bbf''}}$$
are the natural maps of prestacks induced by the identity on $LG \x LG$ (for $p$) and the multiplication $LG \x LG \r LG$ (for $m$).
\xdefi

\rema
\thlabel{explain.convolution}
i) The left adjoint $m_!$ of the functor $m^!$ or, equivalently the left adjoint of the $!$-pullback along the ind-proper map of ind-schemes $(LG \x^{\calP_\bbf} LG/\calP_{\bbf''})^\et \r (LG /\calP_{\bbf''})^\et = \Fl_{\bbf''}$ exists by \inter{Lem.~\iref{lemm--functoriality.equivariant}, Prop.~\iref{prop--f_!.ind-Artin}}.\smallskip\\
ii) The exterior product $M_1 \boxtimes M_2 \in \DM(\calP_{\bbf'} \bsl LG / \calP_\bbf \x \calP_\bbf \bsl LG / \calP_{\bbf''})$ exists by virtue of the construction in \inter{Prop.~\iref{prop--boxtimes}}, which gives a functor
$$\boxtimes\co \DM(\calP_{\bbf'} \bsl \Fl_\bbf) \t \DM(\calP_\bbf\bsl \Fl_{\bbf''}) \r \DM(\calP_{\bbf'} \bsl \Fl_\bbf \x \calP_\bbf\bsl \Fl_{\bbf''}),\eqlabel{boxy}$$
using the descent equivalence $\DM(\calP_{\bbf'} \bsl LG / \calP_\bbf) = \DM(\calP_{\bbf'} \bsl \Fl_\bbf)$.

\xrema

For clarity, we will momentarily denote the functor in \refeq{boxy} by $\boxtimes^\R$ and the convolution product stemming from this choice by $\star^\R$.
Alternatively, we may consider
$$\boxtimes^\L\co \DM(\Fl_{\bbf'}^\opp / \calP_\bbf) \t \DM(\Fl_\bbf^\opp / \calP_{\bbf''}) \r \DM(\Fl_{\bbf'}^\opp / \calP_\bbf \x\Fl_\bbf^\opp / \calP_{\bbf''})$$
where $\Fl_{\bbf'}^\opp=(\calP_{\bbf'}\bsl LG)^\et$ and likewise for $\Fl_\bbf^\opp$. The resulting convolution product functor is denoted by $\star^\L$.

\prop
\thlabel{convolution.HoDM.independent}
On the level of the homotopy categories, the two functors $\star^\R$ and $\star^\L$ are naturally isomorphic, i.e., one has $\star^\R \cong \star^\L$ as functors
$$ \Ho\big(\DM(\calP_{\bbf'} \bsl LG / \calP_\bbf)\big) \x \Ho\big(\DM(\calP_\bbf \bsl LG / \calP_{\bbf''})\big) \lr \Ho\big(\DM(\calP_{\bbf'} \bsl LG / \calP_{\bbf''})\big).$$
\xprop

\pf
Fix $v \le w$ in $W_{\bbf'}\bsl W / W_\bbf$ and consider the diagram
$$\xymatrix{
\Fl^{\opp, \le v}_{\bbf'} / \calP_\bbf \ar@{^{(}->}[d] &
\calP_{\bbf'} \bsl LG^{\le v} / \calP_\bbf \ar[l]_\approx \ar[r]^\approx \ar@{^{(}->}[d]^{i_{v, w}} &
\calP_{\bbf'} \bsl \Fl_\bbf^{\le v} \ar@{^{(}->}[d] \\
\Fl^{\opp, \le w}_{\bbf'} / \calP_\bbf \ar@{^{(}->}[d] &
\calP_{\bbf'} \bsl LG^{\le w} / \calP_\bbf \ar[l]_\approx \ar[r]^\approx \ar@{^{(}->}[d] &
\calP_{\bbf'} \bsl \Fl_\bbf^{\le w} \ar@{^{(}->}[d] \\
\Fl^{\opp}_{\bbf'} / \calP_\bbf &
\calP_{\bbf'} \bsl LG / \calP_\bbf \ar[l]_\approx \ar[r]^\approx &
\calP_{\bbf'} \bsl \Fl_\bbf,}$$
where $LG^{\leq w}=\calP_{\bbf'}w\calP_\bbf$ denotes the scheme-theoretic image of $\calP_{\bbf'}\x \calP_\bbf\to LG$, $(p,p')\mapsto p\cdot \dot{w}\cdot p'$ where $\dot w\in LG(S)$ is any representative of $w$, and likewise for $LG^{\leq v}$.
Note that this agrees with the preimage of $\Fl_\bbf^{\leq w}$ under the quotient map $LG\to \Fl_\bbf$, resp.~the preimage of $\Fl_{\bbf'}^{\opp, \leq w}$ under $LG\to \Fl_{\bbf'}^{\opp}$.
The labels $\approx$ at the horizontal arrows indicate maps of prestacks which become equivalences after étale sheafification and therefore descent equivalences upon applying $\DM$ (\thref{properties.DM}).

By \thref{DM!.placid.prestacks}, we have an exterior product, denoted by $\boxtimes$, for motives on placid prestacks such as the top middle term.
Under the descent equivalence, it is compatible with the exterior product $\boxtimes^\L$ for motives on prestacks of the form as in the top left term, and similarly with the top right term.
Of course, the same applies for the middle row as well.
Moreover, these identifications are compatible with the pushforwards along the maps (induced by the closed embeddings $i_{v,w}$) between the top and middle row, i.e., there is a natural equivalence
$$\alpha_{v, w}\co \boxtimes^w \circ (i_{v, w})_* \stackrel \cong \r (i_{v, w})_* \circ \boxtimes^v,$$
where $\boxtimes^w$ stands for an exterior product on terms as in the middle row of the diagram, and likewise for $\boxtimes^v$.
For yet another $u \le v$ in $W_{\bbf'}\bsl W / W_\bbf$, this equivalence and the one for $i_{u, v}$ and $i_{u, w}$ are compatible.

We obtain that the equivalence of \ii-categories
$$\DM\big(\Fl^\opp_{\bbf'} / \calP_\bbf\big) \cong \DM\big(\calP_{\bbf'} \bsl \Fl_\bbf\big)\eqlabel{MTM.rl}$$
is compatible with exterior products ($\boxtimes^\L$ and $\boxtimes^\R$, respectively), provided that we pass to the homotopy category and restrict to objects which are supported on some $\Fl^{\le w, \opp}$, resp.~$\Fl^{\le w}$. In particular, this is true for compact objects.
We may drop this compactness condition, since the homotopy category of a compactly generated category, such as the categories $\DM$ on the above prestacks, is again compactly generated by \cite[Rem.~1.4.4.3]{Lurie:HA}, and since the exterior product preserves filtered (homotopy) colimits separately in both variables.
\xpf

\rema
The point of passing to the homotopy categories $\Ho(\DM)$ is that these are \emph{ordinary} categories, as opposed to \ii-categories $\DM$.
For this reason, it is enough to check the compatibility of $\alpha_{v,w}$ for two composable maps, as opposed to verifying higher coherences.
We do not expect this loss of information to be necessary though: a more full-fledged approach would be to establish that $\DM^!$ is a symmetric lax monoidal functor on the \ii-category of ind-placid prestacks, such as $\calP_\bbf' \bsl LG / \calP_\bbf$.
\xrema

Hereafter, we will write $\star$ for $\star^\R$ above.
Since our main interest in this paper lies in the convolution product on the abelian (in particular ordinary) category $\MTM(L^+ G \bsl LG / L^+ G)$, \thref{convolution.HoDM.independent} shows that there is no ambiguity in the definition of the convolution product on this category.

As is well-known, the associativity of the convolution product is a consequence of the base-change formula:

\lemm
\thlabel{associative}
For $A \in \DM(\calP_{\bbf'} \backslash LG / \calP_\bbf)$, $B \in \DM(\calP_\bbf \backslash LG / \calP_\bbf)$ and $C \in \DM(\calP_\bbf \backslash LG / \calP_{\bbf''})$, there is a natural equivalence
$$(A \star B) \star C \cong A \star (B \star C).\eqlabel{star.associative}$$
\xlemm

\pf
For brevity, write $L := LG$ throughout the proof.
By construction in \inter{Prop.~\iref{prop--boxtimes}} (and the associativity of the exterior product for motives on schemes in $\Sch^\ft_S$), the exterior product for three motives on the three prestacks in the lower left entry of the diagram below admits an associativity isomorphism $(A \boxtimes B) \boxtimes C \cong A \boxtimes (B \boxtimes C)$.
Up to the exterior product in the definition of $(A \star B) \star C$, this convolution is computed as $m_! p^! (m \x \id)_! (p \x \id)^!$:
$$\xymatrix{
\calP_{\bbf'} \backslash L \x^{\calP_\bbf} L \x^{\calP_\bbf} L / \calP_{\bbf''} \ar[d]_{\id \x p} \ar[r]^{n := m \x \id} &
\calP_{\bbf'} \backslash L \x^{\calP_\bbf} L / \calP_{\bbf''} \ar[d]_p \ar[r]^m &
\calP_{\bbf'} \backslash L / \calP_{\bbf''} \\
\calP_{\bbf'} \backslash L \x^{\calP_\bbf} L / \calP_\bbf \x \calP_\bbf \backslash L / \calP_{\bbf''} \ar[d]_{p \x \id} \ar[r]^{m \x \id} &
\calP_{\bbf'} \backslash L / \calP_\bbf \x \calP_\bbf \backslash L / \calP_{\bbf''}  \\
\calP_{\bbf'} \backslash L / \calP_\bbf \x \calP_\bbf \backslash L / \calP_\bbf \x \calP_\bbf \backslash L / \calP_{\bbf''}.
}$$
The square in the above diagram is (homotopy) cartesian.
Moreover, the map $m$ is ind-proper, so that proper base change \inter{Prop.~\iref{prop--f_!.ind-Artin}} yields an equivalence
$$p^! (m \x \id)_! \stackrel \cong \lr n_! (\id \x p)^!.$$
Thus, the convolution $(A \star B) \star C$ can be computed by pullback and pushforward along the pictured composite correspondence.
Considering instead the composition of the correspondences computing $A \star (B \star C)$, we obtain the same composition, which yields a zig-zag of equivalences.
\xpf

\subsubsection{Reformulation in terms of schemes}\label{reformulation.convolution}
We now spell out the above definition in terms of ordinary schemes as opposed to prestacks. This relates to the classical definition of the convolution product as in \cite[§2]{Richarz:New}, and is used to show that the convolution product preserves Tate motives (\thref{convolution.Fl.Tate} below).

\defi
\thlabel{twisted.box.product.defi}
Let $\bbf', \bbf,\bbf''$ be facets as in Definition \ref{convolution.product}. We define
\[
\Fl_\bbf\xtw \Fl_{\bbf''}\defined (LG \x^{\calP_\bbf} LG/\calP_{\bbf''})^\et,
\]
which is an ind-proper $S$-ind-scheme. Consider the following commutative diagram of prestacks:
$$\xymatrix{
\calP_{\bbf'} \backslash LG / \calP_{\bbf''}
&
LG / \calP_{\bbf''} \ar[l] \ar[r]^\approx
&
\Fl_{\bbf''}
&
Z \ar[l]
\\
\calP_{\bbf'} \backslash LG \x^{\calP_{\bbf}} LG / \calP_{\bbf''} \ar[d]_p \ar[u]^m
&
LG \x^{\calP_\bbf} LG / \calP_{\bbf''} \ar[l]_u \ar[d] \ar[r]^\approx \ar[u]
&
\Fl_{\bbf} \xtw \Fl_{\bbf''} \ar[d]_{\tilde p} \ar[u]^{\tilde m}
&
X \xtw Y \ar[l]_{\tilde \iota} \ar[d]_{\tilde p} \ar[u]^{\tilde m}
\\
\calP_{\bbf'} \backslash LG / \calP_\bbf \x \calP_{\bbf} \backslash LG / \calP_{\bbf''}
&
LG / \calP_\bbf \x \calP_{\bbf} \backslash LG / \calP_{\bbf''} \ar[l] \ar[r]^\approx
&
\Fl_\bbf \x \calP_{\bbf} \backslash \Fl_{\bbf''}
&
X \x \calP_{\bbf} \backslash Y \ar[l]_\iota \ar[d]_e^\approx
\\
& & &
X \x \calP_{\bbf, i} \backslash Y
}
\eqlabel{convolution}
$$
The left-hand horizontal maps such as $u$ are the natural quotient maps.
By \inter{Lem.~\iref{lemm--DM.G.BarC}}, the !-pullback along such a map can be regarded as forgetting the $\calP_{\bbf'}$-action on some motive.
According to \inter{Thm.~\iref{theo--descent.prestacks}}, the horizontal maps labelled ``$\approx$'' induce descent type equivalences after applying $\DM^!$.
We will use similar equivalences without further comment; for example we identify motives on the double quotient $\calP_{\bbf'} \bsl LG / \calP_\bbf$ with those on $\calP_{\bbf'}\backslash\Fl_\bbf$.
The terms in the right hand column will be explained further below.

For $A\in \DM(\calP_{\bbf'}\backslash\Fl_\bbf)$, $B\in \DM(\calP_{\bbf}\backslash \Fl_{\bbf''})$, the \emph{twisted box product} is defined as
\[
A\ttw B \defined u^! p^! (A \boxtimes B) \in \DM(\Fl_\bbf\xtw\Fl_{\bbf''}).
\]
\xdefi

Let $\tilde{m}\co \Fl_\bbf\xtw \Fl_{\bbf''}\to \Fl_{\bbf''}$ be the map of ind-schemes induced by multiplication, i.e., the map $m$ above is the non-sheafified version obtained by passing to the left-$\calP_{\bbf'}$-quotients. By virtue of the following lemma, we will denote $\tilde{m}$ simply by $m$ (it will be clear from the context which version we mean).

\lemm
\thlabel{star.twiddle}
For $A\in \DM(\calP_{\bbf'}\bsl \Fl_\bbf)$, $B\in \DM(\calP_{\bbf}\bsl \Fl_{\bbf''})$, the object $\tilde{m}_!(A\ttw B)\in \DM(\Fl_{\bbf''})$ is the non-equivariant object underlying $A\star B=m_!(p^!(A\boxtimes B))$.
\xlemm

\pf
We only need to show that the functor $m_!$ commutes with the forgetful map to its non-equivariant version $\tilde{m}_!$.
This is precisely the characterization of $m_!$ in \inter{Lem.~\iref{lemm--functoriality.equivariant}}, see also \thref{explain.convolution} i) for its existence.
\xpf

Both functors, $\str \star \str $ and $\str \ttw \str$ preserve colimits separately in each variable.
They therefore factor over the Lurie tensor product $\DM(\calP_{\bbf'} \bsl LG / \calP_\bbf) \t \DM(\calP_\bbf \bsl LG / \calP_{\bbf''})$.
Since categories of motives are compactly generated \inter{Lem.~\iref{lemm--compact.motives.Ind.Artin}}, the functors are therefore determined by their values on compact objects.
Suppose then that $A$ and $B$ are compact objects, so they are supported on closed, finite type subschemes $X\subset \Fl_\bbf$, $Y\subset \Fl_{\bbf''}$ which are finite unions of Schubert schemes (these are the objects in the right vertical column in \refeq{convolution}).
The right-most vertical maps in the diagram, such as $\iota$ are induced by the closed embeddings of these subschemes.
In this case, $A\ttw B$ admits the following description:

Let $\calP_{\bbf}=\lim_{i\geq 0}\calP_{\bbf,i}$ as in \inter{Lem.~\iref{lemm--affine.proj}} and denote $\calU_{\bbf,i}:=\ker(\calP_\bbf\to \calP_{\bbf,i})$. We let $X_{i}\subset \Fl_{\bbf,i}:=(LG/\calU_{\bbf,i})^\et$ (resp.~$X_\infty\subset LG$) be the finite type (resp.~non-finite type) $S$-scheme defined by the preimage of $X$ under the canonical projection $\Fl_{\bbf,i}\to\Fl_\bbf$ (resp.~$LG\to \Fl_\bbf$). The left-$\calP_{\bbf}$-action on $Y$ factors through some $\calP_{\bbf,i}$ for $i>\!\!>0$.
We write $X \xtw Y := X_\infty \x^{\calP_\bbf} Y$, which is equivalent to $X_i \x^{\calP_{\bbf, i}} Y$ since $(X_\infty/\calP_\bbf)^{\et}=X$ and $(X_\infty/\calU_{\bbf,i})^{\et}=X_i$.
The vertical map $e$ labelled $\approx$ in the above diagram yields an equivalence upon applying $\DM^!$ (this stems from $\A^1$-invariance, using that $\calU_{\bbf, i}$ is split pro-unipotent, see \inter{Prop.~\iref{prop--DM.G.homotopy.invariant}}).
In particular, we can regard $B \in \DM(\calP_\bbf \backslash Y)$ as an object in $\DM(\calP_{\bbf, i} \backslash Y)$.
By the support setup, we can write $A \boxtimes B$ as $\iota_! (A_0 \boxtimes B_0)$ with $A_0 = \iota^* A$ etc., so that
$$A \ttw B = (e \circ \tilde p)^! (A \boxtimes B) = \tilde \iota_! ((e \circ \tilde p)^! A_0 \boxtimes B_0).$$
Note that the schemes $X$, $Y$, $\calP_{\bbf, i}$ and $Z$ intervening in the correspondence $X \x \calP_{\bbf, i} \backslash Y \stackrel{e \circ p} \gets X \xtw Y \stackrel{\tilde m} \r Z$, are of finite type over $S$ (unlike the remaining terms in the diagram).
We also see that the above definition of $\ttw$ agrees with the definition of $\ttw$ used for example in \cite[Lem.~2.20, Rmk.~2.21]{Richarz:New}.

Finally, writing $Z := \tilde m(X \xtw Y)$ (scheme-theoretic image, again a finite type $S$-scheme), $A \star B$ has as its underlying non-equviarant object $\tilde m_!( A \ttw B)$, which is, by proper base change, the !-pushforward along $Z \subset \Fl_\bbf$, of $\tilde m_! (e \circ \tilde p)^! (A_0 \boxtimes B_0)$.

\lemm
\thlabel{twiddle.basic}
i\textup{)} Let $X$,$Y$ be as above. Then $1_X\ttw 1_Y=1_{X \xtw Y}\in \DM(X \xtw Y)$.\smallskip\\
ii\textup{)} If $\xi: X \r X'$ \textup{(}resp.~$\upsilon: Y \r Y'$\textup{)} is an inclusion of finite type $\calP_\bbf'$-equivariant \textup{(}resp.~$\calP_\bbf$-equivariant\textup{)} subschemes of $\Fl_\bbf$ \textup{(}resp.~of $\Fl_\bbf''$\textup{)} then
$$(\xi \xtw \upsilon)_! (A \ttw B) = \xi_! A \ttw \upsilon_! B.$$
\xlemm
\pf The maps $X_i\x^{\calP_{\bbf,i}}Y\to X\times \calP_{\bbf,i}\bsl Y\gets X\x Y$ induce forgetful maps
\[
\DM(X_i\x^{\calP_{\bbf,i}}Y)\gets \DM(X\times \calP_{\bbf,i}\bsl Y)\to \DM(X\x Y),
\]
under which $1_X\boxtimes 1_Y=1_{X\x Y}$ corresponds to $1_{X \xtw Y}$ under $!$-pullback.
For the second statement note that the map $r := e \circ \tilde p$ in \refeq{convolution} is a $\calP_{\bbf, i}$-torsor, in particular a smooth map (of finite type $S$-schemes).
Therefore $r^!$ commutes with the exterior product and with the !-pushforward along the embeddings $X \x \calP_{\bbf, i} Y \r X' \x \calP_{\bbf, i} Y'$ and $\xi \xtw \upsilon: X \xtw Y \r X' \xtw Y'$.
\xpf


\subsubsection{Compatibility with the $\ell$-adic realization}\label{compatibility.realization.convolution}
We denote by $\D^\bd_{\cstr, \calP_\bbf}(\Fl_{\bbf''}, \Ql)$ the category of $\calP_\bbf$-equivariant $\ell$-adic sheaves on $\Fl_{\bbf''}$.
In \cite{MirkovicVilonen:Geometric} (see also \cite[§2]{Richarz:New} and \cite[\S10.2]{PappasZhu:Kottwitz}), the convolution product for $\ell$-adic sheaves
$$\star_\ell\co \D^\bd_\cstr(\Fl_\bbf, \Ql) \x \D^\bd_{\cstr, \calP_\bbf}(\Fl_{\bbf''}, \Ql) \r \D^\bd_\cstr(\Fl_{\bbf''}, \Ql)$$
is defined as follows: consider the diagram
$$\Fl_\bbf \x \Fl_{\bbf''} \stackrel p \longleftarrow LG \x \Fl_{\bbf''} \stackrel q \lr \Fl_\bbf \xtw \Fl_{\bbf''} = (LG \x^{\calP_\bbf} LG/\calP_{\bbf''})^\et \stackrel m \lr \Fl_{\bbf''}.\eqlabel{convolution.terms}$$
Then $A_1 \star_\ell A_2 := m_* (A_1 \ttw_\ell A_2)$, where $A_1 \ttw_\ell A_2$ is the unique object in $\D^\bd_\cstr(\Fl_\bbf \xtw \Fl_{\bbf''},\Ql)$ such that $p^* (A_1 \boxtimes_\ell A_2) = q^* (A_1 \ttw_\ell A_2)$, using the $\calP_\bbf$-equivariance of $A_2$.

\prop
\thlabel{convolution.DM.ell}
Under the $\ell$-adic realization functor $\rho_\ell$ (cf.~\inter{Thm.~\iref{theo--f!.Artin}\iref{item--realization.calX}})
the convolution product corresponds to the convolution product $\star_\ell$ considered in the context of the $\ell$-adic Satake equivalence, i.e., there is a natural isomorphism
$$\rho_\ell (A \star B) \cong \rho_\ell(A) \star_\ell \rho_\ell(B)$$
for $A \in \DM(LG / \calP_{\bbf})^\comp$ and $B \in \DM\big(\calP_\bbf \bsl LG / \calP_{\bbf''}\big)^\comp$.
\xprop

\pf
The twisted box product $\ttw$ for motives (\thref{twisted.box.product.defi}) is formed using descent along $!$-pullbacks.
We need to compare its $\ell$-adic realization with $\ttw_\ell$, formed using $*$-pullbacks.

Since the objects $A, B$ are compact, hence supported on $S$-finite type closed subschemes $X\subset \Fl_\bbf$, $Y\subset \Fl_{\bbf''}$, we can replace the maps $p$, $q$ in \refeq{convolution.terms} by the diagram of $\calP_{\bbf,i}$-torsors
\[
X \x Y\stackrel p \longleftarrow X_i\x Y \stackrel q \lr X \xtw Y,
\]
for some $i\geq 0$ as in \S\ref{reformulation.convolution} above.
Let $G:=\calP_{\bbf,i}$ which is a smooth affine $S$-group scheme acting on $Z:=X_i\x Y$ either via the torsor $p$ or $q$.
By definition, the category of $G$-equivariant $\ell$-adic sheaves on $Z$ is defined as
$$\eqalign{
\D^\bd_{\cstr, G}(Z, \Ql) & := \lim \left ((\D^\bd_\cstr)^*(\BarC(Z, G),\Ql) \right) \cr
& := \lim \left ( \D^\bd_\cstr(Z,\Ql) \stackrel[p^*]{a^*} \rightrightarrows \D^\bd_\cstr(G \x Z,\Ql)\substack{\rightarrow\\[-1em] \rightarrow \\[-1em] \rightarrow} \cdots \right ),}$$
where we emphasize that the functors in this limit are the $*$-pullbacks along the maps in the bar complex.
The motivic analogue of that category is $\DM^*(G \bsl Z)^\comp := \lim \DM^*(\BarC(Z, G))^\comp$, where again we use $*$-pullbacks to form the limit.
(See also \inter{Rem.~\iref{rema--DM.prestacks}, iv)} for further discussion of the presheaf $\DM^*$.)

The vertices $(G)^{\x n} \x Z$ of the bar construction are separated $S$-schemes of finite type, and the action and projection maps in this diagram are smooth and affine, noting that $G\r S$ is so.
We can therefore use the equivalence of $\DM^*$ with $\DM^!$ applied to the smooth morphisms $p$, $q$, see \thref{DM*vs!}. 
Under this equivalence the functor $p^!$ (resp.~$q^!$) corresponds to $p^*$ (resp.~$q^*$).
Since $p$, $q$ have the same relative dimension $\dim(G/S)$, we can equivalently form $A\ttw B$ using descent along $*$-pullbacks. 
Moreover, the map $m$ is ind-proper, so that $m_* = m_!$.
We conclude using that $\rho_\ell$ is compatible with the six functors.
\xpf

\subsubsection{Convolution product and change of base scheme}

\lemm
\thlabel{f*.convolution}
Let $f\co T \r S$ be a map of schemes satisfying the assumptions in \thref{base.scheme}.
Then, for $M_1 \in \DM(\calP_{\bbf'} \bsl LG / \calP_\bbf)$ and $M_2 \in \DM(\calP_{\bbf'} \bsl LG / \calP_\bbf)$, there is a natural isomorphism
$$f^* (M_1 \star M_2) \cong (f^* M_1) \star (f^* M_2).$$
\xlemm

\pf
All functors involved in the definition of $\star$ are compatible with $f^*$.
For $p^!$, this holds true since (the étale sheafification of) $p$ is a $\calP_\bbf$-torsor, in particular pro-smooth.
\xpf

\subsection{Preservation of Tate motives}\label{sect--DTM.Tate}

In this section, we show that the convolution product on partial affine flag varieties respects stratified Tate motives.
A key point in this proof is the well-knwon distinguished triangle \refeq{convolution.Fl.Tate:2} below (cf.~\cite[App.]{KazhdanLusztig:Hecke}), which is a geometric incarnation of the following formula for the multiplication in the Iwahori-Hecke algebra over a finite field $\Fq$ given in \cite[IV, \S2, Ex.~24]{Bourbaki:Lie456}:
$$\phi_s \star \phi_s = (q-1)\cdot \phi_s + q\cdot \phi_e,\eqlabel{Iwahori}$$
where $\phi_s$ is the characteristic function of the Iwahori double coset $\calB(\Fq) s \calB(\Fq)$ for a simple reflection $s$, and $\phi_e$ is the characteristic function of the base point. Here $\calB:=\calP_\bba$ denotes the standard Iwahori subgroup associated with the choice of the alcove $\bba$.

\theo
\thlabel{convolution.Fl.Tate}
Let the base scheme $S$ be as in \thref{base.scheme}.
For any three facets $\bbf',\bbf,\bbf''$ contained in the closure of $\bba$, the convolution product restricts to a functor
$$\star\colon \DTM (\calP_{\bbf'} \backslash LG / \calP_\bbf) \x \DTM (\calP_\bbf \backslash LG / \calP_{\bbf''}) \r \DTM(\calP_{\bbf'} \backslash LG / \calP_{\bbf''}).$$
In particular, taking $\bbf'=\bbf=\bbf'' = 0$, there is a convolution product on $\DTM(L^+ G \bsl LG / L^+G)$.
\xtheo
Let $A \in \DTM (\calP_{\bbf'} \backslash LG / \calP_\bbf)$, $B\in \DTM (\calP_\bbf \backslash LG / \calP_{\bbf''})$ and denote $\Fl_{\bbf''}:=(LG/\calP_{\bbf''})^\et$.
By \inter{Def.~\iref{defi--DTM.G}}, we have to show that the non-equivariant motive underlying $A\star B$ in $\DM(\Fl_{\bbf''})$ is Tate.
Further, we may assume that both objects $A, B$ are compact.
By \inter{Thm.~\iref{theo--equivariant.DTM.flag}}, the motives $A$ and $B$ are constructed in finitely many steps from the generators
\[
\iota_{w,!}1 \in \DTM (\calP_{\bbf'} \backslash LG / \calP_\bbf),
\]
where $w \in W_{\bbf'} \bsl W / W_{\bbf}$ and $\iota_w\co \calP_{\bbf'} \backslash LG^w / \calP_\bbf \to \calP_{\bbf'} \backslash LG / \calP_\bbf$ denotes the inclusion of the stratum\footnote{Formally, $LG^w$ is the scheme-theoretic image of the map $\calP_{\bbf'}\x\calP_\bbf\to LG$, $(b,p)\mapsto b\cdot \dot w\cdot p$ which is well defined because the source is quasi-compact and the target is an ind-scheme.} $LG^w:=\calP_{\bbf'}\dot w\calP_\bbf$ where $\dot w\in LG(S)$ is any representative of $w$.
We have to show that the non-equivariant motive underlying $\iota_{w,!}1 \star \iota_{w',!}1$ is Tate for any $w\in W_{\bbf'}\bsl W/W_\bbf$, $w'\in W_{\bbf}\bsl W/W_{\bbf''}$, i.e., by \thref{star.twiddle}, that
\begin{equation}\label{really.to.show}
\tilde m_! (\iota_{w,!}1 \ttw \iota_{w',!}1) \in \DTM(\Fl_{\bbf''}).
\end{equation}
We show \eqref{really.to.show} in several steps starting with the following key case.

\prop
\thlabel{convolution.Fl.Tate.Iwahori}
If $\calP_{\bbf'}=\calP_\bbf=\calP_{\bbf''}=:\calB$ is the standard Iwahori subgroup, then \thref{convolution.Fl.Tate} holds.
\xprop
\pf
In this case, $W_{\bbf'}=W_\bbf=W_{\bbf''}=1$, so that $w,w'\in W$.
\smallskip\\
{\it First case.} Assume that $w=w'=s$ is a simple reflection. Then there is an isomorphism
\[
\tau=(p,\tilde m)\co \Fl^{\leq s}\xtw \Fl^{\leq s} \overset{\simeq}{\lr} \Fl^{\leq s}\x \Fl^{\leq s},
\eqlabel{iso.Fls.proof}
\]
where $p\co (x,y)\mapsto x$ is the projection on the first factor, and $\tilde m\co (x,y)\mapsto x\cdot y$ is the multiplication map, as above.
Note that the image of $\tilde m$ lies inside $\Fl^{\le s}$ by standard properties of Tits systems \cite[Ch.~IV, \S2.1, (2)]{Bourbaki:Lie456} which are applicable in view of \cite[\nopp 5.2.12]{BruhatTits:Groups2} (see also the discussion in \cite[\S 1.1]{Richarz:Iwahori}). Further, $\tau$ being a closed immersion (being a proper monomorphism) between integral schemes of the same dimension, it must be an isomorphism.
Under this isomorphism, we have for the strata $\tau(\{e\}\xtw\{e\})=\{e\}\x\{e\}$, $\tau(\{e\}\xtw\Fl^s)=\{e\}\x\Fl^s$ and $\tau(\Fl^s\xtw \{e\})=\Delta(\Fl^s)$, where $\Delta\co \Fl^{\leq s}\r \Fl^{\leq s}\x\Fl^{\leq s}$ denotes the diagonal. We conclude that
\[
\tau(\Fl^s\xtw\Fl^s) \,=\, (\Fl^s \x \Fl^{\leq s})\backslash \Delta(\Fl^s).
\]
Identifying $\Fl^{\leq s}\simeq \bbP^1$ such that $\{e\}\simeq \{\infty\}$ and $\Fl^s\simeq \bbA^1$, we get a commutative diagram of $S$-schemes
$$
\xymatrix{
\Fl^{\le s}
&
\im (\tilde m) \ar@{=}[l] \ar[r]^\cong
&
\P^1
\\
\Fl^{\le s} \xtw \Fl^{\le s} \ar[u]^{\tilde m} 
&
\Fl^{s} \xtw \Fl^{s} \ar[l]_(.4){\iota_{s \xtw s}} \ar[r]^(.4)\cong \ar[u]^{\tilde m} 
&
(\A^1 \x \P^1) \setminus \Delta(\A^1) \ar[u]^{q}
}
$$
where $\iota_{s\xtw s}$ and $\iota_s$ are the inclusion of the open strata and $q$ is the projection onto the second factor.
Writing $a := \tilde m \circ \iota_{s \xtw s}$ and using \thref{twiddle.basic} for $\iota_{s \xtw s}=\iota_s \xtw \iota_s$, we have to prove
\[
\tilde m_! (\iota_{s,!}1 \ttw \iota_{s,!}1)\,=\,
\tilde m_! (\iota_s \xtw \iota_s)_! 1 \,=\,
a_!1 \in \DTM(\Fl^{\le s}),
\]
or equivalently that
$$M := q_! 1 \in \DTM(\P^1),$$
where the Tateness of the motive is with respect to the stratification of $\P^1\simeq \Fl^{\leq s}$ by $\{\infty\} \sqcup \A^1\simeq \Fl^e \sqcup \Fl^s$.
To check this, we use the localization sequence, noting that $q^{-1}(\{\infty\})\simeq \bbA^1$ and $q^{-1}(\bbA^1)\simeq \bbA^1\x \Gm$ which gives
$$
(\iota_{\A^1})_! \iota_{\A^1}^* M = \iota_{\A^1,!} (1(-1)[-2] \oplus 1[-1]) \lr M \lr \iota_{\infty,!}i_{\infty}^* M = \iota_{\infty_!} 1(-1)[-2].
\eqlabel{convolution.Fl.Tate:2}
$$
The outer terms are in $\DTM(\P^1)$, hence so is the middle.
This finishes the first case.\smallskip\\
{\it Second case.} Let $w,w' \in W$, and assume $l(w w')=l(w)+l(w')$ for the Bruhat length. Then the composition of the open inclusion $\iota_{w \xtw w'}\co \Fl^{w}\xtw \Fl^{w'}\r \Fl^{\leq w}\xtw \Fl^{\leq w'}$ with the multiplication map $m$ is an isomorphism
\[
m\circ j_{w\xtw w'}\co \Fl^{w}\xtw \Fl^{w'}\overset{\simeq }{\lr} \Fl^{ww'}.
\eqlabel{iso.Flw.proof}
\]
This implies that $\iota_{w,!}1\star \iota_{w',!}1\simeq \iota_{ww',!}1\in \DTM(\Fl)$ is Tate and finishes the second case.\smallskip\\
{\it Third case.} Let $w,w'\in W$, and assume that $w'=s$ is a simple reflection. If $l(ws)=l(w)+1$, then we conclude $\iota_{w,!}1\star \iota_{s,!}1\in \DTM(\Fl)$ by the second case. If $l(ws)=l(w)-1$, we write $v:= ws$. Since $s^2=1$, we have $w=vs$, and by construction $l(vs)=l(v)+1$. In particular, $\iota_{w,!}1\simeq \iota_{v,!}1\star \iota_{s,!}1$ by the second step. Applying $\iota_{v,!}1\star (\str)$ to the localization sequence \refeq{convolution.Fl.Tate:2}, we get a cofiber sequence
\begin{equation}\label{convolution.Fl.Tate:3}
(\iota_{w,!}1(-1)[-2]\oplus \iota_{w,!}1[-1])\lr \iota_{v,!}1\star (\iota_{s,!}1\star \iota_{s,!}1)\lr \iota_{v,!}1(-1)[-2]\lr
\end{equation}
where $\iota_{v,!}1\star (\iota_{s,!}1\star \iota_{s,!}1)\simeq (\iota_{v,!}1\star \iota_{s,!}1)\star \iota_{s,!}1\simeq \iota_{w,!}1\star \iota_{s,!}1$ by \thref{associative}. Hence, $\iota_{w,!}1\star \iota_{s,!}1$ is an extension of Tate motives on $\Fl$, and thus defines an object in $\DTM(\Fl)$. This finishes the third case. \smallskip\\
{\it General case.} Let $w,w'\in W$ be arbitrary. Fix a reduced expression $w'=s_1\cdot\ldots\cdot s_n$ where $s_i$ are simple reflections and $n=l(w')$. By the second case, we have $\iota_{w',!}1\simeq \iota_{s_1,!}1\star \ldots \star \iota_{s_n,!}1$ where we omit the parenthesis in view of \thref{associative}. By repeated use of the third case, we conclude that $\iota_{w,!}1\star \iota_{w',!}1\in \DTM(\Fl)$ is Tate. This finishes the general case, and the theorem follows.
\xpf

\rema
\thlabel{motivation.formula}
If $k = \Fq$, \refeq{convolution.Fl.Tate:2} gets mapped by the $\ell$-adic realization to
$$\iota_{s,!} \Ql(-1)[-2] \oplus \iota_{s,!} \Ql[-1] \r \iota_{s,!} \Ql \star \iota_{s,!} \Ql \r \iota_{e,!} \Ql(-1)[-2].$$
Taking the alternating trace of the geometric Frobenius, we obtain the identity \refeq{Iwahori} using the relation $\on{trace} (\on{Frob} | \Ql(-1)) = q$.
\xrema

\rema The method used in the proof of \thref{convolution.Fl.Tate} works more generally for not necessarily split reductive groups $G$ defined over $k\rpot{\varpi}$ which are residually split, i.e., $\Fl^{s}\simeq \bbA^1_k$ whenever $s\in W$ is a simple reflection. However, we will not need these non-split cases in this manuscript.
\xrema

\prop
\thlabel{convolution.Fl.Tate.2x.Iwahori}
\thref{convolution.Fl.Tate} holds in the case that $\calP_{\bbf'}=\calP_{\bbf''}=\calB$ is the standard Iwahori subgroup of $LG$, and any facet $\bbf$ in the closure of the standard alcove \textup{(}so that $\calB$ is a subgroup of $\calP_\bbf$\textup{)}.
\xprop

\pf
Let $\calP := \calP_\bbf$.
For any $w\in W/W_\bbf$, $w'\in W_{\bbf}\bsl W$, consider the diagram
$$\xymatrix{
L^w \x^\calB L^{w'}/\calB \ar[d]^{p_\calB} \ar@/^2pc/[rr]^{m_\calB} \ar[r]^b
&
L^w \x^{\calP} L^{w'}/\calB \ar[d]^{p} \ar[r]^{m} \ar@{..>}@/^1pc/[l]^s & L/\calB
\\
L^w / \calB \x \calB \backslash L^{w'} / \calB \ar[r]^a
&
L^w/\calP \x \calP \backslash L^{w'} / \calB,
}$$
where $L^w:=LG^w=\calB w\calP\subset LG$ (resp.~$L^{w'}:=LG^{w'}=\calP w'\calB\subset LG$). To construct the map $s$ with $b\circ s=\id$, it suffices to construct a section to the composition of quotient maps
\[
L^w\times L^{w'}\to L^w\times^\calB L^{w'} \to L^w\times^\calP L^{w'},
\]
which is equivariant for right $\calB$-action on the second factor.
It follows from \inter{Prop.~\iref{prop--cells.flag}} that there exists a closed subscheme $U\subset \calB$, a finite direct product of some affine root groups (depending on $w$), such that the multiplication $U\dot{w}\x \calP\to L^w$ is an isomorphism.
Here $\dot{w}\in W$ is any representative of $w\in W/W_\bbf$. Thus, $L^w\x^\calP L^{w'}\simeq U\dot{w}\x L^{w'}\subset L^w\x L^{w'}$ is the desired section.

In particular, the adjunction map $b_!b^!\to \id$ has the section $\id\simeq b_!s_!s^!b^!\to b_!b^!$, so that the motive $M:=p^!(\iota_{w,!}1\boxtimes \iota_{w',!}1) \in \DM(L^w \x^{\calP} L^{w'}/\calB)$ is a direct summand of $b_!b^!M$. Thus, the motive $m_!M=\iota_{w,!}1\star \iota_{w',!}1$ is a direct summand of the motive
\[
m_!b_!b^!M\simeq (m\circ b)_! (p\circ b)^!(\iota_{w,!}1\boxtimes \iota_{w',!}1) \simeq m_{\calB,!}p_\calB^! \big(a^!(\iota_{w,!}1\boxtimes \iota_{w',!}1)\big).
\]
The map $\Fl^w_\bba \x \Fl^{w'}_\bba \r \Fl^w_\bbf \x \Fl^{w'}_\bba$ is stratified with respect to the stratification by $\calB \x \calP$-orbits.
Hence $a^!(\iota_{w,!}1\boxtimes \iota_{w',!}1)$ is stratified Tate.
Therefore the motive $m_!b_!b^!M$ is stratified Tate by \thref{convolution.Fl.Tate.Iwahori}, and so is $\iota_{w,!}1\star \iota_{w',!}1$ as a direct summand. This proves \eqref{really.to.show}.
\xpf

\begin{proof}[We now prove \thref{convolution.Fl.Tate} in the general case.]
By definition, the category $\DTM(\calP_{\bbf'} \backslash LG / \calP_{\bbf''})$ consists of those $\calP_{\bbf'}$-equivariant motives on $\Fl_{\bbf''} = (LG / \calP_{\bbf''})^\et$ whose underlying non-equivariant motive is a Tate motive with respect to the stratification by $\calP_{\bbf'}$-orbits.
Such an orbit is the form $X = (\calP_{\bbf',i} / H)^\et$, where $\calP_{\bbf'}\to \calP_{\bbf',i}$ is a smooth $S$-affine quotient (i.e., of finite type) and $H\subset \calP_{\bbf',i}$ is a smooth closed subgroup scheme with connected fibers over $S$.  Let $e\co S \r X$ be the unit section.
The composition
\[
\DM_{\calP_{\bbf'}}(X) \simeq \DM_{\calP_{\bbf',i}}(X) \stackrel{e^!} \r \DM_H(S),
\]
which is an equivalence by \inter{Lem.~\iref{lemm--DM.G/H}}, restricts to an equivalence $e^!\co \DTM_{\calP_{\bbf'}}(X) \simeq \DTM_H(S)$ on Tate motives by \inter{Prop.~\iref{prop--DTM.G/H}}.
Thus, for a $\calP_{\bbf'}$-equivariant motive $M$ on $\Fl_{\bbf''}$, the following are equivalent:
\begin{itemize}
\item $M$ is Tate with respect to the stratification by $\calP_{\bbf'}$-orbits;
\item $M$ is Tate at the base point of each $\calP_{\bbf'}$-orbit;
\item $M$ is Tate with respect to the (finer) stratification by $\calB$-orbits.
\end{itemize}
We may therefore assume that $\calP_{\bbf'}=\calB$ is the standard Iwahori.
Using that $\DTM(\calB \bsl LG / \calP_{\bbf''})$ consists of $\calP_{\bbf''}$-equivariant motives whose underlying motive on $\Fl^\opp = (\calB \bsl LG)^\et$ is stratified Tate \inter{\iref{theo--equivariant.DTM.flag}}, we similarly reduce to the case that $\calP_{\bbf''}=\calB$ and therefore deduce \thref{convolution.Fl.Tate} in general from \thref{convolution.Fl.Tate.2x.Iwahori}.
\end{proof}

\section{Purity of Tate motives}
\label{sect--purity}

Throughout \refsect{purity}, the base scheme $S$ is as in \thref{base.scheme}.


\subsection{Intersection motives are pure}

In this section, we show that the intersection motives $\IC_w$ for the stratification of $\Fl_\bbf$ given by the $\calP_{\bbf'}$-orbits (for arbitrary facets $\bbf, \bbf'$ contained in the closure of the base alcove $\bba_0$) are pure.
This will be proven by lifting the corresponding fact for $\ell$-adic intersection complexes to motives over $\Fp$, which is then extended to more general base schemes using the results of \refsect{changing.base}.

\theo
\thlabel{IC.pure}
The intersection motives $\IC_w \in \MTM(\Fl_\bbf)$ of the $\calP_{\bbf'}$-orbits are pure of weight $\dim \Fl_\bbf^w$.
(Here purity refers to the weight structure for motives on ind-schemes, established in \cite[Theorem~2.4.2]{RicharzScholbach:Intersection}.)
\xtheo

\pf
Let $S = \Spec \Fp$ first. Pick some prime $\ell$ invertible in $\Fp$.
The $\ell$-adic realization $\rho_\ell: \DM(\Fl_\bbf)^\comp \r \D^\bound_\cstr(\Fl_\bbf, \Ql)$ takes values in the subcategory $\D_{\et, \mix}(\Fl_\bbf, \Ql)$ of mixed complexes.
With respect to the standard weight structure on that category (and the motivic weight structure on $\DM$), the functor is weight-exact.
(This follows from the definition of the motivic weight structure and the fact that in the realization these functors preserve weights, see \cite[Prop.~3.6.1.2]{Bondarko:WeightsForRelative} for details.)
Its restriction to $\DTM(\Fl_\bbf)$ is t-exact and conservative by \inter{Lem.~\iref{lemm--Tate.conservative}}.
It therefore \emph{creates} the $t$-structure and the weight structure (\thref{weight.creation}).
Now recall the notation \eqref{orbit.inclusion}.
Since $\rho_\ell$ also commutes with $(\iota_w)_!$ and $(\iota_w)_*$, the motive $\IC_w$ is mapped under $\rho_\ell$ to the $\ell$-adic intersection complex on the Schubert variety $\Fl_\bbf^{\leq w}$, i.e., to
$\rho_\ell(\IC_w) = (i_{w})_*(j_w)_{!*} (\Ql[\dim(\Fl^w_\bbf)])$.
By \cite[Ex.~III.10.3]{KiehlWeissauer:Weil}, it is pure of weight $+\dim \Fl^w_\bbf$, hence so is $\IC_w$ itself.

For general $S$, consider the zig-zag $S \r \Spec \Z \gets \Spec \Fp$.
By \thref{f!.nice} and \refeq{f*.IC}, the purity of $\IC_{w, S} \in \MTM(\Fl_{\bbf, S})$ is equivalent to the one of $\IC_{w, \Z} \in \MTM(\Fl_{\bbf, \Spec \Z})$, which in turn is equivalent to the one of $\IC_{w, \Fp} \in \MTM(\Fl_{\bbf, \Spec \Fp})$.
\xpf

\coro
\thlabel{weight.grading}
There is a functorial weight filtration for any object $M \in \MTM (\Fl_\bbf)$
$$0 = M_{-\infty} \subset \dots \subset M_0 \subset M_1 \subset \dots M_\infty = M$$
such that $\gr^\W_i(M) := M_i / M_{i-1} \in \MTM(\Fl_\bbf)$ is pure of weight $i$ (for $-\infty < i< \infty$).
If $M$ is compact, this filtration is finite, and the $\gr^\W_i(M)$ are compact.
\xcoro

\pf
Let $\MTM_i := \{ \IC_w(n), \dim \Fl_\bbf^w = i + 2n \}$.
The collection of these subcategories forms a semi-orthogonal family in the sense of \cite[Def.~1.1.4]{Bondarko:WeightHearts}: $\Hom_{\DTM(\Gr)}(\MTM_i, \MTM_j[s])$ vanishes for $s < 0$ by the existence of the perverse motivic t-structure.
It vanishes for $s > i-j$ for weight reasons, noting that $\MTM_j[s]$ consists of objects which are pure of weight $j+s$ by \thref{IC.pure}.
The claim for compact objects then follows from \cite[Thm.~1.2.1, (ii'), (iv')]{Bondarko:WeightHearts}.
In general, any object $M \in \MTM$ is the filtered colimit of its compact subobjects $M^{(n)} \subset M$, so that we obtain a (possibly infinite) weight filtration of $M$ by taking the colimit of the weight filtrations of the $M^{(n)}$.
\xpf

\subsection{Convolution preserves weights}

In this section, we show that the convolution product preseves the subcategories of motives of weight $\le n$ and $\ge n$.
To prove this, we need the following lemma:

\lemm\thlabel{exterior.product.weights}
For separated schemes $X_1, X_2$ of finite type over a perfect field $k$, the exterior product
$$\boxtimes : \DM(X_1) \x \DM(X_2) \r \DM(X_1 \x_k X_2)$$
is weight-exact.
\xlemm

\pf
By the description of the $\le 0$- and $\ge 0$-part of the weight structure in \cite[Rem.~1.17]{Hebert:Structure}, we have to show that exterior products of motives of the form $(f_i)_! 1(n)[2n]$, where $f_i: T_i \r X_i$ is a proper map and $T_i$ is regular for $i=1,2$, is again pure of weight 0.

Since $k$ is perfect, a finite type $k$-scheme is regular iff it is smooth over $k$ \StP{0B8X}.
Hence $T_1 \x_k T_2$ is again regular.
We conclude using the formula $f_{1, !} 1 \boxtimes f_{2, !} 1 = (f_1 \x f_2)_! 1$ whose proof is straightforward using the base change formula for $f^*$ vs. $g_!$ and the projection formula (see \inter{Synopsis \iref{syno--motives}, vii) and x)}).
\xpf

\rema
If $k$ is an imperfect field, \thref{exterior.product.weights} still holds by virtue of the equivalence
$$\DM(X_i) \cong \DM(X_i \x_k \hat k)$$
for some perfect closure $\hat k$ of $k$.
The stronger corresponding statement for $\SH[\frac 1 {\operatorname{char }k}]$ instead of $\DM$ is shown in \cite[Cor.~2.1.7]{ElmantoKhan:Perfection}.
The equivalence is also shown in \cite[\S2.2]{RicharzScholbach:MotivicWitt}.
This equivalence is compatible with $\boxtimes$ and !-pushforwards, so we obtain our claim in this case.
\xrema

In order to state that the convolution product preserves weights, we need to talk about weights on equivariant motives.
The idea is simple: a $G$-equivariant motive is declared to be of weights $\ge 0$ or $\le 0$ if its corresponding underlying non-equivariant motive has the corresponding property.
The following definition makes this precise. Note that we only define a pair of full subcategories of $\DM(G \bsl X)$. We do not claim they constitute a weight structure, i.e., we do not assert the existence of weight truncation triangles.

\defi
\thlabel{weights.DM.G}
Suppose a smooth group scheme $G$ acts on a scheme $X$ (both supposed to be separated and of finite type over $S$).
We define
$$\eqalign{\DM(G \bsl X)^{\wt \le 0} & := \lim \DM^{!, \wt \le 0}(\BarC(G, X)) \cr
\DM(G \bsl X)^{\wt \ge 0} & := \lim \DM^{!, \wt \ge 0}(\BarC(G, X)),}$$
i.e., we apply $\DM(-)^{\wt \le 0}$ to each term in the bar construction (see around \inter{Def.~\iref{defi--DM.G}}), with transition functors given by !-pullback.
\xdefi

\rema
The smooth transition maps preserve the subcategory of objects of weight $\ge 0$ and also $\le 0$, so the limits make sense.

The definition is independent of the choice of the presentation of $G \bsl X$: if $G \bsl X = G' \bsl X'$, then $M \in \DM(G \setminus X)$ is of weights $\le 0$, say, if it is so on $X$ (under $!$-pullback) and therefore on $G' \x X = G \x X'$.
Here and in the following we use the standard weight preservation properties under smooth pullback and (in  \refeq{DM.colimitissimo} below) also under proper pushforward \cite[Thm.~3.8]{Hebert:Structure}.
The projection map $G \x X' \r X'$ is smooth and surjective, hence $M|_{X'}$ is of weights $\le 0$ by \thref{weight.creation}.
\xrema

The following weight preservation property will be central to the stability of the Satake category $\Sat_G$ under convolution (see \thref{Satake.convolution}).
We only consider compact objects since it is enough for our purposes, but the statement could be extended to arbitrary ones, at the expense of a more lengthy discussion of weights in that case.
Recall from \inter{Cor.~\iref{coro--Ind.Artin.co.limit} and Prop.~\iref{prop--DM.G.homotopy.invariant}} that
$$\DM(\calP_{\bbf} \backslash \Fl_\bbf)^\comp = \colim_{(t_{ij})_!} \DM(\calP_{\bbf} \backslash \Fl_{\bbf, i})^\comp = \colim_{(t_{ij})_!} \DM(\calP_{\bbf, i} \backslash \Fl_{\bbf, i})^\comp,$$
where $\Fl_\bbf = \colim \Fl_{\bbf, i}$ is a presentation as an ind-scheme (with transition maps denoted by $t_{ij}$) and $\calP_{\bbf, i}$ is an appropriate finite-type quotient of $\calP_\bbf$ acting on $\Fl_{\bbf, i}$, and $\colim$ denotes the colimit in the \ii-category of \ii-categories, which (see loc.~cit.) can in this case just be thought of as the union of the above \ii-categories, as $i$ grows.
Using that $(t_{ij})_!$ is weight-exact, we define
$$\DM(\calP_{\bbf} \backslash \Fl_\bbf )^{\comp, \wt \le 0} := \colim_{(t_{ij})_!} \DM(\calP_{\bbf, i} \bsl \Fl_{\bbf, i})^{\comp,\wt \le 0}\eqlabel{DM.colimitissimo}$$
and likewise for $\ge 0$.

\prop
\thlabel{convolution.DTM.weights}
The convolution product for Tate motives is weight-exact, i.e., for any objects $A, B \in \DTM(\calP_\bbf \backslash \Fl_\bbf)^{\comp, \wt \le 0}$, their convolution $A \star B$ is also of weights $\le 0$ and likewise with $\ge 0$.
\xprop

\pf
We first assume $S$ is (the spectrum of) a perfect field $k$ or just $S = \Spec \Fp$, in which case we show the stronger statement that the convolution product on $\DM(\calP_\bbf \backslash \Fl_\bbf)^\comp$ (as opposed to $\DTM$) preserves weights.
We use the notation in the diagram \refeq{convolution} and the discussion around it.
The functor $- \boxtimes - : \DM(X) \x \DM(\calP_{\bbf, i} \setminus Y) \r \DM(X \x \calP_{\bbf, i} \setminus Y)$ preserves weights by \thref{weights.DM.G} and \thref{exterior.product.weights}.
Similarly, $(e \circ \tilde p)^!$ preserves weights: $e \circ \tilde p$ is a smooth map of prestacks, i.e., it admits a smooth covering on which the map is a smooth map of finite type $k$-schemes. Hence $(e \circ \tilde p)^!$ preserves weights by the same argument as in \thref{weights.DM.G}.
Finally, $\tilde m$ is proper so that $\tilde m_! $ preserves weights again.
Hence the non-equivariant motive underlying $A \star B$, namely $\tilde m_! (A \ttw B)$ has the same weights as $A \boxtimes B$, so that the same holds for $A \star B$ itself.

For general $S$, we consider parallely the structural map $f: T := S \r \Spec \Z$ and the closed immersion $f: T := \Spec \Fp \r \Spec \Z$.
The functor $f^* : \DTM(\calP_{\bbf, \Z} \backslash LG_\Z / \calP_{\bbf, \Z}) \r \DTM(\calP_{\bbf, T} \backslash LG_T / \calP_{\bbf, T})$ creates the weight structure by \thref{f!.nice}.
We conclude by using \thref{f*.convolution}.
\xpf

\section{Mixed Tate motives on the affine Grassmannian}
\label{sect--tate.aff.grass}

In this section, we endow the category $\MTM(L^+G\backslash LG/L^+G)$ with a Tannakian structure, cf.~\thref{Tannaka_Cat}.
A recurrent idea is that the conservativity of the $\ell$-adic realization functor restricted to stratified Tate motives allows us to lift many statements from the $\ell$-adic to the motivic setting.

\syno
\thlabel{MTM.basics}
Throughout \refsect{tate.aff.grass}, the base scheme $S$ is as in \thref{base.scheme}.
We fix a split reductive group $G\to S$ and a Borel pair $T\subset B\subset G$ over $S$ where $T$ is a split maximal torus contained in the Borel subgroup $B$.
 We start by listing some basic properties as needed in the following, see~\S\ref{sect--motive.affine.flag.basics} for more details and the references cited there.
\begin{enumerate}
\item \label{item--MTM.basics.1}
The affine Grassmannian $\Gr=\Gr_G\to S$ is the quotient of ordinary \'etale sheaves $(LG/L^+G)^\et$ with base point denoted by $e\in \Gr(S)$.
For each dominant cocharacter $\mu\in X_*(T)^+$, the locally closed immersion of the $L^+G$-orbit of $\varpi^\mu\cdot e$ is denoted by
\[
\iota_\mu\co \Gr^\mu \;\stackrel {j_\mu} {\hookrightarrow} \; \Gr^{\leq \mu}\;\stackrel {i_\mu}{\hookrightarrow}\; \Gr.
\]
Each Schubert scheme $\Gr^{\leq \mu}\to S$ is proper, and the open orbit $\Gr^\mu\subset \Gr^{\leq \mu}$ is $S$-smooth, fibrewise dense with geometrically connected fibers of relative dimension $\lan2\rho,\mu\ran\in \bbZ_{\geq 0}$.
For $\la,\mu\in X_*(T)^+$, we have $\Gr^\la\subset \Gr^{\leq \mu}$ if and only if $\la \leq \mu$ in the dominance partial order on $X_*(T)^+$.
For details see \inter{Exam.~\iref{exam--length.function.exam}} and \inter{Lem.~\iref{lemm--orbit.flag}}.

\item \label{item--MTM.basics.2}
Throughout, we work with the prestack double quotient $L^+G\backslash LG/L^+G\to S$, cf.~\inter{\iref{sect--DTM.double}}.
For each point $s\in S$, there is a homeomorphism of topological spaces
\[
 |L^+G_s\backslash LG_s/L^+G_s|\simeq X_*(T_s)^+\overset{{\text{$S$ connected}}}{=}X_*(T)^+,
\]
where $X_*(T)^+$ is endowed with the topology given by the dominance partial order ``$\leq$'', i.e., for each $\mu\in X_*(T)^+$ the subset $\{\la\,|\,\la\leq \mu\}$ is closed.

\item \label{item--MTM.basics.3}
The étale descent equivalence $\DM(L^+G\bsl LG/L^+G)=\DM(L^+G\bsl \Gr)$ (cf.~\thref{properties.DM}) will be used freely throughout.
We have the full subcategory
\[
\DTM(L^+G\bsl LG/L^+G)\subset \DM(L^+G\bsl LG/L^+G)
\]
of stratified Tate motives.
For each $\mu\in X_*(T)^+$, there is an equivalence $\DTM(L^+G\bsl \Gr^\mu)=\DTM(H_\mu\bsl S)$ where $H_\mu=L^+G\cap (\varpi^\mu L^+G\varpi^{-\mu})$ denotes the stabilizer, cf.~ \inter{Prop.~\iref{prop--DTM.G/H}}.
Here $H_\mu\to S$ is a fibrewise connected, strictly pro-algebraic closed subgroup of $L^+G$ by \inter{Lem.~\iref{lemm--orbit.flag} ii)}, that is, it can be written as a sequential limit of fiberwise connected, $S$-smooth, $S$-affine group schemes with surjective transition maps.

\item
\label{item--MTM.basics.4}
The category of stratified Tate motives admits a non-degenerate $t$-structure whose heart is the full subcategory of mixed Tate motives
\[
\MTM(L^+G\bsl LG/L^+G)\subset \DTM(L^+G\backslash LG/L^+G).
\]
This category is abelian, $\bbQ$-linear, and the forgetful functor
\[
\MTM(L^+G\backslash LG/L^+G)\;\lr\; \MTM(\Gr)
\]
is fully faithful and induces a bijection on simple objects by \thref{motivic.t.structure}.
(In fact, it is an equivalence by \thref{equivalence} below.)
The simple objects are the intersection motives $\IC_\mu(n)$ for $\mu\in X_*(T)^+$, $n\in \bbZ$. These are pure of weight $\lan2\rho,\mu\ran -2n\in \bbZ$ by \thref{IC.pure}.
The closed set $\{\la\,|\,\la\leq \mu\}\subset X_*(T)^+$ identifies with the closure of the support of $\IC_\mu(n)$ for any $n\in \bbZ$, cf.~\inter{Lem.~\iref{lemm--double.orbit}, \iref{lemm--orbit.flag}}.
Furthermore, for each $\mu\in X_*(T)^+$, there are equivalences $\MTM(\Gr^\mu)=\MTM(L^+G\bsl \Gr^\mu)=\MTM(H_\mu\bsl S)=\MTM(S)$ by \thref{orbit.simply.conn} and \inter{Cor.~\iref{coro--MTM.G.trivial}}.

\item \label{item--MTM.basics.5}
The convolution product defines a functor
\[
\str\star\str\co \DTM(L^+G\bsl LG/L^+G)\x \DTM(L^+G\bsl LG/L^+G)\lr \DTM(L^+G\bsl LG/L^+G)
\]
on the category of stratified Tate motives, cf.~\thref{convolution.Fl.Tate}.
(We show in \thref{convolution.Tate.perverse} below that the convolution product preserves the subcategory $\MTM(L^+G\backslash LG/L^+G)$.
Moreover, the two possible approaches in defining the convolution product discussed in \thref{convolution.HoDM.independent} are isomorphic on this category.)

\item
\label{item--MTM.basics.6}
For each geometric point $f\co \bar{s}\to S$, the fibre of the $\ell$-adic realization functor (cf.~\thref{realization.functor.coro})
\[
\rho_{\ell,\bar s}:=f^*\circ \rho_\ell\co \MTM(\Gr)^\comp\;\to\; \Perv(\Gr_{\bar s},\Ql)
\]
is exact, conservative and faithful.
Under the realization, the motivic convolution product corresponds to the classical convolution product in the $\ell$-adic setting (\thref{convolution.DM.ell}).
\end{enumerate}
\xsyno

\subsection{Indecomposable objects}
The affine Grassmannian admits a decomposition into open and closed sub-ind-schemes $\Gr=\sqcup_{\tau\in \pi_1(G)}\Gr_{\tau}$ where $\pi_1(G)$ is the algebraic fundamental group, i.e., the finitely generated abelian group given by the quotient of the cocharacter lattice by the coroot lattice.
Within each component $\Gr_{\tau}$, every Schubert cell $\Gr^\mu\to S$ has either even- or odd-dimensional fibers: indeed if $\la, \mu \in X_*(T)^+$ with the same class in $\pi_1(G)$, then $\lan2\rho,\la\ran\equiv \lan2\rho,\mu\ran$ in $\bbZ/2$ because $\rho$ takes integer values on the coroot lattice.
This defines the locally constant function
\begin{equation}\label{parity.fct}
 X_*(T)^+\;\to\; \bbZ/2,\;\;\; \mu\mapsto \lan2\rho,\mu\ran \mod 2,
\end{equation}
where the source is equipped with the topology given by the dominance order as in \thref{MTM.basics} \refit{MTM.basics.2}.
An object $A\in \MTM(\Gr)$ is said to have \emph{constant parity} $p(A)\in \bbZ/2$ if the restriction of \eqref{parity.fct} to the closure of its support is a constant function, i.e., the object is supported either on a union of even components, or is supported on a union of odd components.
The following result is a direct consequence of our discussion and the Kazhdan-Lusztig parity vanishing.

\coro
\thlabel{parity.coro.grass}
Let $A\in \MTM(\Gr)$.
\begin{enumerate}
\item \label{item--parity.coro.grass.1}
There exists a canonical decomposition $A=A_{\on{even}}\oplus A_{\on{\odd}}$ into objects of even and odd constant parity.
\item \label{item--parity.coro.grass.2}
If $A\in \MTM(\Gr)$ is of constant parity $p(A)\in \bbZ/2$, then
\[
{^\cl\H}^i(A) = 0,\;\;\;\;\text{whenever \;\; $i\not \equiv p(A)\mod 2$},
\]
where ${^\cl\H}^i$ denotes the truncation with respect to the classical motivic $t$-structure on $\Gr$.

\item \label{item--parity.coro.grass.3}
If $\iota\co X\subset  \Gr$ is a finite union of Schubert schemes, then
\[
{\motH}^i(\iota^*A) = 0,\;\;\;\;\text{whenever \;\; $i\not \equiv 0\mod 2$},
\]
where ${\motH}^i$ denotes the truncation with respect to the perverse motivic $t$-structure on $X$.
\end{enumerate}
\xcoro

\pf Part i) is immediate from the definitions using that each connected component of $\Gr$ is either of constant even or constant odd parity.
Part ii) and iii) for $\ell$-adic sheaves are certainly well-known, and hence are immediate from the conservativity of $\rho_{\ell,\bar s}$ (\thref{MTM.basics} \refit{MTM.basics.6}).
Let us give an argument by reduction to \thref{KL.vanishing.theo}.
The category $\MTM(\Gr)$ is compactly generated, so $A$ is the filtered colimit of its compactly generated subobjects.
Moreover, $^\cl \H$ and $\motH$ commute with filtered colimits, so we may assume $A$ is compact.
For ii), we use that the length function on $X_*(T)^+=W_0\bsl W/W_0$ is computed as $l(\mu)=\lan2\rho,\mu\ran$ which is also the relative dimension of each $\Gr^\mu\to S$.
Hence, for each $\mu\in X_*(T)^+$, the vanishing of ${^\cl\H}^i(\iota_\mu^*A)\in \DTM_{L^+G}(\Gr^\mu)=\DTM_{H_\mu}(S)$ (\thref{MTM.basics} \refit{MTM.basics.3}) in degrees $i$ as above follows from \thref{KL.vanishing.theo}.
Here we use that the forgetful functor $\DM_{H_\mu}(S)\to \DM(S)$ is conservative.
By the compactness of $A$, its support is a finite-type subscheme of $\Gr$.
We can therefore invoke the localization property of $\DM$ to see that the condition
\[
{^\cl\H}^i(\iota_\mu^*A)=\iota_\mu^*{^\cl\H}^i(A)=0
\]
for all $\mu\in X_*(T)^+$ implies the vanishing ${^\cl\H}^i(A)=0$.
\nts{A priori, the localization property states that for a jointly surjective pair $(i, j)$, the functors $(i^* , j^*)$ are jointly conservative. For an ind-scheme $X = \colim X_i$, we additionally use $\DM(X) = \lim \DM(X_i)$. Thus !-restriction to all $X_i$ is conservative. }
For iii), the argument is similar now taking the dimension shifts in the construction of the perverse motivic $t$-structure into account.
\xpf

The next lemma adresses the interplay of Tate motives on the base scheme $S$ and intersection motives.
For $L \in \MTM(S)$ and $\mu\in  X_*(T)^+$, we write
$$\IC_{L, \mu} \defined (i_\mu)_*(j_{\mu})_{!*} L[\dim \Gr^\mu] \in \MTM(\Gr)$$
for the intersection motive twisted by the motive $L$ (more precisely by its $*$-pullback along the projection $\Gr^\mu \r S$).
We have $\IC_{1_S(n),\mu}=\IC_\mu(n)$ for any $n\in\bbZ$.

\lemm
\thlabel{IC.star.skyscraper}
For any $L\in \MTM(S)$, $\mu\in X_*(T)^+$, there is a canonical isomorphism
$$\IC_{L,0} \star \IC_{1_S,\mu} \simeq \IC_{L, \mu}.$$
\xlemm

\pf
Let $i\co S=\Gr^{\leq 0}\to \Gr$ denote the inclusion of the base point, so that $\IC_{L,0} =i_*L$.
In the notation of \refeq{convolution}, both the map $\tilde m$ and the map $e \circ \tilde p$ are isomorphisms when restricted to $\Gr^{\leq 0} \xtw \Gr^{\le \mu}=\Gr^{\le \mu}$.
Thus $\IC_{L,0} \star \IC_{1_S,\mu}= i_* L \ttw \IC_\mu=p_\mu^*L\otimes \IC_\mu=\IC_{L,\mu}$ where $p_\mu\co \Gr^{\le \mu}\to S$.
The last equality is checked using the characterization of $\IC$-motives in \inter{Lem.~\iref{lemm--middle.extension}}.
\xpf



\coro
\thlabel{extensions.grass}
Let $L, L'\in\MTM(S)$ and $\la,\mu\in X_*(T)^+$. Then
\[
\Ext^1_{\MTM(\Gr)}(\IC_{L,\la},\IC_{L',\mu})=
\begin{cases}
\Ext^1_{\MTM(S)}(L,L'), &\text{if $\la=\mu$} \\
0, & \text{else}
\end{cases}
\]
In particular, if $\MTM(S)$ is semisimple \textup{(}e.g.~$S=\Spec(\bbF_q)$\textup{)}, then $\MTM(\Gr)$ is semisimple as well.
\xcoro
\pf
This is immediate from \thref{parity.coro.grass} \refit{parity.coro.grass.3}, and we refer to \cite[Prop.~1 ff.]{Gaitsgory:Central} and \cite[Prop.~3.1]{Richarz:New} for more details. Here is a sketch for the reader's convenience: First assume $\mu=\la$, and denote $\IC:=\IC_{L,\mu}$, $\IC':=\IC_{L',\la}$. Let $\Gr^\mu\overset{j}{\to} \Gr^{\leq \mu}\overset{i}{\leftarrow} \Gr^{\leq\mu}\backslash \Gr^\mu$. We have a long exact localization sequence
\[
\dots\to \Hom(\IC,i_!i^!\IC'[1])\to \Hom(\IC,\IC'[1])\to \Hom(\IC,j_*j^*\IC'[1])\to \Hom(\IC,i_!i^!\IC'[2])\to\dots.
\]
We have the following isomorphisms:
$$\Hom(\IC',j_*j^*\IC'[1])=\Hom(j^*\IC,j^*\IC'[1])=\Ext^1_{\MTM(\Gr^\mu)}(L,L')=\Ext^1_{\MTM(S)}(L,L'),$$
where the last one is a consequence of the equivalence $\MTM(S)\simeq \MTM(\Gr^\mu)$ (\thref{MTM.basics} \refit{MTM.basics.4}).
It is therefore enough to show that the outer groups in the above exact sequence vanish.
The Kazhdan-Lusztig parity vanishing implies that $i^*\IC$, resp.~$i^*\IC'$ lives in perverse degree $\leq -2$. Indeed, by general properties it lives in degree $\leq 0$, degree $0$ vanishes because it is an $\IC$-sheaf, degree $-1$ vanishes by \thref{parity.coro.grass} \refit{parity.coro.grass.3}. By duality $i^!\IC'$ lives in perverse degrees $\geq 2$, and taking the shifts $[1]$, resp.~$[2]$ into account, we see that the outer groups vanish by the axioms of a $t$-structure. This implies the corollary in the case $\mu=\la$. Now let $\mu<\la$, and denote $i\co \Gr^{\leq \mu}\to \Gr^{\leq \la}$ the closed embedding. Again the group
\[
\Hom(i_*\IC_{L,\mu},\IC_{L',\la}[1])=\Hom(\IC_{L,\mu},i^!\IC_{L',\la}[1])
\]
vanishes for $t$-structure reasons as above using the parity vanishing. The case $\la<\mu$ is similar. Now if both $\mu\not\leq \la$ and $\la\not\leq \mu$, there are no extensions between $\IC$-sheaves.
This is shown in loc.~cit.~without appealing to parity vanishing.
\xpf

The category of compact objects $\MTM(\Gr)^\comp$ is both Noetherian and Artinian (\thref{motivic.t.structure}), i.e., each object has finite length. Thus, it is a Krull-Remak-Schmidt category by \cite{Krause:KrullSchmidt}, so that each object $A\in \MTM(\Gr)^\comp$ admits a direct sum decomposition into indecomposable objects $A=A_1\oplus\ldots\oplus A_n$ which is unique up to permutation of the factors.

\coro \thlabel{parity.decomposition}
Let $A\in \MTM(\Gr)^\comp$. Then there exist $L_1,\ldots,L_n\in \MTM(S)^\comp$ indecomposable, and $\mu_1,\ldots, \mu_n\in X_*(T)^+$ such that $A\simeq \oplus_{i=1,\ldots, n}\IC_{L_i,\mu_i}$. Further, each $\IC_{L_i,\mu_i}$ is simple if and only if $L_i\simeq 1_S(n_i)$ for some $n_i\in \bbZ$.
\xcoro
\pf
By the discussion above, we may assume that $A$ is indecomposable in which case we have to show $A\simeq \IC_{L,\mu}$ for some necessarily indecomposable $L\in \MTM^\comp(S)$ and $\mu\in X_*(T)^+$. We proceed by induction on the length $l(A)$. The condition $l(A)=1$ is equivalent to $A$ being simple, and thus $A\simeq \IC_{L,\mu}$ with $L=1(n)$ for some $n\in \bbZ$ (\thref{MTM.basics} \refit{MTM.basics.4}). Let $l(A)\geq 2$, and let $0\not= A'\subset A$ be a subobject of length $l(A')=1$.
In particular, $A'$ is indecomposable.
The quotient $A/A'$ is also indecomposable.
By induction $A'\simeq \IC_{L',\mu}$, $A/A'\simeq \IC_{L,\la}$, and thus $[A]\in \Ext^1_{\MTM(\Gr)}(\IC_{L,\la},\IC_{L',\mu})$. As $A$ is indecomposable, the class $[A]\not = 0$ which by \thref{extensions.grass} implies that $\la=\mu$ and that $A$ is of the desired form.
\xpf

The following result is similar to \cite[Prop.~2.1]{MirkovicVilonen:Geometric}.

\coro \thlabel{equivalence}
The forgetful functor
\[
\MTM(L^+G\bsl LG/L^+G)\;\overset{\simeq}{\lr}\;\MTM(\Gr)
\]
is an equivalence of $\bbQ$-linear abelian categories.
\xcoro
\pf
The functor is fully faithful by \inter{\iref{theo--equivariant.DTM.flag} iii)}. As every object in $\MTM(\Gr)$ is isomorphic to a direct sum of objects $\IC_{L,\mu}$, it is essentially surjective as well.
\xpf

\subsection{Tensor structure} \label{sect--tensor.structure.mtm}
In this section, we show that the convolution product on $\DM(L^+ G \bsl LG / L^+ G)$ preserves the subcategory $\MTM(L^+G\bsl LG/L^+G)$.

\lemm \thlabel{convolution.Tate.perverse}
For $A,B \in \MTM(L^+G\bsl LG/L^+G)$ one has $A\star B \in \MTM(L^+G\bsl LG/L^+G)$, i.e., the category $\MTM(L^+G\bsl LG/L^+G)$ is stable under the convolution product.
\xlemm

\pf
By \thref{convolution.Fl.Tate}, convolution preserves Tateness, i.e., $A\star B \in \DTM(L^+G\bsl LG/L^+G)$. The remaining property $A\star B \in \MTM(L^+G\bsl LG/L^+G)$, equivalently ${^m \H}^i(A\star B)=0$ for all $i\not = 0$ is shown using the isomorphism
\[
\rho_\ell\big({\motH}^i(A\star B)\big)={^p \H}^i\big(\rho_\ell(A\star B)\big)={^p \H}^i\big(\rho_\ell(A)\star_\ell \rho_\ell(B)\big).
\]
For the classical $\ell$-adic convolution functor $\star_\ell$ and the rightmost isomorphism see around \thref{convolution.DM.ell}.
We now use that at least over a separably closed base field the convolution product of perverse equivariant sheaves is again perverse \cite[Thm.~2.1]{Richarz:New}.
We then conclude using the conservativity of the composite $\rho_{\ell,\bar s}$, cf.~\thref{MTM.basics} \refit{MTM.basics.6}.
\xpf

The following proposition is a subtle part of the geometric Satake equivalence in the different settings \cite{Ginzburg:Perverse, BeilinsonDrinfeld:Quantization, MirkovicVilonen:Geometric, Richarz:New, Zhu:Affine, BaumannRiche:Satake}.
Here we benefit from the existence of the symmetric monoidal structure in these settings and the faithfulness of the $\ell$-adic realization in \thref{MTM.basics} \refit{MTM.basics.6} to check the required compatibilities between the commutativity and associativity constraints.

\prop
\thlabel{constraints}
Let $A,B,C\in \MTM(L^+G\bsl LG/L^+G)$. There exist functorial equivalences
\[
c_{A,B}\co A\star B\simeq B\star A, \;\;\;\text{and}\;\;\; a_{A,B,C}\co (A\star B)\star C\simeq A\star(B\star C),
\]
called commutativity and associativity constraints which are uniquely determined by the following two properties:
\begin{enumerate}
\item[i\textup{)}] The isomorphisms are colimit-preserving in each argument.
\item[ii\textup{)}] For any geometric point $\bar{s}\to S$, the constraints map under the composition of functors \textup{(}\thref{equivalence}, \thref{MTM.basics} \refit{MTM.basics.6}\textup{)}
\[
\MTM(L^+G\bsl LG/L^+G)^c\overset{\simeq}{\lr} \MTM(\Gr)^c\overset{\rho_{\ell,\bar s}}{\lr} \Perv(\Gr_{\bar{s}},\Ql)
\]
to the usual constraints used in geometric Satake as, e.g., in \cite[Prop.~2.21]{Zhu:Affine}.
\end{enumerate}
In particular, the category $\MTM(L^+G\bsl LG/L^+G)\simeq \MTM(\Gr)$ is a symmetric monoidal tensor category with respect to these constraints.
\xprop
\pf
{\it Uniqueness.} By property i) it is enough to characterize the constraints on the subcategory of compact objects. For any $\bar{s}\to S$ as above, the functor $\rho_{\ell,\bar s}\co \MTM(L^+G\bsl LG/L^+G)^c\to \Perv(\Gr_{\bar{s}},\Ql)$ is faithful by \thref{MTM.basics} \refit{MTM.basics.6}. This implies uniqueness.\smallskip\\
{\it Existence.} Note that once $c_{A,B}$, $a_{A,B,C}$ with properties i) and ii) exist, these constraints have to satisfy the axioms required for a symmetric monoidal category (hexagon axiom etc.). This follows from the corresponding identities for the $\ell$-adic Satake equivalence, and the faithfulness of the functor in ii). It remains to construct the constraints. The associativity constraint was constructed in \thref{associative}. For the commutativity constraint, we use the categorical analogue of Gelfand's trick whose construction is explained in \cite[\S5.3.8]{BeilinsonDrinfeld:Quantization} and \cite[\S2.4.3]{Zhu:Affine}. The use of prestacks simplifies the construction a little bit: Fix a pinning $(G,B,T,X)$. Define the anti-involution $\theta\co G\to G, g\mapsto (g^*)^{-1}=(g^{-1})^*$ where $(\str)^*$ denotes the Cartan involution.
The latter is characterized by the fact that it maps a dominant cocharacter $\la\in X_*(T)^+$ to $-w_0\la$, where $w_0$ is the longest element in the finite Weyl group.
By functoriality, we obtain an anti-involution on $L:=LG$ preserving $L^+:=L^+G$, and thus an equivalence of prestacks, still denoted by
\[
\theta\co L^+\bsl L/L^+ \stackrel {\simeq}{\longrightarrow} L^+\bsl L/L^+.
\]
For all $A,B\in \MTM(L^+G\bsl LG/L^+G)$ we construct a canonical isomorphism $\theta^!(A\star B)\simeq (\theta^!B)\star (\theta^!A)$ as follows:
there is a (homotopy) Cartesian diagram of prestacks
\[
\xymatrix{
L^+\bsl L\x^{L^+} L/L^+ \ar[r]^{m} \ar[d]^{\on{sw}\circ \theta\tilde\x\theta} & L^+\bsl L/L^+ \ar[d]^\theta \\
L^+\bsl L\x^{L^+} L/L^+ \ar[r]^m & L^+\bsl L/L^+,
}
\]
where $\on{sw}$ is induced from the switch $L\x L\to L\x L$, $(g_1,g_1)\mapsto (g_2,g_1)$.
Hence, we obtain
\[
\theta^!(A\star B)=\theta^!m_!p^!(A\boxtimes B)=m_!(\on{sw}\circ \theta\tilde\x\theta)^!p^!(A \boxtimes B)= m_!(\theta\tilde\x\theta)^!p^!(B \boxtimes A) \stackrel{\textup{($*$)}} =m_!p^!(\theta^!B\boxtimes\theta^!A)=(\theta^!B)\star (\theta^!A).
\]
The isomorphism labelled ($*$) follows from $(\theta \x \theta)^! (A \boxtimes B) = (\theta^! A) \boxtimes (\theta^! B)$, which holds since $\theta$ is a placid map and $\DM^!$ is symmetric lax monoidal as a functor on placid prestacks with placid maps.
Next we define an isomorphism of (plain) endofunctors on $\MTM(L^+G\bsl LG/L^+G)$ denoted by
\[
e\co \theta^!\stackrel \simeq \longrightarrow \id.
\]
For this we fix a square root $i\in \bbC$ of $-1$, and work temporarily with coefficients in $\bbQ(i)$. We define $e$ on each indecomposable object $\IC_{L,\mu}$ (\thref{parity.decomposition}) to be the map corresponding to $i^{(2\rho,\mu)}\cdot \id$ under
\[
\Hom_{\bbQ(i)}(\theta^!\IC_{L,\mu},\IC_{L,\mu})=\Hom_{\bbQ(i)}(L[d_\mu],L[d_\mu]),
\]
obtained by restriction to the open orbit $L^+G\backslash LG^\mu/L^+G\subset L^+G\backslash LG^{\leq \mu}/L^+G$ (this orbit is invariant under $\theta$).
Since $i^{(2\rho,\mu)}\cdot \id$ is a central endomorphism, one checks that $e$ is functorial. 
We leave the details to the reader.

Next define
\[
c_{A\star B}'\co \;A\star B\;\stackrel {e_{A\star B}}{\longleftarrow} \;\theta^!(A\star B)\;\simeq\; (\theta^!B)\star(\theta^!A)\;\stackrel {e_B\star e_A}{\longrightarrow}\; {B\star A},
\]
which is invariant under the Galois automorphism $i\mapsto -i$, and thus defined over $\bbQ$.
Finally, if both objects $A,B$ have constant parity, we define $c_{A\star B}:=(-1)^{p(A)p(B)}c_{A\star B}'$ where $p$ denotes the parity function from \eqref{parity.fct}.
For general objects $A, B$, not necessarily of constant parity, we use \thref{parity.coro.grass} \refit{parity.coro.grass.1} to extend $c_{A\star B}$ linearly.
See also \cite[\nopp 5.3.21]{BeilinsonDrinfeld:Quantization} and \cite[Rmk.~6.2 ff.]{MirkovicVilonen:Geometric} for slicker formulations.
It is a difficult theorem which is proven in \cite[Prop.~2.21]{Zhu:Affine}, relying on \cite{LusztigYun:Analogue}, that the $\ell$-adic realization $\rho_\ell(c_{A\star B})$ is the (modified) commutativity constraint coming from the fusion interpretation of the convolution product.
This finishes the construction of the constraints.
\xpf

\rema
If $S$ is the spectrum of a field, then the above commutativity constraint can also be constructed as e.g.~in \cite{BeilinsonDrinfeld:Quantization} and \cite{MirkovicVilonen:Geometric} (see also \cite{Gaitsgory:Central}) by using the fusion interpretation of the convolution product and the motivic nearby cycles functor constructed in \cite{Ayoub:Motivic}, \cite[§10]{Ayoub:Realisation}.
\xrema

\subsection{Tannakian structure}
In this section, we show that the category $\MTM(L^+G\bsl LG/L^+G)\simeq \MTM(\Gr)$, which admits a symmetric monoidal structure with respect to the convolution product $\star$ by \thref{constraints}, has in fact a Tannakian structure with fibre functor being the global motivic cohomology functor.

\subsubsection{The fiber functor}

The fiber functor is a motivic analogue 
of the augmentation map for the spherical Hecke algebra.

\defi
\thlabel{fiber.functor}
The \emph{fiber functor} is the composition
$$\omega\co \MTM(L^+G\bsl LG/L^+G) \stackrel {\sigma^!} \lr \MTM (\Gr) \stackrel{\epsilon_!} \lr \DTM(S) \stackrel{\gr^\cl} \lr \MTM(S) \stackrel{\gr^\W} \lr \MTM(S)^{\wt = 0} = \Vect_\Q$$
of the forgetful functor $\sigma^!$ (which is an equivalence of categories by \thref{equivalence}), the pushforward along the structural map $\epsilon\co \Gr \r S$
(which preserves Tate motives by \inter{Lem.~\iref{lemm--Tate.proper.descent}}, using that the stratification of $\Gr$ by $L^+G$-orbits is cellular), followed by the grading functors for the classical motivic (which agrees in this case with the perverse  motivic) $t$-structure and the weight structure (\thref{weight.grading}):
$$\eqalign{
\gr^\cl & : M \mapsto \bigoplus_i {^\cl\H^i}(M), \cr
\gr^\W & : M \mapsto \bigoplus_i (\gr^W_{2i} M) (i),}$$
and finally the equivalence \cite{Levine:Tate} of pure Tate motives of weight 0 with the category of $\Q$-vector spaces (here we use that $S$ is connected).
\xdefi

As a consequence of the Kazhdan-Lusztig parity vanishing one obtains:

\coro\thlabel{parity.coro.grass.last}
Let $A\in \MTM(L^+G\bsl LG/L^+G)$ be of constant parity $p(A)\in \bbZ/2$ \textup{(}\thref{parity.coro.grass}\textup{)}.
Then ${^\cl\H^i}(\epsilon_!\sigma^!A)=0$ whenever $i\not \equiv p(A) \mod 2$.
\xcoro

\pf
This is immediate from the conservativity of the $\ell$-adic realiztaion as in \thref{MTM.basics} \refit{MTM.basics.6} using the well-known statement for $\ell$-adic sheaves.
Here is an argument by reduction to \thref{parity.coro.grass} according to which $A$ has only, say, even classical cohomology, i.e., $^\cl \H^i(A) = 0$ for $i$ odd.

As in the proof of \thref{parity.coro.grass}, we may assume that $A$ is compact.
The functors $(\iota_{\mu})_!$ and $(\iota_{\mu})^*$ between $\DTM(\Gr_G)$ and $\DTM(\Gr_G^\mu)$ are exact with respect to the classical motivic t-structure since the corresponding statement is true for $\ell$-adic sheaves.
Moreover, by localization, $\sigma^! A$ is an iterated extension of the $A_\mu := (\iota_{\mu})_! (\iota_{\mu})^* \sigma^! A$ where $\mu$ runs over the finitely many dominant cocharacters in the support of $A$.
Using the long exact cohomology sequence for the $^\cl \H$-cohomologies of $\epsilon_! \sigma^! A$, we may thus replace $A$ by $A_\mu$.
We now use that the Iwahori stratification of $\Gr^\mu$ consists of affine spaces.
By the same localization argument, we may replace $\Gr^\mu$ by such a stratum $\A^n_S$.
We are left to showing that for the structural map $\A^n_S \stackrel p \r S$, $p_!\co\DTM(\A^n_S) \r \DTM(S)$ preserves parity vanishing.
In fact $p_!$ is an equivalence of categories, being the adjoint of $p^!\co \DTM(S) \r \DTM(\A^n_S)$ which is fully faithful by homotopy invariance of $\DM$ and essentially surjective by definition of $\DTM$.
\xpf

It is a classical fact
that the augmentation map from parahoric Hecke algebras to the coefficient field respects the multiplicative structure.
The corresponding fact is also well-known in a categorified situation.
The proof below is thus similar to, say, \cite[Prop.~2.20]{Zhu:Affine}.

\prop
\thlabel{pushforward.symmetric.monoidal}
Let $\bbf$ be a facet and $\calP := \calP_\bbf$ its associated parahoric subgroup of $LG$. Then the pull-push along the correspondence
$$\calP \bsl LG / \calP \stackrel \sigma \longleftarrow LG / \calP \stackrel \epsilon \lr S$$
yields a functor
$$\omega: \DM (\calP \bsl LG / \calP) \stackrel {\sigma^!} \lr \DM(LG / \calP) \stackrel {\epsilon_!} \lr \DM(S).$$
At least on the level of homotopy categories, this functor has a natural monoidal structure with respect to the convolution product on $\DM (\calP \bsl LG / \calP)$ and the ordinary tensor product on $\DM(S)$.
\xprop

\pf
The left adjoint $\epsilon_!$ of $\epsilon^!$ exists since the étale sheafification $\Fl := \Fl_\bbf = (LG / \calP)^\et$ is an ind-scheme, so the pushforward along the structural map to $S$ exists by \inter{Thm.~\iref{theo--motives.Ind-schemes}}.
The unitality of $\omega$ is clear, since the monoidal unit is the skyscraper motive supported at the base point of $\Fl$.
To show the monoidality of $\omega$,
we abbreviate 
$\tilde X := LG \x LG / \calP$ (all products are products of prestacks over $S$) and 
$X := \calP \backslash LG \x LG / \calP = \calP \bsl \tilde X$.
The group $\calP^\opp \x \calP$ acts on $\tilde X$ (and $X$) by acting on the right on the first, and on the left on the second factor.
We write $\Delta := \{(p^{-1}, p)| p \in \calP\} \subset \calP^\opp \x \calP$ for the ``diagonal'' subgroup.
The structural map $X \r S$ gives a map $\Delta \bsl X \r \Delta \bsl S$.
On the other hand, the composition of the mulitplication and structural map $LG / \calP \stackrel{\epsilon} \r S$ yields another map
$$\Delta \bsl X \stackrel{m} \r \calP \bsl LG / \calP \stackrel {\epsilon'} \r \calP \bsl S.$$
These two maps agree: this can be checked after precomposing with the epimorphism $\tilde X \r (\calP \x \Delta) \bsl \tilde X = \Delta \bsl X$, where it boils down to using that the structural map $\tilde X \r S$ agrees with $\tilde X \stackrel{m} \r LG / \calP \stackrel \epsilon \r S$, since $S$ is is the final object, and in particular is acted upon trivially by all copies of $\calP$.

This shows that the following diagram of prestacks is cartesian, as soon as we omit the dotted map. Note that the map $\epsilon' \x \epsilon'$ arises as $(\calP^\opp \x \calP) \bsl (X \r S)$.
Once we do include the dotted map, the small bottom left square is still cartesian (but the top left square does not commute):
$$\xymatrix{
X \ar[d]^{\pr} \ar[r]
&
\Delta \backslash X \ar[d]^m \ar[r]^(.4)p
&
\calP \bsl LG / \calP \x \calP \bsl LG / \calP \ar[dd]^{\epsilon' \x \epsilon'}
\\
LG / \calP \ar[d]^\epsilon \ar@{.>}[r]^\sigma
&
\calP \backslash LG / \calP \ar[d]^{\epsilon'}
\\
S \ar[r]^{\sigma_S}
&
\Delta \backslash S \ar[r]
&
\calP \bsl S \x \calP \backslash S
}$$
We then have $\omega(A \star B) = \epsilon_! \sigma^! m_! p^! (A \boxtimes B)$.
As was already noted in the proof of \thref{star.twiddle}, $\epsilon_!$ exchanges with the forgetful functor $\sigma^!$ by ind-proper base change. Similarly, $\sigma_S^!$ exchanges with $(\epsilon' \circ m)_!$, so the above is equivalent to $(\epsilon \x \epsilon)_! (\sigma \x \sigma)^! (A \boxtimes B)$.
The !-pullback along the map $X \r (\calP^\opp \x \calP) \bsl X$ commutes with $\boxtimes$ by construction of $\boxtimes$, see \thref{explain.convolution}.
Furthermore, $\boxtimes$ also commutes with the !-pushforward along the structural map $X \r S$.
After reducing this claim to the case of finite-type $S$-schemes (instead of the ind-finite type ind-scheme $X$), this is a consequence of the projection formula.
Hence the above object is equivalent to
$\omega(A) \t \omega(B)$.
\xpf

\theo
\thlabel{Tannaka_Cat}
The category of compact objects $\MTM(L^+G\bsl LG/L^+G)^\comp=\MTM(\Gr)^c$ \textup{(}\thref{equivalence}\textup{)}, endowed with the convolution product, the constraints from \thref{constraints} and $\omega$ as fibre functor
is a neutral Tannakian category over $\bbQ$ \textup{(}\cite[Ch.~II, Def.~2.19]{DeligneMilneOgusShih:Hodge}\textup{)}.
\xtheo

\pf
We check the conditions in \cite[Ch.~II, Prop.~1.20]{DeligneMilneOgusShih:Hodge} using the functor $\omega$ (\thref{fiber.functor}).
\begin{enumerate}
\item {\it The functor $\omega$ has a monoidal structure.}
By \thref{pushforward.symmetric.monoidal}, it remains to observe that $\gr^\perv$ and $\gr^\W$ have monoidal structures:
both functors and the respective tensor products preserve colimits, so we may consider the subcategories of compact objects instead.
Then we are in the situation of \cite{Levine:Tate}: the characterization of $\DTM(S)^{\le 0, \comp}$ (resp.~$\DTM(S)^{\ge 1, \comp}$) as the subcategories generated under extensions by the objects $1_S(n)[m]$ with $n \in \Z$ and $m \ge 0$ (resp. $m \le -1$) immediately shows the monoidality of $\gr^\cl$ restricted to the subcategory of complexes concentrated in either even or odd degrees.
Now we take the parity vanishing of \thref{parity.coro.grass.last} into account.
The monoidality of $\gr^\W$ holds by \cite[Thm.~1.4.v]{Levine:Tate}.

\item {\it The functor $\omega$ is $\bbQ$-linear, exact and faithful.} The functor $\omega$ is clearly $\Q$-linear. To check the exactness, suppose
$$0 \r M' \r M \r M'' \r 0$$
is a short exact sequence in $\MTM(L^+G\bsl LG/L^+G)= \MTM(\Gr)$.
To show it maps to an exact sequence under $\omega$, we may assume $M'$ and $M''$ are indecomposable (since $\Ext^1_{\MTM(\Gr)}(\str,\str)$ commutes with finite direct sums in each variable), i.e., by \thref{parity.decomposition} of the form $M' = \IC_{L', \mu'}$, $M''=\IC_{L'', \mu''}$ for some $\mu',\mu''\in X_*(T)^+$ and some indecomposable motives $L', L''\in \MTM(S)$.
By \thref{extensions.grass}, the extension splits unless $\mu' = \mu''=:\mu$ in which case $M\simeq \IC_{L, \mu}$, where $L$ is an extension in $\MTM(S)$ of $L'$ and $L''$.
We conclude that $\omega$ is exact since $\omega (\IC_{L, \mu}) \simeq \omega(L\star \IC_{1,\mu}) \simeq \omega(L) \t \omega(\IC_{1, \mu})$ using i) above.
The faithfulness of $\omega$ follows from the conservativity of $\omega$, which in its turn follows from the conservativity of the $\ell$-adic realization at some geometric point $\ol s \r S$ (\thref{MTM.basics} \refit{MTM.basics.6}) and the conservativity of the fiber functor in the $\ell$-adic situation.

\item {\it The constraints constructed in \thref{constraints} give the usual constraints in $\Vect_\Q$ after applying $\omega$.} 
This is immediate from the $\ell$-adic case, see \cite[Prop.~2.21]{Zhu:Affine} using that the realization is faithful.
We stress that one needs to change the natural commutativity constraints one has on complexes:
this is possible by the parity vanishing, see also the passage from $c_{A\star B}'$ to $c_{A\star B}$ in the proof of \thref{constraints}.

\item {\it Neutral object.} Clearly, the skyscraper at the base point $\IC_0 \in \MTM(L^+G\bsl LG/L^+G)$ satisfies $\End(\IC_0)=\Q$ and $\omega(\IC_0)$ is 1-dimensional.
\item {\it Any $M\in \MTM(L^+G\bsl LG/L^+G)$ with $\dim_\bbQ \omega(M)=1$ admits a dual object $M^{-1}$ such that $M \star M^{-1} = \IC_0$.} If $\dim_\bbQ \omega(M)=1$, then $M$ is indecomposable by the faithfulness of $\omega$, i.e.,  $M\simeq \IC_{L, \mu}$ with $L$ indecomposable. Since $\omega(\IC_{L, \mu})=\omega( L) \t \omega(\IC_{1,\mu})$, the motive $L\in \MTM(S)$ is also pure, and hence $\t$-invertible.
Moreover, $\dim_\Q \omega(\IC_\mu) = \dim_{\Ql} \omega_\ell(\IC_{\mu, \ell})$, where the subscripts denote the corresponding functors in the $\ell$-adic realization.
By the $\ell$-adic Satake equivalence (see e.g.~\cite[\S 9]{BaumannRiche:Satake}\footnote{The reference uses constructible sheaves in the analytic topology over $\bbC$, but the same argument works by invoking \cite[Cor.~6.9]{HainesRicharz:TestFunctions}.}, or \cite[Cor.~3.5]{Richarz:New}), this implies $\dim(\Gr^{\le \mu}/S) = 0$, i.e., $\IC_\mu$ is dualizable with respect to $\star$. Namely its dual is $\IC_\mu\star\IC_{-\mu}\simeq \IC_0$.
\end{enumerate}
\xpf

\section{The dual group}\label{sect--dual.grp.sect}

In this final section we determine the Tannaka dual of the categories $\MTM(L^+G_S\bsl LG_S/L^+G_S)$ and the so-called Satake category
$$\Sat_{G, S} \subset \MTM(L^+G_S\bsl LG_S/L^+G_S),\eqlabel{Sat.MTM}$$ which can be thought of as the semi-simplification of the latter category.
We show in \thref{Satake} that the Tannaka dual of $\Sat_{G, S}$ is Deligne's modified Langlands dual group $\widehat G_1/\bbQ$ as constructed in \cite[\S 2]{FrenkelGross:Rigid}.
In particular, this group is independent of the (connected) base scheme $S$.

For $S = \Spec \Fq$, the inclusion \refeq{Sat.MTM} is an equivalence.
For more general bases $S$, we show in \thref{full.Tannaka} that the Tannakian group of $\MTM(L^+G_S\bsl LG_S/L^+G_S)$ is the semi-direct product of $\widehat{G}_1$ with a pro-unipotent affine group scheme coming from extensions between Tate motives on $S$.

Throughout \S\ref{sect--dual.grp.sect}, the base scheme $S$ is as in \thref{base.scheme}. Also recall \thref{MTM.basics}.

\subsection{The Satake category}
\label{sect--Satake.category}

\defi
The {\em Satake category} $\Sat_G=\Sat_{G,S}$ is the full subcategory of $\MTM(L^+G\bsl LG/ L^+G)$ generated by means of arbitrary direct sums (as opposed to allowing extensions) by the intersection motives $\IC_\mu(n)$, $\mu \in X_*(T)^+$, $n \in \Z$.
\xdefi

\lemm
\thlabel{Hom.IC}
For $L, L' \in \MTM(S)$ and $\la, \mu \in X_*(T)^+$, we have natural identifications
$$\Hom_{\MTM(\Gr_G)}(\IC_{L, \la}, \IC_{L', \mu}) = \begin{cases}
\Hom_{\MTM(S)}(L, L') & \la = \mu, \cr
0 & \la \ne \mu.
\end{cases}$$
\xlemm

\pf
For $\la = \mu$, this is a standard property of intermediate extensions, see \cite[Cor.~III.5.11]{KiehlWeissauer:Weil}.
To show the vanishing in case $\la \ne \mu$, we may assume $L$ is a simple object of $\MTM(S)$.
In this case, $\IC_{L, \mu}$ is also simple, so any non-zero morphism would need to be an isomorphism, which is impossible if $\la \ne \mu$.
\xpf

\coro
The category $\Sat_G$ is abelian. Its subcategory of compact objects $\Sat_G^\comp$ is semi-simple.
\xcoro
\qed

\lemm
\thlabel{Satake.convolution}
The full subcategory $\Sat_G\subset \MTM(L^+G\bsl LG/L^+G)$ is stable under the convolution product.
\xlemm

\pf
We have to show that $M := \IC_\mu \star \IC_\lambda$ is a direct sum of some intersection motives of the form $\IC_\kappa(n)$.
(A priori we only know it is a successive extension of twists of some $\IC_\kappa$.)
By \thref{parity.decomposition}, $M=\oplus_{(L,\mu)} \IC_{L,\mu}$ with $L$ indecomposable.
The intersection motives $\IC_\mu$ and $\IC_\la$ are pure by \thref{IC.pure}, hence so is $M$ by \thref{convolution.DTM.weights}.
Therefore, each direct summand $\IC_{L,\mu}$ is also pure.
Let $j\co \Gr_\mu\to \Gr^{\leq \mu}$ be the open stratum. Since $j^*=j^!$, the motive $j^*\IC_{L,\mu}=L[d_\mu]$ is also pure, which implies $L$ is pure.
Since $L$ is also indecomposable, it is of the form $L=1_S(n)$ for some $n\in \bbZ$, hence $\IC_{L,\mu}=\IC_\mu(n) \in \Sat_G$.
\xpf

\coro
\thlabel{Tannaka.Sat}
The subcategory $\Sat_G^\comp \subset \Sat_G$ spanned by the compact objects in the Satake category has the following properties:
\begin{enumerate}
\item
$\Sat_G^\comp$ is a neutral Tannakian subcategory.
\item For any map $f\co T\to S$ of connected schemes as in \thref{base.scheme}, there is an equivalence of neutral Tannakian categories
\[
f^*\co \Sat_{G, S}^c\stackrel \simeq \longrightarrow \Sat_{G,T}^c,
\]
having the property $f^*\IC_{\mu,S}(n)=\IC_{\mu,T}(n)$ for all $\mu\in X_*(T)^+$, $n\in\bbZ$.
\end{enumerate}
\xcoro
\pf Part i) is immediate from \thref{Satake.convolution} and \thref{Tannaka_Cat}. For ii), we use \thref{base.change.t.exact} which gives $f^*\IC_{\mu,S}(n)=\IC_{\mu,T}(n)$, so that $f^*$ is an equivalence of $\bbQ$-linear abelian categories. The compatibility of $f^*$ with the convolution product was checked in \thref{f*.convolution}. Also $\omega(\IC_{\mu,S}(n))=\omega(\IC_{\mu,T}(n))$
is immediate from \thref{fiber.functor}. The rest is clear from the characterization of the constraints in \eqref{constraints}.
\xpf

Now fix a pinning $(G,B,T,X)$, and denote by $(\widehat{G},\widehat{B},\widehat{T},\widehat{X})$ the dual group in the sense of Langlands formed over $\bbQ$.
By definition $\widehat{G}$ is a split reductive $\bbQ$-group with split maximal torus $\widehat{T}$, and Borel subgroup $\widehat{B}$.
Denote by $\widehat T_\ad$ the image of $\widehat T$ under the map $\widehat G\to \widehat G_\ad$ to the adjoint group.
Then we may view the half sum $\rho$ of the roots in $B$ (=coroots in $\widehat B$) as a cocharacter $\rho\co \bbG_{m,\bbQ}\to \widehat T_\ad\subset \widehat{G}_\ad$.
We let $\bbG_{m,\bbQ}$ act through $\rho$ by inner automorphisms on $(\widehat{G},\widehat{B},\widehat{T})$ from the right.
Colloquially speaking, this action is given by the formula $g\cdot \la=\rho(\la)^{-1}g\rho(\la)$.
We consider the semi-direct product $\widehat{G}_1:=\widehat{G}\rtimes \bbG_{m,\bbQ}$ which is again a split reductive $\bbQ$-group with Borel pair $ \widehat{T}\times \bbG_{m,\bbQ}=:\widehat{T}_1\subset \widehat{B}_1:=\widehat{B}\rtimes \bbG_{m,\bbQ}$. 

For each $\mu\in X_*(T)^+$, and $n\in \bbZ$ we get an irreducible algebraic $\widehat{G}_1$-representation \cite[Ch.~II.5]{Jantzen:Representations}
\[
V_\mu(n)\defined \Ind_{\widehat{B}_1^\opp}^{\widehat{G}_1}(\mu_n),
\]
where $\widehat{B}_1^\opp\subset \widehat G_1$ denotes the Borel opposite to $\widehat{B}_1$, and $\mu_n\co \widehat{B}_1^\opp\to \widehat{T}_1\to \bbG_{m,\bbQ}$ is the composition of the projection with the character $(\mu,n)\in X_*(T)^+\x \bbZ=X^*(\widehat{T}_1)^+$.
Then $V_\mu(n)$ is the representation of $\widehat G_1$ of highest weight $(\mu,n)$.
We denote by $\Rep_{\bbQ}^\fd(\widehat{G}_1)$ the category of algebraic $\widehat{G}_1$-representations on finite-dimensional $\bbQ$-vector spaces.
This category is semi-simple with simple objects the highest weight representations as above.


\rema
\thlabel{epsilon.remark}
The split reductive group $\widehat G_1$ is Deligne's modified Langlands dual group constructed in \cite{FrenkelGross:Rigid}, see also \cite{Deligne:Letter2007}.
More precisely, one checks that the map $(g,\la)\mapsto (g\cdot (2\rho)(\la),\la^2)$ induces a short exact sequence of $\bbQ$-group schemes
\[
1\to \Bmu_2\to \widehat G\x \bbG_{m,\bbQ}\to \widehat G_1\to 1,
\]
where $\Bmu_2\simeq \bbZ/2$ is the constant subgroup scheme generated by the element $(\epsilon,-1)$, $\epsilon:=(2\rho)(-1)$.
It follows that the semi-direct product $\widehat G_1=\widehat G\rtimes \bbG_{m,\bbQ}$ is (canonically) a direct product if $\epsilon=1$.
The latter condition is also equivalent to $\rho$ being a cocharacter of $\widehat T$ (as opposed to $\widehat T_\ad$).
For example, this is the case if $G$ is simply connected, so that $\widehat{G}_\ad=\widehat{G}$ is adjoint.
We note that the difference of $\widehat G$ versus $\widehat G_1$ relates to the notions of $L$-algebraic versus $C$-algebraic as introduced by Buzzard and Gee in \cite{BuzzardGee:Conjectures}.
For further discussion and examples we refer to \cite[Prop.~5.39 ff.]{BuzzardGee:Conjectures}.
\xrema

\theo
\thlabel{Satake}
There is an equivalence of Tannakian categories
\[
\big (\Sat_G^\comp,\star,\omega\big)\simeq \big(\Rep_{\bbQ}^\fd(\widehat{G}_1),\t, v\big),
\]
where $v\co \Rep_{\bbQ}^\fd(\widehat{G}_1)\to \on{Vec}_\bbQ$ denotes the forgetful functor. The intersection motives $\IC_\mu(n)$ correspond to the irreducible $\widehat{G}_1$-representations $V_\mu(n)$ for $(\mu,n)\in X_*(T)^+\x \bbZ=X^*(\widehat{T}_1)^+$.
\xtheo
\pf
We denote by $\Aut^\star_{\Sat_G}(\omega)$ the affine $\bbQ$-group scheme of tensor automorphisms of $\omega$ provided by the neutral Tannakian category $(\Sat_G^c,\star,\omega)$, cf.~\cite[Ch.~II, Thm.~2.11]{DeligneMilneOgusShih:Hodge}. The Satake category $\Sat_{G}=\Sat_{G,S}$ is independent from the (connected) base scheme $S$ by \thref{Tannaka.Sat} ii), i.e., for any map $T\to S$ of schemes as in \thref{base.scheme} we have
\[
\Aut^\star_{\Sat_{G,S}^\comp}(\omega) = \Aut^\star_{\Sat_{G,T}^\comp}(\omega).
\]
As $G$ is split, it is defined over $\bbZ$, and in particular over any scheme as is the affine Grassmannian.
By \thref{Tannaka.Sat}, we may assume that $S=\Spec(\bbF_p)$ is a finite field, and we denote $\widetilde G_1:=\Aut^\star_{\Sat_G^\comp}(\omega)$.
For any $\ell\not = p$, consider the $\ell$-adic realization
\[
\rho_\ell\co \Sat_{G}^\comp\to \Perv_{L^+G}(\Gr,\bbQ_\ell).
\]
Let $\Sat^\comp_{G, \ell}$ be the essential image of $\rho_\ell$, i.e., the full subcategory of $\Perv_{L^+G}(\Gr,\bbQ_\ell)$ consisting of the objects isomorphic to $\rho_\ell (\bigoplus_i \IC_{\mu_i}(n_i)) = \bigoplus_i \IC_{\ell, \mu_i}(n_i)$.
Under the $\ell$-adic convolution product it is a Tannakian category; its fiber functor, denoted by $\omega_\ell$, is defined the same way as $\omega$.
For any $M, N \in \Sat^\comp_{G, \ell}$, \thref{Hom.IC} and analogous computations in the $\ell$-adic context give natural isomorphisms
$$\Hom_{\Sat_G^\comp}(M, N) \t_\Q \Ql = \Hom_{\Sat^\comp_{G, \ell}}(\rho_\ell (M), \rho_\ell(N)),$$
so that
$$\Aut^\star_{\Sat^\comp_G}(\omega) \t_\Q \Ql = \Aut^\star_{\Sat^\comp_{G, \ell}}(\omega_\ell).$$

By the $\ell$-adic geometric Satake equivalence (in particular \cite[Prop.~A.6]{RicharzZhu:Ramified}, \cite[Rmk.~2.10]{HeinlothNgoYun:Kloosterman}, \cite[\nopp 5.5.14]{Zhu:Introduction}), we deduce an isomorphism
\begin{equation}\label{iso.groups}
\widetilde G_1\otimes_\bbQ \bbQ_\ell\simeq \widehat{G}_1\otimes_\bbQ \bbQ_\ell.
\end{equation}
This holds for all prime numbers $\ell\not = p$. Since we have a zig zag $\bbF_p\gets \bbZ\to \bbF_{p'}$ for any pair of primes $p, p'$, we see that \eqref{iso.groups} also holds for $\ell=p$. Using \eqref{iso.groups} for one prime number $\ell$, it follows that $\widetilde G_1$ is (geometrically) connected and reductive: we use fpqc descent for the extension $\bbQ_\ell/\bbQ$ for the properties `(geometrically) connected', `of finite type' and `smooth' \cite[04KV, 02KZ, 02VL]{StacksProject}, and further that the unipotent radical commutes with arbitrary field extensions over perfect fields for the property `reductive' \cite[Prop.~1.1.9]{ConradGabberPrasad:PseudoReductive}.
If $\widetilde G_1$ is also split, then we have $\widetilde G_1\simeq \widehat{G}_1$ over $\bbQ$ by the isomorphism theorem (\cite[Thm.~6.1.17]{Conrad:Groups}).

However, the condition \eqref{iso.groups} for all primes $\ell$ is not enough to ensure that $\widetilde G_1$ is split (\cite{Gross:Groups}), and we argue as follows.
By construction, the functor $\omega\co \Sat_{G}^\comp\to \on{Vec}_\bbQ$ is equipped with a $\bbZ$-grading coming from the grading on ${^\cl\H}^*$ in \thref{fiber.functor}.
This defines a cocharacter $2\widetilde\rho\co \bbG_{m,\bbQ}\to \widetilde G_1$ via the Tannakian formalism, and we denote by $\widetilde T_1\subset \widetilde G_1$ its centralizer.
Under the $\ell$-adic realization we have $2\widetilde\rho=2\rho$ by \cite[\nopp 5.3.20]{Zhu:Introduction}. 
Then, for all primes $\ell$, we have
\begin{equation}\label{iso.groups.tori}
\widetilde T_1\otimes_\bbQ\Ql\simeq \widehat T_1\otimes_\bbQ\Ql,
\end{equation}
compatible with the isomorphism \eqref{iso.groups}.
Using \eqref{iso.groups.tori} for one prime, this implies as above that $\widetilde T_1$ is a commutative reductive group over $\bbQ$, and thus must be a torus. 
As $\widehat T_1\subset \widehat G_1$ is a maximal torus \eqref{iso.groups.tori} implies that $\widetilde T_1\subset \widetilde G_1$ is a maximal torus (as this can be checked over the algebraically closed overfield $\bar\bbQ_\ell$).
However, the $\Gal(\bar \bbQ/\bbQ)$-Galois representation $X_*(\widetilde T_{1,\bar \bbQ})$ is trivial at all primes $\ell$, and hence must be trivial by Minkowski's theorem.
This shows that $\widetilde T_{1}$ is a maximal split torus of $\widetilde G_1$.
\xpf

\subsection{The full Tannakian group}

\defi
Let $\widetilde G_S$ be the Tannaka dual group of the Tannakian category $\MTM(L^+G_S \bsl LG_S / L^+G_S)^\comp$ (\thref{Tannaka_Cat}).
\xdefi

\exam
\thlabel{finite.field.exam}
For $S = \Spec \Fq$, the inclusion $\Sat_G^\comp \r \MTM(L^+G \bsl LG / L^+G)^\comp$ is an equivalence.
This follows from \thref{extensions.grass}, and the semisimplicity of $\MTM(\Fq)^\comp$ (which follows from $\Hom_{\DTM(\Fq)}(1, 1(i)[n])=0$ unless $n=i=0$ because the higher algebraic $K$-theory of finite fields is torsion by Quillen's computation \cite{Quillen:FiniteField}).
Therefore
$$\widetilde G_{\Fq} = \widehat G_1.$$
\xexam

We now exhibit the relation of the Satake category $\Sat_G^\comp$ and $\MTM(L^+G\bsl LG/L^+G)^\comp=\MTM(\Gr_G)^\comp$ (\thref{Tannaka_Cat}) over other bases.
In short, the category $\MTM(\Gr_G)^\comp$ arises by amalgamating Tate motives 
on the base $S$ together with the intersection motives arising from the presence of the group $G$.

Recall from \cite{Levine:Tate} that $\MTM(S)^\comp$ is a neutral Tannakian category with fibre functor
\[
\MTM(S)^\comp \stackrel{\gr^\W} \r \MTM(S)^{\comp, \W=0} = \Vect_\Q,
\]
which is the special case of \thref{Tannaka_Cat} for the trivial group.
We denote its Tannaka dual by $\scrG_S$. (It is the same as $\widetilde G_S$ for the trivial group $G=1$).
The semi-simplification of $\MTM(S)^\comp$ is just $\MTM(S)^{\pure, \comp}$, which is equivalent to $\Z$-graded $\Q$-vector spaces. Equivalently, its Tannaka dual is $\GmX \Q$.
Therefore $\scrG_S$ sits in a split exact sequence
$$1 \r \scrU_S \r \scrG_S \r \Gm \r 1,$$
where $\scrU_S$ is the pro-unipotent radical of $\scrG_S$.
We have the following commutative diagram of neutral Tannakian categories
$$\xymatrix{
\MTM(S)^{\pure, \comp} \ar@{^{(}->}[r] \ar[d]^{e_*} & \MTM(S)^\comp \ar[d]^{e_*} \\
\Sat_G^\comp \ar@{^{(}->}[r] & \MTM(\Gr_G)^\comp,}
$$
where $e\co S \r \Gr_G$ is the inclusion of the base point.
It induces a commutative diagram of the corresponding Tannaka dual groups
\begin{equation}\label{full.tannaka.group}
\xymatrix{
\widetilde G_S \ar[r] \ar[d] & \scrG_S \ar[d] \\
\widehat G_1 \ar[r] & \bbG_{m,\bbQ}.
}
\end{equation}

\theo
\thlabel{full.Tannaka}
The diagram \eqref{full.tannaka.group} is Cartesian, i.e., it induces an isomorphism $\Q$-group schemes
$$\alpha\co \widetilde G_S \overset{\simeq}{\lr} \widehat G_1 \x_{\Gm} \scrG_S = \scrU_S \rtimes \widehat G_1.$$
In other words, there is an equivalence of Tannakian categories
$$
\MTM(L^+G\bsl LG/L^+G)^\comp=\MTM(\Gr_G)^\comp \;=\; \Rep_\bbQ^\fd(\scrU_S \rtimes \widehat G_1)
$$
between the category of compact mixed Tate motives on the double quotient $L^+ G \bsl LG / L^+ G$ over $S$ \textup{(}or equivalently, compact mixed Tate motives on $\Gr_G$, cf.~\thref{Tannaka_Cat}\textup{)} and the category of finite-dimensional representations of $\scrU_S \rtimes \widehat G_1$ \textup{(}regarded as a pro-algebraic group over $\Q$\textup{)}.
\xtheo

\pf
To check $\alpha$ is an isomorphism, we consider a $\Q$-algebra $R$ and evaluate the $R$-points of the two group schemes.
Recall that by definition that $\widetilde G(R)$ consists of the tensor automorphisms of the fiber functor $\omega\co \MTM(\Gr_G)^\comp \to \Vect_\Q$, i.e., families of additive $R$-linear automorphisms $g_X\co \omega(X) \t R \r \omega(X) \t R$ satisfying natural compatibility relations (see \cite[Ch.~II, §2]{DeligneMilneOgusShih:Hodge}), namely (for all  $X, X_1, X_2 \in \MTM(\Gr_G)^\comp$, $\tau\co X_1 \r X_2$)
\begin{enumerate}
\item
$g_1  =  \id_R$,
\item
$g_{X_1 \star X_2}  =  g_{X_1} \t g_{X_2}$, and
\item
$g_{X_2} \circ \tau  =  (\omega(\tau)\t1)  \circ g_{X_1}$.
\end{enumerate}

The map $\widetilde G(R) \r \widehat G_1(R)$ sends such a family $g := (g_X)_{X \in \MTM(\Gr)}$ to the collection of automorphisms, where $X$ only lies in the $\star$-subcategory $\Sat_G \subset \MTM(\Gr_G)$.
Likewise for $e_* \MTM(S) \subset \MTM(\Gr_G)$.
The injectivity of $\alpha$ therefore follows from \thref{parity.decomposition}: any motive in $\MTM(\Gr_G)$ is (isomorphic to one) of the form
\[
X = \bigoplus_{i=1}^n \IC_{L_i, \mu_i} = \bigoplus_{i=1}^n (\underbrace{e_* L_i}_{\in e_* \MTM(S)} \star \underbrace{\IC_{\mu_i}}_{\in \Sat_G}).
\]
Thus, $g_X = \bigoplus_i g_{e_* L_i} \t g_{\IC_{\mu_i}}$ is trivial if $g|_{\Sat_G} = \id$ and $g|_{e_* \MTM(S)} = \id$.

Given families of automorphisms $h := (h_X)_{X \in \Sat_G}$ and $h' := (h'_Y)_{Y \in e_* \MTM(S)}$ such that $h_X = h'_X$ for $X \in i_* \MTM(S)^{\pure}$, the surjectivity of $\alpha(R)$ requires to glue these families of automorphisms to one on $\MTM(\Gr_G)$.
To show this we use that there is an isomorphism
$$\IC_\mu(n) \star e_* L \simeq  \IC_{\mu'}(n') \star e_* L'$$
(if and) only if $\mu = \mu'$ and $L' \simeq L(n-n')$.
Indeed the former follows from support considerations, the latter follows by restricting the isomorphism of motives to $\Gr_G^{ \mu } \subset \Gr_G^{\le \mu}$.
For a motive $X$ as above, we can therefore define $\lambda_X := h_{\IC_\mu(n)} \t h'_{e_* L}$ independently of the presentation of $X$ as a convolution product.
We extend this additively.
For a morphism $\tau\co X := \bigoplus_i \IC_{\mu_i, L_i} \r \bigoplus_{i'} \IC_{\mu'_{i'}, L'_{i'}} =: X'$, we obtain $g_{X'} \circ \tau = (\omega(\tau) \t 1) \circ g_X$ by the description of the Hom group in \thref{Hom.IC} and the corresponding functoriality property for $h'$.
\xpf

\subsection{Extension to Ind-Categories}
The equivalences in \thref{Satake} and \thref{full.Tannaka} admit the following extensions to not necessarily compact objects. Such a statement can be useful in contexts when one wants to invoke adjoint functor theorems.

\coro
\thlabel{Satake.Ind.objects}
There are equivalences of symmetric monoidal abelian $\bbQ$-linear categories
$$\eqalign{
\Sat_G & = \Rep_\bbQ(\widehat G_1), \cr
\MTM(L^+G_S\bsl LG_S/L^+G_S)=\MTM(\Gr_{G, S}) & = \Rep_\bbQ(\scrU_S \rtimes \widehat G_1),}$$
where $\Rep_\bbQ$ denotes the category of not necessarily finite-dimensional representations.
\xcoro

\pf
Any representation of a flat group scheme, in particular any representation of a pro-algebraic group such as $\scrU_S \rtimes \widehat G_1$ is locally finite \cite[§I.2.13]{Jantzen:Representations}.
As finite-dimensional representations are compact objects, the category $\Rep_\bbQ(\scrU_S \rtimes \widehat G_1)$ is compactly generated by its subcategory $\Rep_\bbQ(\scrU_S \rtimes \widehat G_1)^\fd$.
On the other hand, $\MTM(\Gr_G)$ is compactly generated by $\MTM(\Gr_G)^\comp$ by virtue of \inter{Prop.~\iref{prop--MTM.G}}.
Similarly, $\Sat_G$ is compactly generated by definition.
We therefore obtain the claim by applying the ind-completion to the equivalence in \thref{full.Tannaka}.
\xpf

\subsection{From motives to functions}\label{functions.sec}
Let $S=\Spec(\bbF_q)$ be the spectrum of a finite field. 
We consider the (spherical) Hecke ring
\[
\calH_{G}\defined \calC_c\big(L^+G(\bbF_q)\bsl LG(\bbF_q)/L^+G(\bbF_q); \bbZ\big)
\]
where the ring structure is given by the convolution of functions, cf.~\cite[\S2]{Gross:Satake}. 
This ring has a $\bbZ$-basis given by the characteristic functions $c_\mu$, $\mu\in X_*(T)^+$ on the double cosets $L^+G(\bbF_q)\varpi^\mu L^+G(\bbF_q)$.
Taking the trace of geometric Frobenius on motives as in \cite{Cisinski:SurveyCoho} induces a surjective ring morphism 
\begin{equation}\label{trace.eq}
K_0\Sat_G^\comp\;\to\; \calH_G\otimes_\bbZ \bbZ[q^{-1}], \; M\mapsto f_M,
\end{equation}
where $K_0\Sat_G^c$ denotes the Grothendieck ring. 

Recall that the trace of the geometric Frobenius on $\Ql(1)$ is $q^{-1}$, cf.~\cite[(1.2.5) (iv)]{Deligne:Weil2}.
Thus, for the trivial group $G = 1$, the preceding map is the ring homomorphism $\Z[t, t^{-1}] \r \Z[q^{-1}]$ sending the class $t$ of $\Q(1)$ to $q^{-1}$.

In general, if $M$ is any stratified Tate motive on $\Gr_G$ (with respect to the stratification into $L^+G$-orbits), then its restriction $\iota_\mu^*M$ to some Schubert cell $\Gr_G^\mu$ lies in $\DTM(\Gr_G^\mu)$ by construction.
Thus, if $M$ is also compact, then $\iota_\mu^*M$ is a finite successive extension of Tate motives $1_\mu(n)[m]$, $n,m\in \bbZ$.
This immediately implies that the function $f_M$ takes values in $\bbZ[q^{-1}]$ as opposed to $\bbQ$.
Also since $\iota_\mu^*\IC_\mu$ coincides with $1_\mu$ up to a shift, the function $c_\mu$ appears with multiplicity $\pm 1$ in $f_{\IC_\mu}$ which shows that \eqref{trace.eq} is surjective. 

Let $R(\widehat G_1)$ be the Grothendieck ring of the category $\Rep_\bbQ^{\fd}(\widehat G_1)$.
By construction there is a morphism $d\co \widehat G_1\to \bbG_{m,\bbQ}, [(g,\la)]\mapsto \la^2$. 
When viewed as a representation this is nothing but the highest weight representation $V_0(1)$ of $\widehat G_1$.
We denote its class in $R(\widehat G_1)$ by $[V_0(1)]=[d]$.

\coro 
The motivic Satake equivalence in \thref{Satake} induces an isomorphism of rings
\[
\calH_G\otimes_\bbZ \bbZ[q^{-1}]\;\overset{\simeq}{\longrightarrow}\; R(\widehat G_1)/([d^{-1}]-q).
\]
\xcoro
\pf
The motivic Satake equivalence induces an isomorphism on Grothendieck rings 
\[
K_0\Sat_G^c\;\simeq\; K_0\Rep_\bbQ^{\fd}(\widehat G_1)\;=\;R(\widehat G_1)
\]
under which the class $[\IC_0(-1)]$ corresponds to $[V_0(-1)]=[d^{-1}]$.
Hence, it is enough to show that the kernel of \eqref{trace.eq} is the principal ideal generated by $a:=[\IC_0(-1)]-q[\IC_0]$.
As the geometric Frobenius acts on $\bbQ(-1)$ by multiplication with $q$, the class $a$ lies in the kernel.
Conversely, by induction one easily sees that all classes $[\IC_0(-n)]-q^{n}[\IC_0]$, $n\in \bbZ$ lie in the principal ideal generated by $a$.  
An elementary calculation using that the classes $[\IC_\mu]$, $\mu\in X_*(T)^+$ are linearly independent implies the corollary.
\xpf 

\rema
We note that the isomorphism of $\bbQ$-groups $\widehat G\x \bbG_{m,\bbQ}\to \widehat G\x \bbG_{m,\bbQ}$, $(g,\la)\mapsto (g,\la^{-1})$ induces an isomorphism $\widehat G_1\simeq \widehat G_1$, and hence an isomorphism of rings
\[
R(\widehat G_1)/([d^{-1}]-q)\;\simeq\; R(\widehat G_1)/([d]-q).
\]
Under the Tannakian dictionary this normalization corresponds to a sign change in the weight graduation in \thref{fiber.functor}. 
For a concrete example take $G=\on{PGL}_2$ so that $\widehat G=\on{SL}_2$. 
Since $(2\rho)(-1)$ is the diagonal matrix $\on{diag}(-1,-1)$, we get $\widehat G_1=\GL_2$ and $d=\on{det}\co \GL_2\to \bbG_m$.
Thus, the Satake isomorphism with the new normalisation reads in this case
\[
\calH_{\on{PGL}_2}\otimes_\bbZ \bbZ[q^{-1}]\;\overset{\simeq}{\longrightarrow}\; R(\on{GL}_2)/([\on{det}]-q).
\]
If $\mu\in X_*(T)^+=\bbZ_{\geq 0}$, then this isomorphism explicitly is given by 
\[
f_{\IC_\mu}\mapsto [\on{Sym}^\mu\bbQ^2],
\]
where $[\on{Sym}^\mu\bbQ^2]$ denotes the class of the $\mu$-th symmetric power of the standard representation of $\GL_2$.
\xrema

\appendix

\section{Complements on motives}\label{sect--box.product}

In this appendix, we show how the well-known compatibilities (stated in this form in \cite{JinYang:Kuenneth}) 
$$\eqalign{
(f_1 \x f_2)^* (M_1 \boxtimes M_2) & \cong (f_1^* M_1) \boxtimes (f_2^* M_2)\cr
(f_1 \x f_2)_! (M_1 \boxtimes M_2) & \cong (f_{1!} M_1) \boxtimes (f_{2!} M_2)}\eqlabel{boxbox}
$$
of the exterior product $\boxtimes$ of motives with $f^*$ and $f_!$ can be coherently organized.
This will be used to establish the structure of a symmetric lax monoidal 
functor for the presheaf $\DM^!$ when we restrict to placid morphisms between placid schemes, such as the $L^+ G$-torsors $LG^{\leq \mu} \r \Gr^{\leq \mu}$. This is used in the discussion of the convolution product in \refsect{convolution.product}.

Throughout \refsect{box.product}, we assume the base scheme $S$ to be a Noetherian separated scheme of finite Krull dimension. 
We note that all separatedness assumptions (here, on $S$, and below on certain maps of $S$-schemes) could eventually be dropped by proceeding as in \cite[Prop.~2.1.14]{RicharzScholbach:Intersection}.

\subsection{The exterior product}

We equip the categories $\Sch_S^{\ft}$, $\Sch_S$ and $\Cat$, $\Cat_\infty$, the (\ii-)category of all small \mbox{(\ii-)}\-categories, with their cartesian symmetric monoidal structure \cite[§2.4.1]{Lurie:HA}. 

In particular, we consider the cocartesian fibration $(\Sch_S^{\ft})^\x \r \Fin$ taking values in the category of finite pointed sets.
Recall that the objects of $(\Sch_S^{\ft})^\x$ are sequences $(X_1, \dots, X_n)$ with $X_i \in \Sch_S^\ft$ and, among others, the category has morphisms of the form $(X_1, \dots, X_n) \r X_1 \x_S \dots \x_S X_n$, corresponding to $\id_{X_1 \x \dots \x X_n}$.
Let $(\Sch_S^\ft)^{\x, \dual} \r \Fin^\opp$ be the associated cartesian fibration as constructed in \cite{BarwickGlasmanNardin:Dualizing}. 
The \emph{opposite} of this, which is again a cocartesian fibration, encodes the usual symmetric monoidal structure on $(\Sch_S^\ft)^\opp$. 
We will abbreviate the source of this map as $(\Sch_S^\ft)^{\opp, \x}$ or even just $(\Sch_S^\ft)^\opp$.

The subcategory $\AffSch^\ft_S \subset \Sch^\ft_S$ (consisting of affine finite type $S$-schemes, $\Spec R \r S$) is closed under the product since $S$ is by assumption separated.
We further endow $\DGCat_\cont$ and $\PrL$ (presentable \ii-categories with colimit-preserving functors) with the Lurie tensor product, see e.g. \cite[Ch.~1, §6]{GaitsgoryRozenblyum:StudyI}.

\lemm
\thlabel{DM*.slm}
The functor
$$\eqalign{\DM^* & : (\Sch^\ft_S)^\opp \r \DGCat_\cont }$$
admits a natural symmetric lax monoidal structure such that for finite type $S$-schemes $X_1, X_2$, this structure is the exterior product
$$\boxtimes\co \DM(X_1) \t \DM(X_2) \r \DM(X_1 \x X_2).$$
\xlemm

\pf
The functor $\Sch_S^{\ft, \opp} \r \Cat$, $X \mapsto \Sm / X$ is symmetric lax monoidal (with respect to the cartesian monoidal structures on both categories) by means of the exterior product.
The inclusion $\Cat \r \Cat_\infty$ is symmetric monoidal.
The presheaf functor (in the \ii-categorical sense) $\calP\co \Cat_\infty \r \PrL$ is symmetric monoidal \cite[Rem.~4.8.1.8]{Lurie:HA}.
Thus, the composite $X \mapsto \calP(\Sm / X)$ is symmetric lax monoidal. In addition, the (non-full) subcategory $W_X \subset \calP(\Sm / X)$ consisting of the usual $\A^1$-projections and étale hypercoverings are monoidal subcategories, so that the functor $X \mapsto (\calP(\Sm / X), W_X)$ is a symmetric monoidal functor taking values in the \ii-category $\WCat_\infty$ of relative \ii-categories.
The localization functor $\WCat_\infty \r \Cat_\infty$, $(C, W) \mapsto C[W^{-1}]$ is symmetric lax monoidal \cite[Prop.~4.1.7.2, Prop.~4.1.7.4]{Lurie:HA}.

The stabilization process, i.e., turning $\P^1$ into an invertible object, is also a symmetric lax monoidal functor.
This is readily apparent from the description of this process in \cite[§4.1]{Robalo:Theorie}:
abbreviating the notation of loc.~cit. as
$P := \calP(\free^\t(\Delta[0]))^\t$ and
$P_\inv := \calP(\calL^\t_{\free^\t(\Delta[0]),*)} (\free^\t(\Delta[0])))^\t$, let $\CAlg(\PrLt)_\pt$ be the undercategory $\CAlg(\PrLt)_{P/}$.
Its objects are pairs $(C, X)$ consisting of a presentable symmetric monoidal \ii-category $C$ and an object $X \in C$.
Similarly, consider the undercategory $\CAlg(\PrLt)_{\pt, \inv} := \CAlg(\PrLt)_{P_\inv/}$ whose objects consist of similar pairs $(C, X)$, but where $X$ is a \emph{$\t$-invertible} object.
The objects $P$ and $P_\inv$ have natural comonoid structures stemming from the comonoid structure present on any object in a cartesian symmetric monoidal category such as $\Cat_\infty^\x$.
Thus, the undercategories under these two objects have a natural symmetric monoidal structure in such a way that the functor $(C, X) \mapsto X$ is symmetric monoidal.
The natural functor, arising from the map $P \r P_\inv$,
$$\CAlg(\PrLt)_{P_\inv/} \r \CAlg(\PrLt)_{P/}$$
is symmetric monoidal. Hence its left adjoint, which by \cite[Def.~4.1.8]{Robalo:Theorie} is the functor mapping a pointed category $(C, X)$ to $(C[X^{-1}], X)$, is symmetric lax monoidal.
This abstract observation is applied to the functor $\Sch_S^{\ft, \opp} \r \CAlg(\PrLt)_{P/}$, $X \mapsto (\calP(\Sm/X)[\langle \A^1, \et \rangle^{-1}], \P^1_X)$ which is symmetric lax monoidal by the above (and $\P^1_X \x \P^1_Y = \P^1_{X \x Y}$).

The composite, denoted by $\SH^*\co \Sch_S^{\ft, \opp} \r \PrLt$ takes values in $\PrLt_\stb$, the \ii-category of stable presentable symmetric monoidal \ii-categories and colimit-preserving functors and is by the above a symmetric lax monoidal functor.
Finally, $\DM$ arises from composing with the symmetric monoidal functor $\PrL_\stb = \Mod_{\Sp}(\PrL) \stackrel {- \t \Q} \lr \Mod_{\Q}(\PrL_\stb) =: \DGCat_\cont$.
\xpf

We now retrace the construction of $\DM^*_!$ as a functor out of the category of correspondences, by keeping track of the symmetric lax monoidal structure and thus of projection formulas, including all their higher coherences.
We use, in the same vein as Hoyois \cite{Hoyois:Six} and Khan \cite{Khan:Motivic}, the universal property of the category of correspondences. 
An alternative approach for coherently encoding projection formulas avoiding the category of correspondences appears in \cite[§4.5]{AyoubGallauerVezzani:Six}.

Recall the $(\infty, 2)$-category of correspondences 
$\Corr((\Sch_S^\ft)^{\x, \dual})_{\horiz, \vert}^\adm$ defined in \cite[§7]{GaitsgoryRozenblyum:StudyI}.
Here $\horiz$, $\vert$ and $\adm$ are certain subcategories of $(\Sch_S^\ft)^{\x, \dual}$, to be specified below more concretely.
The objects of this category are the objects $X \in (\Sch_S^\ft)^{\x, \dual}$ (which are, in their turn, finite collections of objects in $\Sch_S^\ft$); 1-morphisms from $X$ to $Y$ are spans of the form $Y \stackrel g \gets Z \stackrel f \r X$ with $g \in \vert$ and $f \in \horiz$, and 2-morphisms between such a morphism and another similar correspondence is a map $Z \r Z'$ in $\adm$ that is compatible with the maps to $X$ and $Y$.
This describes the low-dimensional data of this category, we refer to loc.~cit.~for the full definition including the $(\infty, 2)$-categorical structure.

\defi
A map $f$ in a symmetric monoidal \ii-category $\calC^\t \r \Fin$ is called \emph{dormant} if its image in $\Fin$ is an identity map.


We consider the category $\Corr((\Sch_S^\ft)^{\x, \dual})_{\open, \all}^\iso$, where ``$\iso$'' indicates the collection of \emph{dormant} isomorphisms, ``$\open$'' is the subcategory consisting of \emph{dormant} morphisms 
$$(f_1, \dots, f_n)\co (X_1, \dots, X_n) \r (Y_1, \dots, Y_n)$$
where each $f_i$ is an open embedding.
Finally, ``$\all$'' means that no restriction is imposed on the horizontal morphisms.
\xdefi

The extension of the functor $(\Sch_S^\ft)^{\x, \dual} \r \Fin^\opp$ to correspondences, and restricting to the above subcategory gives a functor $\Corr((\Sch_S^\ft)^{\x, \dual})_{\open, \all}^\iso \r \Corr(\Fin^\opp)_{\id, \all}^\iso = \Fin$.
Given a map $\alpha \co \langle m \rangle \r \langle n \rangle$ in $\Fin$ and a cocartesian lift $(X_1, \dots, X_n) \r (Y_1, \dots, Y_m)$ in $(\Sch_S^\ft)^\opp$, on checks that $(X_i) \stackrel \id \gets (X_i) \r (Y_j)$ is a cocartesian lift in $\Corr((\Sch_S^\ft)^{\x, \dual})$.
The Segal condition holds trivially by definition, so that $\Corr((\Sch_S^\ft)^{\x, \dual})_{\open, \all}^\iso$ is a symmetric monoidal \ii-category.

\lemm
The functor $\DM^*$ extends uniquely to a symmetric lax monoidal functor
$$\DM^*_\sharp : \Corr((\Sch_S^\ft)^{\x, \dual})_{\open, \all}^\iso \r \DGCat_\cont.\eqlabel{DM*sharp}$$
\xlemm

\pf
We show the existence of the functor using \cite[Chapter 7, Theorem 3.2.2(b)]{GaitsgoryRozenblyum:StudyI}.
Among the general conditions in Chapter 7, 1.1.1 there, the only requirement to check is the existence of pullbacks in $(\Sch_S^\ft)^{\x, \dual}$ along maps in the subcategory ``$\open$''.
This is clear: given open embeddings $f_i$ and a map as in the horizontal arrow, i.e., a map $\alpha\co \langle n \rangle \r \langle m \rangle$, and $g_i\co Z_i \r \prod_{j \in \alpha^{-1}(i)} Y_j$, the diagram
$$\xymatrix{
(Z_i \x_{\prod_{j \mapsto i} Y_j} \prod_{j \mapsto i} X_i) \ar[r] \ar[d]
&
(X_1, \dots, X_n) \ar[d]^{(f_1, \dots, f_n)} \ar[d] \\
(Z_1, \dots, Z_m) \ar[r]^{(g_1, \dots, g_m)} & (Y_1, \dots, Y_n).
}$$
is cartesian in $(\Sch_S^\ft)^{\x, \dual}$ and the left vertical map is again in ``$\open$''.
In order to check the Beck--Chevalley condition, it suffices to separately consider the case where $(g_i)$ is inert, respectively active, since these form a factorization system.
For inert morphisms, this is clear.
For active morphisms, we may assume that $m=1$ and $n=2$ above, in which case we consider $g\co Z_1 \r Y_1 \x Y_2$.
Then the Beck--Chevalley condition is the assertion that the following diagram commutes, which follows from the construction of $\boxtimes$:
$$\xymatrix{
\DM(Z_1 \x_{Y_1 \x Y_2} X_1 \x X_2) \ar[d] & 
\DM(X_1 \x X_2) \ar[l] \ar[d]^{(f_1 \x f_2)_\sharp} &
\DM(X_1) \t \DM(X_2) \ar[d]^{(f_1)_\sharp \x (f_2)_\sharp} \ar[l] \\
\DM(Z_1) & 
\DM(Y_1 \x Y_2) \ar[l]^{g^*} &
\DM(Y_1) \t \DM(Y_2) \ar[l]^{\boxtimes}.
}\eqlabel{box.sharp.*}$$
The functor $\DM^*_\sharp$ obtained in this way clearly preserves edges that are cocartesian over $\Fin$, thus giving a symmetric lax monoidal functor.
\xpf

\lemm
Let $\sep$, resp.~$\proper$ be the subcategory of $(\Sch_S^\ft)^\x$ spanned by morphisms that are dormant (i.e., map to an identity in $\Fin$), and are componentwise separated (resp.~proper).
The functor $\DM^*_\sharp$ in \refeq{DM*sharp} extends uniquely to a symmetric lax monoidal functor:
$$\DM^*_! \co \Corr((\Sch_S^\ft)^{\x, \dual})_{\sep, \all}^\proper \r \DGCat_\cont.$$
\xlemm

\pf
We apply \cite[Chapter 7, Thm.~5.2.7]{GaitsgoryRozenblyum:StudyI}.
To check its assumptions note that a map in $\proper \cap \open$ is dormant, and componentwise an open embedding. 
Such a map is a monomorphism in $(\Sch_S^\ft)^{\x, \dual}$.
The condition in Chapter 7, 5.2.2 there and also the Beck--Chevalley condition is satisfied since, again, the corresponding properties only need to be checked for maps in $\open$, resp.~$\proper$, which are by definition \emph{dormant}.
Thus the Beck--Chevalley condition reduces to an assertion similar to the commutativity of \refeq{box.sharp.*}, except that $(-)_\sharp$ is replaced by $(-)_!$, in other words, the classical projection formula as recalled in \refeq{boxbox}.
\xpf

\coro
\thlabel{DM*vs!}
Write $\Sm^{\ft}_{S, \sm \cap \sep}$ for the category consisting of smooth separated finite-type $S$-schemes and smooth separated morphisms.
There is a natural isomorphism of functors
$$\Tw : \DM^*|_{\Sm^{\ft}_{S, \sm \cap \sep}} \Rightarrow \DM^!|_{\Sm^{\ft}_{S,  \sm \cap \sep}} : (\Sm^{\ft}_{S,  \sm \cap \sep})^\opp \r \DGCat_\cont,$$
whose evaluation at a map $f: X \r Y$ is the natural transformation
$$\xymatrix{\DM(Y) \ar[rr]^{- \t \omega_Y } \ar[d]^{f^*} & & \DM(Y) \ar[d]^{f^!} \\
\DM(X) \ar[rr]_{- \t \omega_X } \ar@{=>}[urr] & & \DM(X)}$$
stemming from the projection formula.
Here, $\omega_X := p_X^! 1$ with $p_X\co X \r S$ the structural map.
\xcoro

\pf
Any map $X \stackrel f \r Y$ of schemes naturally gives rise to a pair $(Y, X)$, where $Y$ is a comonoid and $X$ is a $Y$-comodule (both with respect to the cartesian monoidal structure) on $\Sch$.
The coaction is given by $X \stackrel{\Delta_X} \r X \x X \stackrel{f\x \id} \r Y \x X$.
Applying the symmetric lax monoidal functor $\DM^*_!$ to this object, we obtain an object in $\calL \calM(\DGCat_\cont)$, namely the commutative algebra object $\DM(Y)$, and the $\DM(Y)$-module $\DM(X)$, where the action arises as
$$\DM(Y) \t \DM(X) \stackrel{\boxtimes} \r \DM(Y \x X) \stackrel{(f \x \id)^*} \r \DM(X \x X) \stackrel{\Delta_X^*} \r \DM(X).$$
This can also be computed as the natural action of $\DM(Y)$, via $f^*$, on the symmetric monoidal category $\DM(X)$ (equipped with its usual $\t$).

The map $f$ also gives rise to a map of comodule objects $(Y, X) \r (Y, Y)$.
Applying the functor $\DM^*_!$ to the map of induced left module objects, namely the pair
$$(\id, f_!)\co (\DM(Y), \DM(X)) \r (\DM(Y), \DM(Y)).$$
(Such an interpretation of the projection formula was observed by Khan \cite{Khan:Motivic}.
Note, however, that the approach to projection formulas laid out in op.~cit.~does not seem to work as is, since the morphisms in Ch.~2, §4.1.6 there cannot be composed.)
This map admits a left adjoint separately for each object in $\calL \calM$, namely $\id$ and $f^!$, respectively.
\nts{More formally, it takes values in the full subcategory $\Alg_{\calL \calM^\t}(\DGCat_\cont^R)$, where the superscript $R$ denotes the \ii-category consisting of the same objects as before, but only functors admitting a left adjoint are allowed.}
By \cite[Cor.~7.3.2.7]{Lurie:HA}, the functor therefore admits a right adjoint relative to $\calL \calM^\t$, still denoted $f^!$.
\nts{More formally, we use the equivalence given by Lurie $\Alg_{\calL \calM^\t}(\DGCat_\cont^R)=(\Alg_{\calL \calM^\t}(\DGCat_\cont))^R$.}

This in particular expresses the existence of a natural map $f^* A \t f^! B \r f^! (A \t B)$ that is functorial in $A$ and $B$.
It follows from the naturality of the construction that it is also functorial in $f$.
The sought-for transformation is defined as the restriction of this map to $B=\omega_Y$.
By relative purity, this map is an isomorphism whenever $f$ is smooth.
Moreover, for $X$ and $Y$ smooth, $\omega_X$ and $\omega_Y$ are $\t$-invertible.
\xpf

\lemm
\thlabel{DM!.slm}
Consider the functor
$$\DM^! : (\Sch_{S, \sep}^\ft)^{\x, \dual} \r \DGCat_\cont$$
obtained from the composite $(\Sch_{S, \sep}^\ft)^{\x, \dual} \r \Corr((\Sch_S^\ft)^{\x, \dual})_{\sep, \all}^\proper \stackrel {\DM^*_!} \r \DGCat_\cont$ by passing to right adjoints.
\begin{itemize}
	\item This functor is symmetric lax monoidal if $S = \Spec k$ is a field.
	\item The restriction of this functor to $(\Sch_{S, \sm \cap \sep}^\ft)^{\x, \dual}$ is symmetric lax monoidal for general $S$.
	Here the subscript $\sm \cap \sep$ refers to the (non-full, symmetric monoidal) subcategory comprising all finite type $S$-schemes, but only smooth separated maps.
	\nts{More formally, we take $\Sch^\x \x_{\Sch} \Sch_{\sm \cap \sep}$. Note that $\Sch_{\sm \cap \sep}$ is a monoidal (non-full) subcatgory of $\Sch$, so that this becomes a symmetric monoidal category.}
\end{itemize}
The fiber over $\langle 1 \rangle$ of $(\Sch_{S, \sm \cap \sep}^\ft)^{\x, \dual}$ identifies with $(\Sch_{S, \sm \cap \sep}^\ft)^\opp$, and we also denote this symmetric lax monoidal functor by $\DM^!\co (\Sch_{S, \sm \cap \sep}^\ft)^\opp \r \DGCat_\cont$.
\xlemm

\pf
The functor $\DM^!$ exists as stated, since the right adjoints happen to preserve colimits as well. (This is well-known and uses the assumptions that $S$ is Noetherian and of finite Krull dimension.)
To check it is a symmetric lax monoidal functor it remains to check that for two maps $f_1$, $f_2$ in $\Sch_S^\ft$ the natural map $(f_1)^! M_1 \boxtimes (f_2)^! M_2 \r (f_1 \x f_2)^! (M_1 \boxtimes M_2)$ is an isomorphism.
If $S$ is a field, this holds by \cite[Prop.~2.3.5]{JinYang:Kuenneth} (note this is nontrivial and uses alterations).
For general $S$, but smooth maps $f_i$, this holds by relative purity and \refeq{boxbox} for $*$-pullbacks.
\xpf

\subsection{Motives on placid prestacks}

For a regular cardinal $\kappa$, recall from \thref{prestack.DM} the category $\AffSch_S^\kappa$. It consists of those affine schemes that can be presented as $\kappa$-small cofiltered limits
$$X = \lim X_i,\eqlabel{X.lim}$$
where the $X_i$ are affine finite type $S$-schemes.

\defi[{\cite[App.~C]{Gaitsgory:Local}, \cite[Def.~4.2.1]{Raskin:D-modules}}]
\thlabel{placid}
An object in $\AffSch_S^\kappa$ is called \emph{placid} if it admits a \emph{placid presentation}, i.e., one such that the transition maps $X_i \r X_j$ in \refeq{X.lim} are smooth (and necessarily affine).
A map between two such placid affine $S$-schemes is called \emph{placid} if for any pair of placid presentations $X = \lim X_i$, $Y = \lim Y_j$, and any $j$, there is some $i$ such that $X \r Y \r Y_j$ factors as $X \r X_i \r Y_j$, where the second map is smooth.
(It follows from \cite[Lemma~4.5.1]{Raskin:D-modules} that this condition only needs to be checked for any fixed presentations of $X$ and $Y$.)
The non-full subcategory of $\AffSch^\kappa_S$ consisting of placid affine schemes and placid maps is denoted by $\AffSch_S^{\kappa,\pl}$.
(Equivalently, \cite[Rem.~4.10.3]{Raskin:D-modules}, $\AffSch_S^{\kappa,\pl}$ can be defined as the pro-category $\Pro_{\kappa\text{-small}}(\AffSch^\ft_{S,\sm})$, the pro-completion of affine schemes of finite type over $S$ with smooth maps.)

We call the category $\PreStk^\pl_S := \Fun((\AffSch_S^{\kappa,\pl})^\opp, \text{\ii-}\Gpd)$ the category of \emph{placid prestacks}.
It is the free completion of $\AffSch_S^{\kappa, \pl}$ under arbitrary colimits.

The restriction of prestacks to $\AffSch_S^{\kappa,\pl}$ induces an adjunction $\PreStk^\pl_S\leftrightarrows \PreStk_S$.
A prestack is called {\it placid} if it lies in the essential image of the functor $\PreStk^\pl_S\r \PreStk_S$.
\xdefi

\exam
The 
groups $\calP_\bbf$ (in particular the positive loop group $L^+G$) are placid affine $S$-schemes, provided that $S$ itself is affine.
So the preimage of any open affine subscheme under the $\calP_\bbf$-torsor $LG^{\le w}\to \Fl^{\le w}$ (see \refsect{DTM.Tate}) is placid affine.

Thus, for any pair of facets $\bbf,\bbf'$, the double quotient of the $\calP_{\bbf'} \x \calP_{\bbf}$-action (from the right and the left) on $LG^{\le w}$ exists as a placid prestack $\calP_{\bbf'} \setminus LG^{\le w} / \calP_{\bbf}$.
More generally, for any quasi-compact closed subscheme $X\subset LG$ the prestack $\calP_{\bbf'} \setminus X / \calP_{\bbf}$ is placid.
\xexam

Our goal is to have $\boxtimes$-products for motives on placid prestacks.
By means of the following two results, this is a formal consequence of the symmetric lax monoidality of $\DM^!$ on $\AffSch_{S, \sm}^\ft \subset \Sch_{S, \sm \cap \sep}^\ft$.

\lemm
\thlabel{Day.convolution}
Let $\calC$ be a small symmetric monoidal \ii-category and $\calD$ be a cocomplete symmetric monoidal \ii-category whose tensor product preserves colimits separately in each variable.
Fix a regular cardinal $\kappa$.
Then the restriction functor $\Fun(\Ind_{\kappa\text{-small}}(\calC), \calD)^\t \r \Fun(\calC, \calD)^\t$
admits a symmetric monoidal left adjoint, where the monoidal structure on the functor categories is given by Day convolution.

Thus, any symmetric lax monoidal functor $F\co \calC \r \calD$ can be left Kan extended to a symmetric lax monoidal functor on the ind-completion $\Ind_{\kappa\text{-small}}(\calC)$.
\nts{The same proof also works for the presheaf category instead of the ind-completion, with \cite[Thm. 5.1.5.6]{Lurie:Higher} instead.
However this is not something we can apply to extend the symmetric lax monoidal functoriality of $\DM^*$ to prestacks since the Lurie tensor product does not commute with limits.}
\xlemm

\pf
This can be proven as in \cite[Cor.~3.8]{Nikolaus:Stable}.
The assumption in loc.~cit. that $\calD$ is accessible is, for the particular situation considered here, not needed since the invokation of the adjoint functor theorem in the proof of loc.~cit. can be replaced by the universal property of the ind-completion \cite[Prop.~5.3.5.10]{Lurie:Higher}.
The last statement follows from the equivalence of symmetric lax monoidal functors and commutative monoid objects in the functor category under the Day convolution \cite[Prop.~2.12]{Glasman:Day}.
\xpf

In the next statement, $\PreStk^{(\pl)}_S$ is equipped with the cartesian symmetric monoidal structure.
\nts{We can alternatively also consider the Day convolution by composing this with the opposite of the symmetric oplax monoidal functor $\PreStk_S^\t \stackrel \id \r \PreStk_S^\x$. The symmetric oplax monoidal structure maps are given, for two prestacks $F$ and $G$ by
$$F \t G = \colim_{X \mapsto F, Y \mapsto G} X \x Y \r F \x G,$$
where the colimit runs over the $X \in \AffSch$ mapping to $F$ etc.}

\coro
\thlabel{DM!*.PreStk.monoidal}
The functor $\DM^*$ on $(\AffSch_S^\ft)^\opp$ admits a natural symmetric lax monoidal extension to $\AffSch_S^\kappa$ and to $\PreStk_S$.
The same is true for $\DM^!$ provided that $S = \Spec k$ is a field.
Finally, the functor $\DM^!$ admits a symmetric lax monoidal extension to $\AffSch_S^{\kappa, \pl}$ and to $\PreStk_S^\pl$.
\xcoro

\pf
First, apply \thref{Day.convolution} to 
$\C=(\AffSch_S^{\ft})^\opp$. 
Second, in order to extend $\DM^*$ from $\AffSch_S$ to a symmetric lax monoidal functor on $\PreStk_S^\x$ we use the argument in \cite[Chapter 9, Prop.~3.2.4]{GaitsgoryRozenblyum:StudyI}, according to which it is enough to observe that for any prestacks $F_1, \dots, F_n$, the map
$$\lim_{X \in \AffSch_S, X \r \prod F_i} \DM(X) \r \lim_{X_i \in \AffSch_S, X_i \r F_i} \DM(\prod X_i)$$
is an equivalence for cofinality reasons.
This shows the claim for $\DM^*$.

The one for $\DM^!$ (and placid prestacks, or arbitrary ones for $S$ being a field) is done the same way, using \thref{DM!.slm} instead.
\xpf

\coro
\thlabel{DM!.placid.prestacks}
The restriction of $\DM^!$ to $(\PreStk_S^\pl)^\opp$ is symmetric lax monoidal.
In particular, for any two placid prestacks $X_1$, $X_2$, there is a natural functor
$$\DM(X_1) \t \DM(X_2) \r \DM(X_1 \x_S X_2).$$
\xcoro

\pf
The restriction of $\DM^!$ to $\PreStk_S^\pl$ is the unique colimit-preserving functor extending the restriction of $\DM^!$ to $(\AffSch_S^\pl)^\opp$.
This latter functor is the unique extension, preserving ($\kappa$-small) cofiltered limits, of the restriction of $\DM^!$ to $(\AffSch_{S, \sm}^{\ft})^\opp$.
We can conclude using \thref{DM!*.PreStk.monoidal}.
\xpf

\ifusebiber
  \printbibliography[title={References}]
\else
  \bibliographystyle{alphaurl}
  \bibliography{bib}
\fi

\end{document}